\documentclass[12pt]{article}
\usepackage{color}
\usepackage[colorlinks,linkcolor=blue,citecolor=blue,anchorcolor=red]{hyperref}
\usepackage{graphicx,amsmath,amsthm,amssymb,dsfont, mathrsfs,MnSymbol, fontenc, bbm}

\def\RR{\mathbb R}

\def\ZZ{\mathbb Z}

\def\QQ{\mathbb Q}

\def\DD{\mathbb D}

\def\AA{\mathbb A}

\def\XX{\mathbb X}
\def\YY{\mathbb Y}

\def\UU{\mathbb U}
\def\SS{\mathbb S}
\def\TT{\mathbb T}

\def\cW{\mathcal W}

\def\fr{\mathfrak R}
\def\rD{\boldsymbol{\rm D}}

\def\cH{\mathcal H}
\def\cT{\mathcal T}
\def\sL{\mathscr L}

\def\cJ{\mathcal J}
\def\fj{\mathcal j}

\def\fF{\mathfrak F}

\def\ff{\mathfrak f}

\def\fj{\mathfrak j}
\def\fa{\mathfrak a}
\def\fb{\mathfrak b}
\def\fh{\mathfrak h}

\def\fm{\mathfrak m}

\def\mm{\mathbbm m}
\def\pp{\mathbbm p}

\def\fd{\mathfrak d}
\def\fr{\mathfrak r}

\def\fO{\mathcal O}
\def\cP{\mathcal P}

\def\cD{\mathcal D}
\def\fD{\mathfrak D}
\def\sD{\mathscr D}
\def\cN{\mathcal N}

\def\cT{\mathcal T}
\def\fT{\mathfrak T}

\def\fU{\mathfrak U}

\def\fs{\mathfrak s}
\def\cS{\mathcal S}

\def\cX{\mathcal X}
\def\cW{\mathcal W}

\def\bx{{\bf x}}
\def\by{{\bf y}}

\def\bn{{\bf n}}
\def\bm{{\bf m}}

\def\bbp{{\mathbbm p}}

\def\bb{{\bf b}}

\def\br{{\bf r}}

\def\sh{\hslash}
\def\balpha{\boldsymbol{\alpha}}

\def\brho{\boldsymbol{\rho}}
\def\bmu{\boldsymbol{\mu}}

\def\bfr{\mathfrak r}
\def\bfc{\boldsymbol{\mathfrak c}}

\def\b1{{\bf 1}}
\def\d1{{\mathbbm{1}}}
\def\one{{\mathbbm{1}}}
\def\mod{{\; \rm mod \;}}
\def\ad{{\rm and}}

\def\with{{\rm with}}
\def\where{{\rm where}}
\def\for{{\rm for}}

\def\card{{\rm card}}
\def\ord{{\rm ord}}

\def\qed{ \ \vrule width.2cm height.2cm depth0cm\smallskip}
\begin{document}
\baselineskip=15pt


\title{On the upper bound of the $L_p$ discrepancy  of Halton's sequence and the
 Central Limit Theorem for Hammersley's net}
\author{Mordechay B. Levin}

\date{}

\maketitle
\begin{abstract}
Let $ (H_s(n))_{n \geq 1} $ be an $s-$dimensional Halton's sequence, and let
${\cH}_{s+1,N}=(H_s(n),n/N)_{n=0}^{N-1}$ be the $s+1-$dimensional Hammersley's point set. Let $  D(\bx,(H_n)_{n=0}^{N-1}  )$ be the local discrepancy of $(H_n)_{n=0}^{N-1}$, and let
 $D_{s,p} ( (H_n)_{n=0}^{N-1}) $ be the $L_p$ discrepancy of  $ (H_n)_{n=0}^{N-1} $.
It is known that $\limsup_{N \to  \infty }
 (\log N)^{-s/2}  D_{s,p}  (  H_{s}(N) )_{n=0}^{N-1} >0$.
 In this paper, we prove that
 $$
           D_{s,p} ((  H_{s}(N) )_{n=0}^{N-1})  =O( \log^{s/2} N) \quad {\rm for} \; \; N \to  \infty .
 $$
I.e., we found the smallest possible order of magnitude of $L_p$ discrepancy  of Halton's sequence.
Then we prove the Central Limit Theorem for Hammersley's net :
\begin{equation}\nonumber 
  D(\bar{\bx},\cH_{s+1,N} ) /D_{s+1,2}(\cH_{s+1,N})
	\stackrel{w}{\rightarrow} \cN(0,1),
\end{equation}
where $\bar{\bx}$ is a uniformly distributed random variable in $[0,1]^{s+1}$.
The main tool is  the theorem on $\bbp$-adic logarithmic forms.
\end{abstract}
%
\newpage

\tableofcontents
\addtocontents{toc}{\protect\setcounter{tocdepth}{2}}

\section{\textcolor{blue}{Introduction}}
\subsection{\textcolor{blue}{Notations and definitions}}
 Let $\cP_{N} = (\beta_{n,N})_{n = 0}^{N-1}$ be an $N$-element point set in the
 $s$-dimensional unit cube $[0,1)^s$. The local {\bf discrepancy} function of $\cP_{N}$ is defined as
\begin{equation}\label{In1}
   D(\bx, \cP_{N}  )= \#\{0 \leq n \leq N-1 \; : \; \beta_{n,N} \in [0,x_1)
    \times \cdots \times [0,x_s) \} -  N x_1 \cdots x_s.
\end{equation}
We define the $\emph{L}_\infty$ and $\emph{L}_p$  discrepancy of an
$N$-point set $\cP_{N}$ as
\begin{equation} \label{In2}
       D_{s,\infty} (\cP_{N}) =
    \sup_{ 0<x_1, \ldots , x_s \leq 1} \; |
    D(\bx,\cP_{N}) |,  \qquad  D_{s,p}(\cP_{N})=\left\| D(\bx,\cP_{N}) \right\|_{s,p}  ,
\end{equation}
\begin{equation} \label{In3}
     \left\| f(\bx) \right\|_{s,p}= \Big( E_s( |f(\bx)|^p)\Big)^{1/p}, \qquad  E_s(f(\bx) )=  \int_{[0,1)^s}  f(\bx) d \bx .
\end{equation}

 {\it Definition 1.} {\it  A sequence $(\beta_n)_{n\geq 0}$ is of {\it
low discrepancy} (abbreviated l.d.s.) if  $D_{s,\infty}
((\beta_n)_{n=0}^{N-1})=O((\log N)^s) $ for $ N \rightarrow \infty $.\\
A sequence  of point sets $(\cP_{N})_{N=1}^{\infty}$ is of
 low discrepancy (abbreviated
l.d.p.s.) if $ D_{s,\infty} (\cP_{N})=O((\log
N)^{s-1}) $ for $ N \rightarrow \infty $.} \\
For examples of such a sequence, see, e.g., \cite{BeCh,Ni}.

In 1954, Roth proved that there exists a constant $ C_s>0 $, such
that
\begin{equation}  \label{In4}
\emph{D}_{s,\infty}((\beta_{n,N}^{(s)})_{n=0}^{N-1})>C_s(\log N)^{\frac{s-1}{2}} \;\;\;\;
    {\rm and} \;\;\;\;  \underline{\lim }{\emph{D}_{s,\infty}((\beta_n^{(s)})_{n=0}^{N-1})(\log
N)^{-s/2}}>0
\end{equation}
for all $N$-point sets $(\beta_{n,N}^{(s)})_{n=0}^{N-1}$ and all sequences
$(\beta_n^{(s)})_{n \geq 0}$.

According to the well-known conjecture (see, e.g., \cite[p.~283]{BeCh},
\cite[p.32]{Ni}), these estimates can be improved to
\begin{equation}   \label{In08}
  {\emph{D}_{s,\infty}((\beta_{n,N}^{(\ddot{s})})_{n=0}^{N-1}) (\log N)^{-\ddot{s}+1}} >C_{\ddot{s}}^{'}
     \;\; {\rm and} \;\;
\underset{N \to \infty }{\underline{\lim }}  (\log N)^{-\dot{s}} \emph{D}_{s,\infty}((\beta_{n}^{(\dot{s})})_{n=1}^{N})>0
\end{equation}
for all $N$-point sets $(\beta_{n,N}^{(\ddot{s})})_{n=0}^{N-1}$ and all sequences $(\beta_n^{(\dot{s})})_{n \geq 0}$ with some $C_{\ddot{s}}^{'} >0$.

In 1972, W. Schmidt proved (\ref{In08}) for $\dot{s}=1 $ and $\ddot{s}=2$.
In 1989,   Beck  proved that  $\emph{D}_{3,\infty}(\cP_{N}) \geq \dot{c} \log N  (\log\log N)^{1/8-\epsilon}$  for $s=3$ and some $\dot{c}>0$. In 2008, Bilyk,  Lacey
and Vagharshakyan (see \cite[p.76]{Bi}) proved   in all
 dimensions $s \geq 3$ that there exists some $\dot{c}(s), \eta >0$ for which the following
estimate holds for all  $N$-point sets :
$\emph{D}_{s,\infty} (\cP_{N})>\dot{c}(s)(\log N)^{\frac{s-1}{2}  +\eta}$.
In \cite{Le1}--\cite{Le3}, Levin proved that \eqref{In08} is true for Hammersley's net, known constructions of $(t,m,s)$ nets  and for  nets obtained from a
 module in a totally real algebraic number field.
It is known that
\begin{equation} \label{In7}
\emph{D}_{s,p}((\beta_{n,N}^{(s)})_{n=0}^{N-1})>C_{s,p}(\log N)^{\frac{s-1}{2}} \;\;\;
    {\rm and} \;\;\;  \underline{\lim }{\emph{D}_{s,p}((\beta_n^{(s)})_{n=0}^{N-1})(\log
N)^{-s/2}}>0
\end{equation}
for all $N$-point sets $(\beta_{n,N}^{(s)})_{n=0}^{N-1}$ and all sequences
$(\beta_n^{(s)})_{n \geq 0}$ with some $C_{s,p} >0$ (see  Roth for $p=2$, Schmidt  for $p>1$ \cite{BeCh}, and Proinov \cite{Pr}).\\

{\it Definition 2.}  A sequence $(\beta_n)_{n\geq 0}$ is of  $L_p$ {\it
low discrepancy} (abbreviated l.d.s.) if  $ D_{s,p}
((\beta_n)_{n=0}^{N-1})=O((\log N)^{s/2}) $ for $ N \rightarrow \infty $.

 A sequence  of point sets $((\beta_{n,N})_{n=0}^{N-1})_{N=1}^{\infty}$ is of
 $L_p$ low discrepancy (abbreviated
l.d.p.s.) if $ D_{s,p} ((\beta_{n,N})_{n=0}^{N-1})=O((\log
N)^{(s-1)/2}) $ for $ N \rightarrow \infty $.\\

The existence of $L_p$ l.d.p.s. was proved by Roth for $p=2$ and by Chen  for $p>1$ \cite{Ch}.
The first explicit construction of $L_p$ l.d.p.s. was obtained by Chen and Skriganov for $p=2$ and by Skriganov  for $p>1$ (see \cite{ChSk, Sk}).
The next explicit construction of $L_p$ l.d.p.s. was proposed by Dick and Pillichshammer
 (see \cite{Di, DiPi,  Ma}).
The first explicit construction of $L_p$ l.d.s. were obtained by Dick, Hinrichs, Markhasin and Pillichshammer \cite{DHMP}. All these explicit constructions was obtained by using $(t,m,s)$ nets. \\

In this paper we obtain a similar result for Halton's sequence.\\
Let $p_1,\ldots , p_s \geq 2 $ be pairwise coprime integers,
 \begin{equation}\label{In5}
 n=\sum_{j\geq 1}e_{p,j}(n) p_i^{j-1},\;  e_{i,j}(n) \in \{0,1, \ldots
 ,p_i-1\}, \;  {\rm and}    \; \phi_i(n)= \sum_{j\geq 1}e_{i,j}(n) p_i^{-j}.
 \end{equation}
Van der Corput    proved that $ (\phi_1(n))_{n\geq 0}$ is the $1-$dimensional l.d.s.\\
The first example of multidimensional l.d.s. was proposed by Halton
\begin{equation}  \label{In6}
  H_s(n)= (\phi_{1}(n),\ldots ,\phi_{s}(n)), \quad n=0,1,2,... \;.
\end{equation}
The first example of multidimensional l.d.p.s. was obtained by Hammersley
\begin{equation}    \label{In6a}
{\cH}_{s+1,N}=
	 (H_s(n), n/N)_{n=0}^{N-1}.
\end{equation}
\subsection{\textcolor{blue}{Main results}}
In this paper we will prove that Halton's sequence is of $L_p$ l.d.s.:\\ \\
{\bf Theorem 1.} {\it Let $s \geq 2, \; p \geq 1$. Then}
\begin{equation}\label{In9}
  D_{s,p} ((  H_{s}(k) )_{k=Q}^{Q+N-1}) =O(\log^{s /2} N )   ,
\end{equation}
where the $O$ constant  is independent of $Q$.\\ \\

%
For the sake of simplicity, we will consider only the case of primes $p_1,\ldots , p_s  $. For $n <0$, we consider \eqref{In5} in the sense of $\bbp-$adic representation.
Note that \eqref{In9} is also true for generalized Halton's sequences  (see, e.g., \cite{Le2}) and for the $s$-dimensional ergodic adding machine \cite{Le2}.

Similarly to \cite{Le4}, in Theorem 2, we prove that the local discrepancy of Hammersley's point net satisfies the Central Limit Theorem  (abbreviated CLT) for $s \geq 3$. This result is not true for $s=2$ because the normalised expectation $  E_{s+1}(D(\bar{\bx}, (\cH_{s+1,N} ) ) /\left\| D(\bar{\bx},\cH_{s+1,N}) \right\|_{s+1,2} $ does not vanish for $N \to \infty $, where $\bar{\bx} =(\bx,x_{s+1})=(x_1,...,x_{s+1})$. The simplest way to avoid this problem is to take
$  D(\bar{\bx}, \cH_{s+1,N} ) - E_{s+1}(D(\bar{\bx}, \cH_{s+1,N} ))$ instead of $D(\bar{\bx},\cH_{s+1,N} )$. But we prefer a different way. In Theorem 3, we get the asymptotic property of $L_p$ discrepancy of Hammersley's point net for $p>0$. Theorem 3 is the corollary of Theorem 1 and Theorem 2. For this reason, we want to prove CLT exactly  for  the discrepancy function. The  normalised expectation of the symmetrized Hammersley's set
$\cH_{s+1,N}^{sym} =(H_s( n), |n|/N)_{-N < n <N}$ vanishes for $N \to \infty $. So for $s=2$ we will take $\cH_{s+1,N}^{sym}$ instead of $\cH_{s+1,N}$.
 \\ \\
{\bf Theorem 2.} {\it Let $s \geq 2$,  $\bar{\bx}$ be a uniformly distributed random variable in $[0,1]^{s+1}$.  Then }
\begin{equation}  \nonumber 
 \frac{ D(\bar{\bx}, \cH_{s+1,N} ) }{ \left\| D(\bar{\bx},\cH_{s+1,N}) \right\|_{s+1,2}}
	\stackrel{w}{\rightarrow} \cN(0,1) \; \; \for \; s \geq 3,\quad \;
  \frac{ D(\bar{\bx},\cH_{3,N}^{sym} ) }{\left\| D(\bar{\bx},\cH_{3,N}^{sym}) \right\|_{3,2}}
	\stackrel{w}{\rightarrow} \cN(0,1).
\end{equation}
\\
{\bf Theorem 3.} {\it Let $s \geq 2$  and $p>0$. Then}
\begin{equation} \nonumber
  \frac{ D_{s+1,p}( \cH_{s+1,N} )}{ D_{s+1,2}( \cH_{s+1,N} )}
	\stackrel{N \rightarrow \infty}{\longrightarrow}   \kappa_p^{1/p}, \; s \geq 3,
			\qquad \quad \frac{ D_{3,p}( \cH_{3,N}^{sym} )}{D_{3,2}( \cH_{3,N}^{sym} )}
	\stackrel{N \rightarrow \infty}{\longrightarrow}   \kappa_p^{1/p},
\end{equation}
\begin{equation}\nonumber
\where \qquad \qquad
\kappa_p= \frac{1}{\sqrt{2 \pi}}\int_{-\infty}^{\infty} |u|^p e^{-u^2/2} d u, \quad \;\;\;\;
			\kappa_{2r} = \frac{(2r)!}{2^r r!} \;\;\; {\rm for \;integer}\; r \geq 1.
\end{equation}

Hence the lower bound \eqref{In7} is optimal for all $p>0$ for Hammersley's point set. \\

We note that Theorem 2 and Theorem 3 are also true for other symmetrization. For example, for
$\ddot{\cH}_{s+1,N}^{sym} = (H_s^{sym}(n), n/(2^sN))_{0 \leq n <2^sN}$, where
$H_s^{sym}(n) =(m_1 +(-1)^{m_1}\phi_1(m_{s+1}),...,m_s +(-1)^{m_s}\phi_s(m_{s+1}) ) $, with
$$n = m_1 +2m_2 +....+2^{s-1}m_s +2^s m_{s+1},\; m_i \in \{0,1\},\; i=1,...,s, \; m_{s+1} \geq 0.$$
For the case $s=1$ see e.g. \cite{ACS,LM}. Similar results are true for multidimensional  von Neumann-Kakutani ergodics
adding machine  (see Remark 3 at the end of the article).\\

{\bf Remark 1.}  The Halton sequence
is of $L_p$-low discrepancy for $s=1$ only for some sort of symmetrisation \cite{KrPi}. It is very difficult to understand why  Halton's sequence is of $L_p$-low discrepancy  for $s\geq 2$ without any symmetrisation. The first idea is that the main tool (Theorem A) may be applied only for $s\geq 2$. But it is not possible that this explanation is complete. Because  we see the same problem for transition  from $s = 2$ to $s \geq 3$. Namely, we do not need any symmetrisation for the validity of  the Central Limit Theorem  for
Hammersley's point set for $s \geq 3$ (see Theorem 2). But we need a symmetrisation for the case $s = 2$. Thus, it remains to say that probably the auto-symmetrisation grows with the increase  of the dimension. \\  \\

{\bf Remark 2.} Let $T\; : \; [0,1)^s \to  [0,1)^s$ be a map, $f \;: \; [0,1)^s \to \RR $ is a function,
 and $\bx_0 \in [0,1)^s$ is a fixed initial condition.
  We say that the ergodic sums $ S_k(f,\bx_0) =f(\bx_0) + f(T(\bx_0))+ \cdots
+ f(T^{k-1}(\bx_0))$ satisfy a {\it temporal
central limit theorem (TCLT) on the orbit of $\bx_0$}, if there exists
 constants $A_N$ and $B_N$ such that
\begin{equation} \nonumber 
  (S_k -A_N)/B_N   \stackrel{w}{\rightarrow} \cN(0,1),
\end{equation}
when $k$ is sampled uniformly from $\{ 1,...,N\}$.

The first example of the TCLT was discovered by Beck (see \cite{Be1, Be2}). Beck considered the transformation $x \to x +
 \sqrt{2} \mod 1$ and the indicator function
of an interval $[0,y)$. In \cite{DS1}, Dolgopyat and Sarig considered a more general notion of the temporal distributional limit
 theorem (see also  \cite{DS2} and references therein).

In \cite{Le6}, we proved an one-parametric temporal CLT
for an $s$-dimensional von Neumann-Kakutani ergodic
adding machine (Halton's sequence) with indicator functions. On the one hand, in this article we get a slightly worse result. Because
Theorem 2  and Theorem 2$^{'}$ (see below) are only $s+1$-parametric spatio-temporal CLT
 (see definition in \cite[p.12]{DS1}).
But, on the other hand, we get an advantage over the results in \cite{Le6}.
Since the results of this paper
 allow us to obtain an estimate for the $L_p$-discrepancy for $0 <p \leq 1$ (see Theorem 3  and Theorem 3$^{'}$).
More precisely. According to \cite[p.104] {Bi}
\begin{equation} \label{In77}
\emph{D}_{s,p}((\beta_{n,N}^{(s)})_{n=0}^{N-1})>C_{s,p}(\log N)^{\frac{s-1}{2}} \;\;\;
    {\rm for} \;\;\;  0 <p \leq 1.
\end{equation}
Hence, Theorem 3 support this conjecture. \\

Now we describe the structure of the paper. In Lemma 1, we get a simple estimate of Fourier's series of truncated discrepancy function of Halton's sequence.
In Lemma 2 and Lemma 3, we find an expression for the upper bound for the $ L_p $ discrepancy of the Holton sequence.
 In Lemma 4, we  minimise the number of terms in this expression.

 In Section 3, we apply the theorem on $\bbp$-adic logarithmic forms to obtain the first estimates of the $L_p$ discrepancy  of Halton's sequence. This is the main chapter of the paper.

 In Section 4, we finish the proof of Theorem 1 and, using the moment's method,
we prove Theorem 2. Next, using the standard tools of  probability theory, we derive Theorem 3 from Theorem 1 and Theorem 2. The main stages of the proof of Theorem 1 are almost the same as in \cite{Le5}.

The main tool of the proof is the theorem on linear forms in $\bbp$-adic logarithm (see \cite{Yu, Bu}).
\\ \\ \\

\section{\textcolor{blue}{Beginning of the proof of Theorems}}

\subsection{\textcolor{blue}{Preliminary lemmas}}
By the moment method, Theorem 2 follows from the following statement: \\  \\
{\bf Lemma A.} {\it Let $s \geq 2,\; \sh \geq 1$. With notations as above
\begin{equation}\label{Beg-0}
 \lim_{ N  \rightarrow \infty}
     \frac{ E_{s+1} \big( D^{\sh}(\bar{\bx}, \ddot{\cH}_{s+1,N} ) \big)}{ \left\| D(\bar{\bx},\ddot{\cH}_{s+1,N}) \right\|^{\sh}_{s+1,2}}  d \bar{\bx}= \begin{cases}
     \frac{\sh!}{2^{\sh/2}(\sh/2)!}  , & {\rm if} \; \sh \; {\rm is \; even},  \\
    0, & {\rm if} \; \sh \; {\rm is \; odd},
  \end{cases}
\end{equation}
where $\ddot{\cH}_{s+1,N}=\cH_{s+1,N}$ for $s \geq 3$ and $\ddot{\cH}_{s+1,N}=\cH_{s+1,N}^{sym}$ for $s =2$.}\\

The proof of the Lemma A is given below.\\

First of all, we need a Fourier series decomposition of the
discrepancy function (see Lemma 1).  We follow the approach of \cite[pp. 29-35]{Ni}. In order to formulate Lemma 1, we need the following notation and definitions.  We will use notation $A \ll B$  equal to $A = O(B)$.

Let
\begin{equation} \label{Beg-0a} 
  \delta_M(a) =   \begin{cases}
    1,  & \; {\rm if}  \;  a \equiv 0  \mod M,\\
    0, &{\rm otherwise},
  \end{cases} \qquad
   \b1(\fT) =   \begin{cases}
    1,  & \; {\rm if}  \;  \fT  \;{\rm is \;true},\\
    0, &{\rm otherwise}
  \end{cases} .
\end{equation}
Let $[y]$ be the integer part of $y$,
\begin{equation} \label{Beg-1a}
   I_M=[-[(M-1)/2],[M/2]] \cap \ZZ, \qquad  I_M^{*}=I_M \setminus \{ 0 \}.
\end{equation}
 Note that  the integers of the interval $I_M$ are a complete set of residues $\mod M$, $M \geq 1$.
By \cite[Lemma 2, p. 2]{Ko}, we have
\begin{equation} \label{Beg-1}
   \delta_M(a) =   \frac{1}{M} \sum_{k\in I_M} e\Big(\frac{ak}{M}\Big),  \qquad \where \quad e(x) =exp(2\pi i x).
\end{equation}
By \cite[Lemma 1, p. 1]{Ko}, we get
\begin{equation} \label{Beg-2a}
  \Big|\frac{1}{M} \sum_{k=0}^{M-1} e\big(k \alpha) \Big| \leq \min\big(1, \frac{1}{2M\llangle \alpha \rrangle} \big) ,\quad \;\; \where \;
	\llangle \alpha \rrangle = \min( \{\alpha\}, 1 - \{\alpha\}).
\end{equation}
Let $\bar{m} = \max(1, |m|) \leq M/2 $.     From \cite[p. 2]{Ko}, we obtain for $R \leq M$ :
\begin{equation} \label{Beg-2}
  \Big|\frac{1}{M} \sum_{k=0}^{R-1} e\big(\frac{m k}{M}\big) \Big| \leq \min\Big(1, \Big|\frac{e(m R/M) -1}{M
	(e(m /M)-1)} \Big|\Big) \leq \frac{|\sin(\pi m R/M)|}{\bar{m}}  \leq \frac{1}{\bar{m}}.
\end{equation}
Let $x_i =0.x_{i,1}x_{i,2}...=\sum_{j \geq 1}  x_{i,j} p_{i}^{-j}$, with $x_{i,j} \in \{0,1,...,p_{i}-1  \}$,
 $i=1,...,s$.
We define the truncation
\begin{equation}  \label{Beg-2-0}
        [x_i]_r =\sum_{1 \leq j \leq r} x_{i,j} p_{i}^{-j} \quad \with \quad r \geq 1.
\end{equation}
If $\bx = (x_1, . . . , x_s)  \in [0, 1)^s$, then the truncation $[\bx]_{\br}$ is defined coordinatewise, that is $[\bx]_{\br}=
( [x_1]_{r_1}, . . . , [x_s]_{r_s})$, where $\br =(r_1,...,r_s)$.
By (\ref{In5}), we have
\begin{equation}\nonumber
    [\phi_{i}(k)]_{r}=[x_i]_{r} \;
			\Leftrightarrow  \; k \equiv    \sum_{1 \leq j \leq r}  x_{i,j} p_{i}^{j-1}  \;\; ({\rm mod} \; p_{i}^{r}).
\end{equation}
 Let  $p_0=p_1 p_2\cdots p_s$, $ P_{\br} =p_1^{r_1}p_2^{r_2} \dots p_s^{r_s}$
  and let $  M_{i,\br}$ be the unique
integer satisfying the two conditions
\begin{equation} \label{Beg-3} 
   M_{i,\br} \equiv
	\big( P_{\br}/p_i^{r_i} \big)^{-1} \mod p_i^{r_i}, \; M_{i,\br} \in [0,p_i^{r_i}), \; 1 \leq i \leq s.
\end{equation}
Applying (\ref{In6}) and the Chinese Remainder Theorem, we get
\begin{equation} \label{Beg-4} 
  [H_s(k)]_{\br} = [\bx]_{\br}  \; \Longleftrightarrow  \; k \equiv   \cX_{\br}  \; ({\rm mod} \; P_{\br}),
\end{equation}
where
\begin{equation}\label{Beg-6} 
 \cX_{\br} \equiv \sum_{i=1}^s  M_{i,\br}
   P_{\br}  p_i^{-r_i} \cX_{i,r_i} \;(\mod \; P_{\br}), \;\; \cX_{i,r_i}=
  \sum_{1 \leq j \leq r}  x_{i,j} p_i^{j-1}, \;  \;
	\cX_{\br} \in [0,P_{\br}),
\end{equation}
$\cX_{i,r_i} \in [0,p_i^{r_i}),$  $1 \leq i \leq s$. \\

Let $n =[\log_2 N]+1$, $\bn =(n,n,...,n)$.
From \cite[p. 29, 30]{Ni}, we get
\begin{equation}  \label{Beg-16}
  \rD(Q,N) :=  D([\bx]_{\bn}, (H(k))_{k=Q}^{Q+N-1}  ) =  D(\bx, (H(k))_{k=Q}^{Q+N-1}  ) + \epsilon s,  \;\;\; |\epsilon| \leq 1.
\end{equation}
Let
\begin{equation}  \label{Beg-18}
  \cX_{\br,\bb} \equiv \sum_{i=1}^s  M_{i,\br}
   P_{\br}  p_i^{-r_i}
  \Big( \sum_{1 \leq j < r_i}  x_{i,j} p_i^{j-1}  + b_i p_i^{r_i-1} \Big) \mod  P_{\br} , \quad
	\cX_{\br,\bb} \in [0,P_{\br}).
\end{equation}
Using \eqref{In1}, we obtain
\begin{equation} \nonumber
 \rD(Q,N)=  \sum_{k=Q}^{Q+N-1} \sum_{r_1,...,r_s =1}^n \sum_{b_1=0}^{x_{1,r_1}-1} \cdots  \sum_{b_s=0}^{x_{s,r_s} -1} \delta_{P_{\br}} (k - \cX_{\br,\bb})- N [x_1]_n \cdots [x_s]_n \\
\end{equation}
\begin{multline}  \label{Beg-20}
    =\sum_{r_1,...,r_s =1}^n  \ddot{\rD}_{Q,N,\br}, \;\;\;
    \ddot{\rD}_{Q,N,\br}:= \sum_{b_1=0}^{x_{1,r_1}-1} \cdots  \sum_{b_s=0}^{x_{s,r_s} -1}  \sum_{k=Q}^{Q+N-1} (\delta_{P_{\br}} (k - \cX_{\br,\bb})-1/P_{\br} ),\\
            |\ddot{\rD}_{Q,N,\br}| \leq x_{1,r_1}\cdots x_{s,r_s}  <p_0.
\end{multline}
For given $\br$, we define
\begin{equation}  \label{Beg-22}
   \ddot{T}(\br) := \{  i \in [1,s] \; | \; r_i >V_1 \}, \qquad V_1= [\log_2^3 n].
\end{equation}
Let $\fs \in [1,s]$. We consider subsets $T_{\fs} \subseteq \{1,...,s\}$ with $\card(T_{\fs}) =\fs$.
It is easy to see that
\begin{equation} \nonumber 
   \sum_{\fs=0}^s  \sum_{T_{\fs} \subseteq \{1,...,s\}, \#T_{\fs}=\fs } \b1 (T_{\fs} =\ddot{T}(\br))=1.
\end{equation}
By \eqref{Beg-20}, we get
\begin{multline}  \label{Beg-23}
 \rD(Q,N)    =\sum_{r_1,...,r_s =1}^n
     \sum_{\fs=0}^s  \sum_{T_{\fs} \subseteq \{1,...,s\}} \b1 (T_{\fs} =\ddot{T}(\br)) \ddot{\rD}_{Q,N,\br} \\
 = \sum_{\fs=0}^s  \sum_{T_{\fs} \subseteq\{1,...,s\} } \fD_{T_{\fs}}(Q,N), \;\;\;\; \with \;\;\;\;
 \fD_{T_{\fs}}(Q,N)=  \sum_{\br \in \UU_{T_{\fs}}}
  \ddot{\rD}_{Q,N,\br},
\end{multline}
where
\begin{equation}   \label{Beg-26}
  \UU_{T_{\fs}}= \{ \br=(r_1,...,r_s) \in [1,n]^s \; | \; r_i > V_1,\; \for \; i \in T_{\fs},\;\;\; r_i \leq V_1,\; \for \; i \not \in T_{\fs}\}.
\end{equation}
From \eqref{Beg-20} and  \eqref{Beg-22}, we obtain
\begin{equation}  \nonumber 
 |\fD_{T_0}(Q,N)| \leq  p_0V_1^s \leq p_0 \log^{3s} n.
\end{equation}
Using  \eqref{Beg-16} and \eqref{Beg-23}, we derive
\begin{multline}  \label{Beg-28}
 |\tilde{\rD}(Q,N)  |  \leq \sum_{\fs=0}^{s-1}  \sum_{T_{\fs} \subseteq\{1,...,s\} }
    | \fD_{T_{\fs}}(Q,N) | = \sum_{\fs=1}^{s-1}  \sum_{T_{\fs} \subseteq\{1,...,s\} }
    | \fD_{T_{\fs}}(Q,N) |  +O( \log^{3s} n),\\
    \with  \qquad \tilde{\rD}(Q,N):= D(\bx, (H(k))_{k=Q}^{Q+N-1}  )  -\fD_{T_{s}}(Q,N).
\end{multline}\\ \\
{\bf Lemma 1.} {\it With the notations as above, we have}
\begin{equation*}
  \fD_{T_{\fs}}(Q,N) =
 \sum_{\br \in U( T_{\fs})}   \sum_{m \in I_{P_{\br}}^{*}}   \varphi_{ \br,Q,N,m} \; \psi_{ \br}(m,\bx)\;
e\Big( \frac{-m}{P_{\br}} \cX_{\br} \Big),
\end{equation*}
\begin{equation*}
\varphi_{ \br,Q,N,m}
  =\frac{e(m(Q+N)/P_{\br}) - e(m Q /P_{\br} )}{P_{\br}(e(m/P_{\br})-1)},  \;\;\;
  |\varphi_{ \br,Q,N,m}| \leq \frac{1}{\bar{m}},
\end{equation*}
 \begin{equation*}
        \psi_{ \br}(m,\bx)  = \prod_{i=1}^s \ddot{\psi}(i,\{-m M_{i,\br}/p_i \}p_i, x_{i, r_i}) ,\quad |\psi_{ \br}(m,\bx)| \leq p_0,
\end{equation*}
\begin{equation} \label{Le2-1}
\ddot{\psi}(i,0,x_{i, r_i})=x_{i, r_i}, \quad
				\ddot{\psi}(i,m', x_{i, r_i}) = \frac{1- e(-m' x_{i, r_i}/p_i)}{e(m'/p_i)-1}\;\;	 \for\; m' \neq 0.
\end{equation} \\

{\bf Proof.}
Similarly to \cite[p. 37-39]{Ni}, we obtain from \eqref{Beg-1}, \eqref{Beg-1a} and \eqref{Beg-20} that
\begin{multline*}
\ddot{\rD}_{Q,N,\br}=    \sum_{b_1=0}^{x_{1,r_1} -1} \cdots  \sum_{b_s=0}^{x_{s,r_s} -1} \sum_{k=Q}^{Q+N-1} \frac{1}{P_{\br}} \sum_{m \in I_{P_{\br}}^{*}} e\Big( \frac{m}{P_{\br}} (k - \cX_{\br,\bb})\Big) \\
=  \sum_{b_1=0}^{x_{1,r_1} -1} \cdots  \sum_{b_s=0}^{x_{s,r_s} -1}  \sum_{m \in I_{P_{\br}}^{*}}
\frac{e(m(Q+N)/P_{\br}) - e(m Q /P_{\br} )}{P_{\br}(e(m/P_{\br})-1)}
\; e\Big( \frac{-m}{P_{\br}} \cX_{\br,\bb} \Big) .
\end{multline*}
According to \eqref{Beg-18}, we get
\begin{equation*}  
  \cX_{\br,\bb} \equiv \cX_{\br} + \sum_{i=1}^s  M_{i,\br}
   P_{\br}  (b_i - x_{i,r_i})/p_i  \mod  P_{\br}.
\end{equation*}
Using  \eqref{Le2-1}, we get
\begin{multline*}
\ddot{\rD}_{Q,N,\br} = \sum_{m \in I_{P_{\br}}^{*}} e\Big( \frac{-m}{P_{\br}}
\cX_{\br} \Big) \varphi_{ \br,Q,N,m}
\sum_{b_1=0}^{x_{1,r_1} -1} \cdots  \sum_{b_s=0}^{x_{s,r_s} -1}
e\Big(  -m\sum_{i=1}^s  M_{i,\br}(b_i-x_{i,r_i})/p_i \Big)\\
 =    \sum_{m \in I_{P_{\br}}^{*}}   \varphi_{ \br,Q,N,m} \; \psi_{ \br}(m,\bx) \;
e\Big( \frac{-m}{P_{\br}} \cX_{\br} \Big).
\end{multline*}
By \eqref{Beg-2} and \eqref{Le2-1}, we obtain
\begin{equation}  \label{Le2-4}
 |\varphi_{ \br,Q,N,m}  | \leq \frac{1}{\bar{m}}, \quad |\ddot{\psi}(i,m, x_{i, r_i})| \leq p_i,\quad |\psi_{ \br}(m,\bx)| \leq p_0, \;\; i=1,...,s.
\end{equation}
Bearing in mind \eqref{Beg-23}, we get the assertion of Lemma 1. \qed \\

Let $q =[\sh/2]$, $\Xi_{\sh}$ be the set of all permutations of the set $ \{ 1,...,\sh \} $,
\begin{multline}  \label{Lem2-6}
 \varpi_{\br,\bm,2} = \b1 \big(\sh=2q, \; \fs=s,\; \exists \eta \in \Xi_{2q}  \; : \;  m_{\eta(2k-1)}/P_{\br_{\eta(2k-1)}}  =- m_{\eta(2k)}/P_{\br_{\eta(2k)}},\\
  k \in [1,q]  \big),\quad \varpi_{\br,\bm,1} = 1- \varpi_{\br,\bm,2}, \quad \nu \in \{ 0,1 \},
\end{multline}
\begin{equation}  \nonumber
  \sD_{T_{\fs,\sh,\nu}}(Q,N)
	=\sum_{\substack{\br_j \in \UU_{ T_{\fs}}  \\ j \in [1,\sh]}}
\sum_{\substack{m_j \in I^{*}_{P_{\br_j}} \\ j \in [1,\sh]}}  \varpi_{\br,\bm,\nu} \prod_{j=1}^{\sh}
  \varphi_{ \br_j,Q,N,m_j} \; \psi_{ \br_j}(m_j,\bx)\;
e\Big( \frac{-m_j}{P_{\br_j}} \; \cX_{\br_j} \Big).
\end{equation}
From Lemma 1, we have
\begin{equation}  \label{Lem2-7}
  \fD_{T_{\fs}}^{\sh}(Q,N)  = \sD_{T_{\fs,\sh,1}}(Q,N) + \sD_{T_{\fs,\sh,2}}(Q,N).
\end{equation}
Let $ \varpi_{\br,\bm,2}=1$ and let $m_j =\ddot{m}_jP_{\balpha_j}$ with $(\ddot{m}_j,p_0)=1$.
Therefore $\mu_k :=\ddot{m}_{\eta(2k)} =- \ddot{m}_{\eta(2k-1)}$ and $\br_{\eta(2k)} -\balpha_{\eta(2k)}=
\br_{\eta(2k-1)} -\balpha_{\eta(2k-1)}$.

Bearing in mind that $ |\varphi_{ \br_j,Q,N,m_j} \psi_{ \br_j}(m_j,\bx) | \leq p_0/\bar{m}_j$, we get
\begin{multline}  \label{Le3-2}
  \sD_{T_{\fs,2q,2}}(Q,N)\ll
  \sum_{\substack{\alpha_{i,j} \in [0,p_0 s n ]  \\i \in [1,s], j \in [1,2q]}}
  \sum_{\substack{\br_j \in \UU_{ T_{s}}  \\ j \in [1,2q]}}
\sum_{\substack{\mu_j \in I^{*}_{P_{\br_j}} \\ j \in [1,q]}}  \frac{\b1(\br_{k} -\balpha_{k}=
\br_{k +q} -\balpha_{k+q}, k=1,...,q )}{ P_{\balpha_1} \cdots  P_{\balpha_{2q}}\mu_1^2 \cdots \mu_q^2 }
 \\
\ll     \sum_{\substack{\alpha_{i,j} \in [0,p_0 s n ] , i \in [1,s], j \in [1,2q]}}
\frac{ n^{qs } }{ P_{\balpha_1} \cdots  P_{\balpha_{2q}}} \ll n^{qs }.
\end{multline}
We will get  a more precise estimate  in $\S4$.
Using  \eqref{Lem2-7} and  \eqref{Le3-2}, in  \S4.1  we show that in order to prove Theorem 1 it is enough to verify that
\begin{equation}  \label{Le3-3}
E_s (\sD_{T_{\fs,\sh,1}}(Q,N))  \ll  n^{\sh s/2 }.
\end{equation}
%
Let $E^{(i)}(f(x_i))= \int_{[0,1)} f(x_i)d x_i  $.
 It is easy to see that
\begin{equation} \label{Le3-9}
 E^{(i)}(f([x_i]_l))=  \frac{1}{p_i^l} \sum_{x_{i,1}=0}^{p_i-1}...\sum_{x_{i,l}=0}^{p_i-1} f(0.x_{i,1} ... x_{i,l}).
\end{equation}
Now, applying Lemma 1, we get: \\ \\
{\bf Lemma 2.} {\it With the notations as above, we get}
\begin{equation} \nonumber
 |E_s( \sD_{T_{\fs,\sh,1}}(Q,N)) | \ll
\sum_{\substack{\br_j \in \UU_{ T_{\fs}}  \\ j=1,...,\sh}} \;
\sum_{\substack{m_j \in I_{P_{\br}}^{*} \\  j=1,...,\sh}} \;
\frac{\varpi_{\br,\bm,1}}{\bar{m}_1 \cdots \bar{m}_{\sh}}  \prod_{i=1}^s |  E^{(i)}(Z_i)|
\end{equation}
and
\begin{equation} \nonumber
 |E_{s+1}( \sD_{T_{\fs,\sh,1}}(-[Nx_{s+1}], 2[Nx_{s+1}])) | \\
 \ll \sum_{\substack{\br_j \in \UU_{ T_{\fs}}  \\ j=1,...,\sh}} \;
\sum_{\substack{m_j \in I_{P_{\br}}^{*} \\  j=1,...,\sh}} \;
\frac{\varpi_{\br,\bm,1}  |\hat{\gamma}^{(\sh)}_{\br,\bm}| }{\bar{m}_1 \cdots \bar{m}_{\sh}} \;
  \prod_{i=1}^s |  E^{(i)}(Z_i)| ,
\end{equation}
with
\begin{multline}  \label{Le3-10}
\hat{\gamma}^{(\sh)}_{\br,\bm}:= \bar{m}_1 \cdots \bar{m}_{\sh} \int_0^1 \prod_{j=1}^{\sh} \varphi_{ \br_j,-[Nx_{s+1}], 2[Nx_{s+1}],m_j} d x_{s+1},\\
\hat{m}_{i}  \equiv  - \sum_{j=1}^{\sh}   m_j M_{i,\br_j}
   p_i^{r^{+}_i -r_{i,j} } \mod p_i^{r^{+}_i}, \quad
       \hat{m}_{i} \in [0,p_i^{r^{+}_i}), \;\;   r^{+}_i :=\max_j  r_{i,j},\; i \in [1,s], \\
   \ad \quad Z_i=  e\Big(   \sum_{i=1}^s \frac{\hat{m}_{i} \cX_{i,r^{+}_i}}{p_i^{r^{+}_i}}
	\Big)   \prod_{j=1}^{\sh}
     \ddot{\psi}(i ,\{-m_{j} M_{i,\br_{j}} /p_i\}p_i , x_{i,r_{i,j}}  ).
\end{multline}
{\bf Proof.}
We will prove  the first relation. The proof of the second relation is similar.
By   \eqref{In2}, \eqref{Le2-4}  and \eqref{Lem2-6}, we obtain
\begin{multline*}  
  |E_s( \sD_{T_{\fs,\sh,1}}(Q,N)) |   \leq
	\sum_{\substack{\br_j \in \UU( T_{\fs})  \\ j=1,...,\sh}}
\sum_{\substack{m_j \in I^{*}_{P_{\br_j}} \\ j=1,...,\sh}} \varpi_{\br,\bm,1} \prod_{j=1}^{\sh}
 |\varphi_{\br_j,Q,N,m_j} \; E_s(Z)|   \\
\leq
	\sum_{\substack{\br_j \in \UU( T_{\fs})  \\ j=1,...,\sh}}
\sum_{\substack{m_j \in I_{P_{\br}}^{*} \\  j=1,...,\sh}}
 \frac{\varpi_{\br,\bm,1}}{ \bar{m}_1 \cdots  \bar{m}_{\sh}} | E_s(Z)|, \quad \where \quad
   Z= \prod_{j=1}^{\sh}
 \psi_{\br_{j}}(m_{j}  ,\bx)
  e\Big(-\sum_{j=1}^{\sh}  \frac{m_j }{P_{\br_j}}\cX_{\br_j}  \Big).
\end{multline*}
From \eqref{Beg-6} and \eqref{Le3-10}, we have
\begin{multline}   
 e\Big(-\sum_{j=1}^{\sh}  \frac{m_j }{P_{\br_j}}\cX_{\br_j}  \Big)
	  =
 e\Big(-\sum_{i=1}^s  \sum_{j=1}^{\sh}  \frac{m_j M_{i,\br_j} \cX_{i,r_{i,j}}}{p_i^{r_{i,j}}}  \Big)
  	  \\
=   e\Big(-\sum_{i=1}^s  \sum_{j=1}^{\sh}   \frac{m_j M_{i,\br_j} p_i^{r^{+}_i -r_{i,j}}
 }{p_i^{r^{+}_i}}
  \cX_{i,r_{i}} \Big)
 =   e\Big(  \sum_{i=1}^s \frac{\hat{m}_{i} \cX_{i,r^{+}_i}}{p_i^{r^{+}_i}}
	\Big).
\end{multline}
In view of \eqref{Le2-1} and  \eqref{Le3-10}, we get
$ E_s(Z)= E^{(1)}(Z_1)  \cdots  E^{(s)}(Z_s)$.

Hence, Lemma 2 is proved. \qed \\

According to \eqref{Le3-3}, we need to calculate the expectations $ E_s( \sD_{T_{\fs,\sh,1}}(Q,N)) $. Taking Lemma 2 into account, it is enough to calculate $ E^{(i)}(Z_i)$. We will do this in Lemma 3. First, we need to rearrange the sequences $r_{i,j}$ in ascending order. Note that the main steps  of calculations in this section are the same as in \cite{Le5}, but more complicated.

Let $\sigma_{i,r}$   be a permutation of the set $ \{1,2,...,\sh \}$
satisfies the condition
\begin{equation}\label{Le3-17}
                r_{i,\sigma_{i,r}(j-1)} \leq  r_{i,\sigma_{i,r}(j)}, \;\;\;\; r_{i,\sigma_{i,r}(\sh)} = r^{+}_i :=\max_j  r_{i,j}, \quad   r_{i,0}: =0, \; \sigma_{i,r}(0): =0,
\end{equation}
$\br_j=(r_{1,j},...,r_{s,j})$ for $j \in [1,\sh], \; i=[1,s] $.

From \eqref{Le3-10},  we derive : $ \hat{m}_{i} \in [0,p_i^{r^{+}_i})  $,
\begin{multline}  \label{Le3-18}
-\hat{m}_{i}  \equiv   \sum_{j=1}^{\sh}   m_j M_{i,\br_j}
   p_i^{r^{+}_i -r_{i,j} }
   \equiv \sum_{j=1}^{\sh}   m_{\sigma_{i,r}(j)} M_{i, \br_{\sigma_{i,r}(j)}}  p_i^{r^{+}_i -r_{i,\sigma_{i,r}(j)}}
   \mod p_i^{r^{+}_i}   .
\end{multline}
Let $\mm_{i,j}$   $(i=1,...,s, \; j=1,...,\sh)$ be the unique integer satisfing the following congruence:
\begin{equation} \label{Le3-20}
\hat{m}_{i}  \equiv \sum_{j=1}^{\sh}  \mm_{i,\sigma_{i,r}(j)}p_i^{r^{+}_i -r_{i,\sigma_{i,r}(j)}  }
  \mod p_i^{r^{+}_i}, \;     \mm_{i,\sigma_{i,r}(j)} \in I_{p_i^{r_{i,\sigma_{i,r}(j) }-r_{i,\sigma_{i,r}(j-1)}} },\; I_{1} = I_{p_i^0}=\{ 0 \}.
\end{equation}

To prove Lemma 3, we apply inequality \eqref{Beg-2a}  separately to the intervals
 $[p_i^{r_{i,\sigma_{i,r}(j-1)}},p_i^{r_{i,\sigma_{i,r}(j)}-1} ] $  $(i,j \geq 1)$.
Note that Lemma 3 is the main part for estimating the $ L_p $ discrepancy of Holton's sequence. \\ \\
%
{\bf Lemma 3.} {\it With notations as above, we have}

\begin{equation}  \nonumber 
     | E^{(i)}(Z_i)| \leq  \prod_{j=1}^{\sh} \frac{4 p_i^{\sh +3}}{\bar{\mm}_{i,\sigma_{i,r}(j)}}
.
\end{equation} \\
{\bf Proof.}
From \eqref{Le3-9}, \eqref{Le3-10} and \eqref{Beg-6}, we get
$  \cX_{i,r_i}=
  \sum_{1 \leq j \leq r}  x_{i,j} p_i^{j-1}, \;  \; $ and
\begin{multline*}
    E^{(i)}(\hat{Z}_i) =      \frac{1}{p_i^{r^{+}_i }}   \sum_{x_{i,k} \in [0,p_i),k \in [1,r^{+}_i]}       e\left(  \frac{\hat{m}_{i} \cX_{i,r^{+}_i}}{p_i^{r^{+}_i}}
	\right)   \prod_{j=1}^{\sh}      \ddot{\psi}(i ,\{-m_{j} M_{i,\br_{j}} /p_i\}p_i , x_{i,r_{i,j}}  ) \\
=      \frac{1}{p_i^{r^{+}_i }}   \sum_{x_{i,k} \in [0,p_i),k \in [1,r^{+}_i]}        e\left(  \frac{\hat{m}_{i} \sum_{k=1}^{r^{+}_i} x_{i,k}p_i^{k-1}}{p_i^{r^{+}_i}}
	\right)    \prod_{j=1}^{\sh}      \ddot{\psi}(i ,\{-m_{j} M_{i,\br_{j}} /p_i\}p_i , x_{i,r_{i,j}}  ).
\end{multline*}
Let's consider the partition  \eqref{Le3-17} of the interval $[1, r_i^{+}]$ :
\begin{multline*}
    E^{(i)}(Z_i) =   \prod_{j=1}^{\sh}  \frac{1}{p_i^{r_{i,\sigma_{i,r}(j) }-r_{i,\sigma_{i,r}(j-1)}}}
    \sum_{\substack{x_{i,k} \in [0,p_i)\\
k \in [r_{i,\sigma_{i,r}(j-1) }+1,r_{i,\sigma_{i,r}(j)}]}} 1\\
 \times       e\left(  \frac{\hat{m}_{i} \sum_{k=r_{i,\sigma_{i,r}(j-1)}+1}^{r_{i,\sigma_{i,r}(j)} } x_{i,k}p_i^{k-1}}{p_i^{r^{+}_i}}	\right)
 \ddot{\psi}(i ,\{-m_{j} M_{i,\br_{j}} /p_i\}p_i , x_{i,r_{i,j}}  )  , \;\;\;\; \with \; r_{i,\sigma_{i,r}(0)}=0.
\end{multline*}
Hence
\begin{equation}   \label{Le4-2}
| E_i(Z_i)|   \leq \prod_{j =1}^{\sh} Z_{i,j},
\end{equation}
where
\begin{multline}   \label{Le4-22}
     Z_{i,j}   =\frac{1}{p_i^{r_{i,\sigma_{i,r}(j) }-r_{i,\sigma_{i,r}(j-1)}} }
\Bigg\vert \sum_{\substack{x_{i,k} \in [0,p_i)\\
k \in [r_{i,\sigma_{i,r}(j-1) }+1,r_{i,\sigma_{i,r}(j)}]}}
        e\Bigg(  \frac{\hat{m}_{i} \sum_{k=r_{i,\sigma_{i,r}(j-1)}+1}^{r_{i,\sigma_{i,r}(j)} } x_{i,k}p_i^{k-1}}{p_i^{r^{+}_i}}
	\Bigg)  \\
 \times \ddot{\psi}(i ,\hat{\mu}_{i,j} , x_{i,r_{i,\sigma_{i,r}(j)}} ) \Bigg\vert \quad \with \; {\rm some} \;\hat{\mu}_{i,j}\in \{1,p_i-1 \}.
\end{multline}
Let
\begin{equation*}
  J=\{ j \in [1,\sh] \; | \; r_{i,\sigma_{i,r}(j-1)} \leq r_{i,\sigma_{i,r}(j) } \leq r_{i,\sigma_{i,r}(j-1)} + \sh +1\},
  \;\; \;\;\sh^{'} =\card(J)  .
\end{equation*}
Let's consider the case  $j \in J$.
From \eqref{Le4-22} and \eqref{Le3-20}, we get
\begin{equation*}
  |\ddot{\psi}(i,m, x_{i, r_i})| \leq p_i, \quad  Z_{i,j} \leq p_i, \quad \ad \quad
     1 \leq \bar{\mm}_{i,\sigma_{i,r}(j)} \leq p_i^{r_{i,\sigma_{i,r}(j)}-r_{i,\sigma_{i,r}(j-1)}} \leq p_i^{\sh +1}.
\end{equation*}
Taking  $Z_{i,j}^{'}=p_i^2 $ for $ r_{i,\sigma_{i,r}(j)} \leq r_{i,\sigma_{i,r}(j-1)}+\sh +1 $,
 we have
\begin{equation} \label{Le2-60}
  Z_{i,j} \leq p_i <  Z_{i,j}^{'}=p_i^2  \leq \frac{p_i^{r_{i,\sigma_{i,r}(j)}-r_{i,\sigma_{i,r}(j-1)}+2}}{\bar{\mm}_{i,\sigma_{i,r}(j)}  }
      \leq  \frac{p_i^{\sh+3}}{\bar{\mm}_{i,\sigma_{i,r}(j)} } \;.
\end{equation}
Let's consider the case $j \in [1, \sh]\setminus J$.    (i.e.,$ r_{i,\sigma_{i,r}(j-1) } \leq r_{i,\sigma_{i,r}(j)} -2$).

  Bearing in mind that $|\ddot{\psi}(i,m, x_{i, r_{i,j}})| \leq p_i  $, we obtain from \eqref{Le4-22} that
\begin{multline}   \label{Le4-22a}
          Z_{i,j}  \leq Z_{i,j}^{'}, \qquad \qquad \where\\
     Z_{i,j}^{'}   =\frac{1}{p_i^{r_{i,\sigma_{i,r}(j) }-r_{i,\sigma_{i,r}(j-1)}-2}}
\left| \sum_{\substack{x_{i,k} \in [0,p_i)\\
k \in [r_{i,\sigma_{i,r}(j-1) }+1,r_{i,\sigma_{i,r}(j)}-1]}}
        e\Bigg(  \frac{\hat{m}_{i} \sum_{k=r_{i,\sigma_{i,r}(j-1)}+1}^{r_{i,\sigma_{i,r}(j)} } x_{i,k}p_i^{k-1}}{p_i^{r^{+}_i}}
	\Bigg) \right| .
\end{multline}

Taking into account \eqref{Le4-2} and \eqref{Le2-60}, we get that in order to prove Lemma~3 it is enough to verify the following inequality
\begin{equation}  \label{Le2-25}
 Z_{i,j}^{'} \leq
 \frac{p_i^{\sh+3}}{\bar{\mm}_{i, \sigma_{i,r}(j)}} \quad  \for \quad r_{i,\sigma_{i,r}(j-1) } \leq r_{i,\sigma_{i,r}(j)} -\sh -2.
\end{equation}
%
%
 Applying  \eqref{Beg-2a},  we derive
\begin{equation} \label{Le2-70}
   Z_{i,j}^{'} \leq p_i \min \Bigg(1,   \frac{1}{2p_i^{r_{i,\sigma_{i,r}(j)} -r_{i,\sigma_{i,r}(j-1)}-1}
		\big\llangle \hat{m}_{i}  /p_i^{r^{+}_i -r_{i,\sigma_{i,r}(j-1) }} \big\rrangle   } \Bigg).
\end{equation}
 If $|\mm_{i, \sigma_{i,r}(j)}| \leq 2 \sh$, then we will use the trivial estimate
$$
 Z^{'}_{i,j} \leq p_i \leq 2 p_i\sh/\bar{\mm}_{i,\sigma_{i,r}(j)} \leq  p_i^{\sh+3}/\bar{\mm}_{i,\sigma_{i,r}(j)}.
$$
Now let's consider the case $ |\mm_{i, \sigma_{i,r}(j)}| > 2 \sh$.

By   \eqref{Le2-70}, we get that in order to prove \eqref{Le2-25} it is enough to verify that
\begin{equation}\label{Le2-80}
   \left\llangle \frac{\hat{m}_{i} } {p_i^{r^{+}_i -r_{i,\sigma_{i,r}(j-1) }} }
    \right\rrangle  \geq
     \frac{|\mm_{i,\sigma_{i,r}(j)} | /2}{p_i^{ r_{i,\sigma_{i,r}(j) } -r_{i,\sigma_{i,r}(j-1)}} }.
\end{equation}
From \eqref{Le3-20}, \eqref{Beg-1a} and the previous conditions,  we obtain $r_{i,\sigma_{i,r}(j) } \geq r_{i,\sigma_{i,r}(j-1)}+\sh +2 $,
$2\sh <|\mm_{i, \sigma_{i,r}(j)}| \leq  \frac{1}{2} p_i^{r_{i,\sigma_{i,r}(j) } -r_{i,\sigma_{i,r}(j-1) }}$, $j \in [1,h] \setminus J$ and
\begin{equation}  \nonumber
  \left\llangle \frac{\hat{m}_{i} }{p_i^{r^{+}_i -r_{i,\sigma_{i,r}(j-1) }} }
    \right\rrangle = \left\llangle   \sum_{\ell=1}^{\sh}   \frac{ \mm_{i, \sigma_{i,r}(\ell)} p_i^{r^{+}_i -r_{i,\sigma_{i,r}(\ell)}}}  {p_i^{r^{+}_i -r_{i,\sigma_{i,r}(j-1) }} }
    \right\rrangle  =
\left\llangle   \sum_{\ell=1}^{\sh}   \frac{ \mm_{i, \sigma_{i,r}(\ell)}}{p_i^{r_{i,\sigma_{i,r}(\ell)}} -r_{i,\sigma_{i,r}(j-1) }}
    \right\rrangle
\end{equation}
\begin{equation}  \nonumber
=\left\llangle   \sum_{\ell=j}^{\sh}   \frac{ \mm_{i, \sigma_{i,r}(\ell)} }{p_i^{r_{i\sigma_{i,r}(\ell)}} -r_{i,\sigma_{i,r}(j-1) }}
    \right\rrangle   =
\left\llangle   \frac{ \mm_{i, \sigma_{i,r}(j)} } {p_i^{r_{i,\sigma_{i,r}(j)} -r_{i,\sigma_{i,r}(j-1) }} } +    \sum_{ \ell =j+1}^{\sh}  \frac{ \mm_{i, \sigma_{i,r}(\ell)}}  {p_i^{r_{i,\sigma_{i,r}(\ell)} -r_{i,\sigma_{i,r}(j-1) }} }
    \right\rrangle
\end{equation}
\begin{equation}  \nonumber
=\left\llangle   \frac{ \mm_{i, \sigma_{i,r}(j)} + \varepsilon_{i} }{p_i^{r_{i,\sigma_{i,r}(j)} -r_{i,\sigma_{i,r}(j-1) }} }     \right\rrangle \quad \with \;\;
\varepsilon_{i} =
\sum_{ \ell =j+1}^{\sh}   \frac{ \mm_{i, \sigma_{i,r}(\ell)} } {p_i^{(r_{i,\sigma_{i,r}(\ell)}-r_{i,\sigma_{i,r}(j-1) }) - (r_{i,\sigma_{i,r}(j)} -r_{i,\sigma_{i,r}(j-1) })} }.
\end{equation}
It is easy to see that $|\varepsilon_{i}| \leq \sh$.
Therefore
\begin{equation}  \nonumber
 \left\llangle \frac{\hat{m}_{i} } {p_i^{r^{+}_i -r_{i,\sigma_{i,r}(j-1) }} }
    \right\rrangle = \left\llangle   \frac{ \mm_{i, \sigma_{i,r}(j)} + \varepsilon_{i}} {p_i^{r_{i,\sigma_{i,r}(j)} -r_{i,\sigma_{i,r}(j-1) }} }     \right\rrangle
    \geq
    \frac{ |\mm_{i,\sigma_{i,r}(j)}| -\sh}{p_i^{r_{i,\sigma_{i,r}(j)}-r_{i,\sigma_{i,r}(j-1)}}}
 \geq
     \frac{ |\mm_{i,\sigma_{i,r}(j)}|/2 }{p_i^{r_{i,\sigma_{i,r}(j)}-r_{i,\sigma_{i,r}(j-1)}}}.
\end{equation}
Hence, statements  \eqref{Le2-80}, \eqref{Le2-25} and Lemma 3 are proved.\qed\\  \\
Using Lemma 2 and Lemma 3, we obtain:
\begin{equation}  \label{Le5-3}
  |E_s( \sD_{T_{\fs,\sh,1}} (Q,N)) | \ll
\sum_{\substack{\br_j \in \UU_{ T_{\fs}} \\  j=1,...,\sh}} \;
\sum_{\substack{ m_j \in I^{*}_{P_{\br_j}} \\  j=1,...,\sh}} \;
\frac{ \varpi_{\br,\bm,1} }{\bar{m}_1 \cdots \bar{m}_{\sh}}
 \prod_{i=1}^s  \prod_{j=1}^{\sh} \frac{1}{\bar{\mm}_{i, \sigma_{i,r}(j)}}.
\end{equation} \\
We can apply our main tool, Theorem A, only if $m_i$  and  $\fm_{i,j} $ $(\mm_{i,j}) $ are small enough. The purpose of Lemma 4 is to show that the part of the discrepancy with large $m_i$  and  $\fm_{i,j} $ $(\mm_{i,j}) $   is negligible. So, we will consider only the case $ \max_{i,j}( |m_j|, |\fm_{i,j}|) \leq n^{4 \sh s} $.
 This will allow us to apply Theorem A  in the future.\\
\subsection{\textcolor{blue}{The main lemma}}
{\bf Lemma 4.}  {\it With the notations as above}
\begin{equation}  \label{Le5-1}
  |E_s (\sD_{T_{\fs,\sh,1}}(Q,N)) | \ll   1 + \cD_{T_{\fs}},
\end{equation}
with
\begin{equation} \label{Le5-1a}
  \cD_{T_{\fs}}=
\sum_{\substack{m_j \in I^{*}_{n^{4 \sh s}}, \; |\fm_{i,j}| \leq n^{4 \sh s}\\  i=1,...,s,\; j=1,...,\sh}} \; \sum_{\substack{\br_j \in \UU_{ T_{\fs}}\\ j=1,...,\sh}}\;
\frac{\varpi_{\br,\bm,1} \zeta_{\br,\fm} }{\bar{m}_1 \cdots \bar{m}_{\sh}}  \;
 \prod_{i=1}^s  \prod_{j=1}^{\sh} \frac{1}{\bar{\fm}_{i, j}}
       \delta_{p_i^{r^{+}_i}} (\tilde{m}_{i} )
\end{equation}
and
\begin{equation} \nonumber
  |E_{s+1} (\sD_{T_{\fs,\sh,1}}(-[Nx_{s+1}], 2[Nx_{s+1}])) | \ll   1 + \tilde{\cD}_{T_{\fs}},
\end{equation}
with
\begin{equation} \nonumber
   \tilde{\cD}_{T_{\fs}}=
\sum_{\substack{m_j \in I^{*}_{n^{4 \sh s}}, \; |\fm_{i,j}| \leq n^{4 \sh s}\\  i=1,...,s,\; j=1,...,\sh}}\;
 \sum_{\substack{\br_j \in \UU_{ T_{\fs}}\\ j=1,...,\sh}}\;
\frac{\varpi_{\br,\bm,1} \zeta_{\br,\fm}\; |\hat{\gamma}^{(\sh)}_{\br,\bm}| }{\bar{m}_1 \cdots \bar{m}_{\sh}}  \;
\;
\prod_{i=1}^s  \prod_{j=1}^{\sh} \frac{1}{\bar{\fm}_{i, j}}
       \delta_{p_i^{r^{+}_i}} (\tilde{m}_{i} ) ,
\end{equation}\\
where
\begin{equation}  \label{Le5-50}
 \tilde{m}_{i}   = \tilde{m}_{i}( \br,\bm,\fm)=
  \sum_{j=1}^{\sh}  (\fm_{i,j}  +  m_{j} M_{i, \br_j}) p_i^{ r^{+}_i -  r_{i,j}  }
\end{equation}
and
\begin{equation}  \label{Le5-0}
 \zeta_{\br,\fm} = \prod_{i=1}^s \b1( \max_{1 \leq j \leq \sh} |\fm_{i,j}| \leq  p_i^{r^{+}_i}),
\end{equation}
\begin{equation}  \nonumber
\hat{\gamma}^{(\sh)}_{\br,\bm}:= \bar{m}_1 \cdots \bar{m}_{\sh} \int_0^1 \prod_{j=1}^{\sh} \varphi_{ \br_j,-[Nx_{s+1}], 2[Nx_{s+1}],m_j} d x_{s+1}.
\end{equation}
\\
{\bf Proof.} We will prove  the first relation. The proof of the second relation is similar.

  From \eqref{Le3-10}, \eqref{Le3-20} and \eqref{Le5-50}, we derive
%
%
\begin{multline*} 
   \tilde{m}_{i}( \br,\bm,\mm)=
  \sum_{j=1}^{\sh}  (\mm_{i,j}  +  m_{j} M_{i, \br_j}) p_i^{ r^{+}_i -  r_{i,j}  }\\
= \sum_{j=1}^{\sh}  (\mm_{i,\sigma_{i,r}(j)} +  m_{\sigma_{i,r}(j)} M_{i, \br_{\sigma_{i,r}(j)}})p_i^{ r^{+}_i -  r_{\sigma_{i,r}(j)}  } \\
= \hat{m}_{i} + \sum_{j=1}^{\sh}   m_j M_{i,\br_j}
   p_i^{r^{+}_i -r_{i,j} } \equiv 0  \mod p_i^{r^{+}_i}.
\end{multline*}
%
Using \eqref{Le3-20},we get that $ \mm_{i,\sigma_{i,r}(j)} \in I_{p_i^{r_{i,\sigma_{i,r}(j) }-r_{i,\sigma_{i,r}(j-1)}} }$ are unique for given $m_{i,j} \; (j \geq 1)$. Together with \eqref{Le5-3}, this implies that
\begin{multline}  \label{L5a-4}
    |E_s (\sD_{T_{\fs,\sh,1}}(Q,N)) | \ll
\sum_{\substack{\br_j \in \UU_{ T_{\fs}} \\  j=1,...,\sh}} \;
\sum_{\substack{ m_j \in I^{*}_{P_{\br_j}} \\  j=1,...,\sh}} \;
\frac{ \varpi_{\br,\bm,1} }{\bar{m}_1 \cdots \bar{m}_{\sh}}
 \prod_{i=1}^s  \prod_{j=1}^{\sh} \frac{1}{\bar{\mm}_{i, \sigma_{i,r}(j)}}\\
    \ll \sum_{\substack{\br_j \in \UU_{ T_{\fs}} \\  j=1,...,\sh}} \;
\sum_{\substack{ m_j \in I^{*}_{P_{\br_j}} \\  j=1,...,\sh}}  \;\frac{ \varpi_{\br,\bm,1}  }{\bar{m}_1 \cdots \bar{m}_{\sh}} \;
\sum_{\substack{\fm_{i,\sigma_{i,r}(j)} \in I_{p_i^{r_{i,\sigma_{i,r}(j) }-r_{i,\sigma_{i,r}(j-1)}} }  \\  i=1,...,s,\;j=1,...,\sh}}
 \prod_{j=1}^{\sh} \prod_{i=1}^s \frac{1}{\bar{\fm}_{i,j}}
   \delta_{p_i^{r^{+}_i}} (  \tilde{m}_{i}( \br,\bm,\fm) ).
\end{multline}
%
%
From \eqref{Beg-1a}, we obtain $\# I^{*}_{P_{\br}} \leq P_{\br} \leq (p_1 \cdots p_s)^n$.
Taking into account that the variables $\fm_{i,j} \; (i,j \geq 1)$ are uniquely determined from the variables $m_{i,j} \; (i,j \geq 1)$  by the conditions
 $ \delta_{p_i^{r^{+}_i}} (   \tilde{m}_{i}( \br,\bm,\fm) ) =1$, we get that
$   \cD_{T_{\fs}}  \ll n^{\sh(s+1)}$ (see \eqref{Le5-1}, \eqref{Le5-1a} and \eqref{L5a-4} ). We have that the part of the right hand of \eqref{L5a-4},  satisfying to the condition $|\fm_{i,j}| >n^{4 \sh s}$ for some
$(i,j) \in [1,s] \times [1,\sh]$, is equal to $O(n^{\sh(s+1)-4\sh s})$.
Therefore, we can only consider the case $ |\fm_{i,j}| \leq n^{4 \sh s}$:
\begin{multline*}  
    |E_s (\sD_{T_{\fs,\sh,1}}(Q,N)) | \ll n^{-s} +	 \sum_{\substack{\br_j \in \UU_{ T_{\fs}} \\ j=1,...,\sh}} \;\;
\sum_{\substack{ m_j \in I^{*}_{P_{\br_j}} \\  j=1,...,\sh}} \;\;
\sum_{\substack{|\fm_{i,j}| \leq n^{4 \sh s} \\  i=1,...,s,\;j=1,...,\sh}} \;
\frac{ \varpi_{\br,\bm,1} \zeta_{\br,\fm} }{\bar{m}_1 \cdots \bar{m}_{\sh}}   \\
 \times   \prod_{j=1}^{\sh} \prod_{i=1}^s \frac{1}{\bar{\fm}_{i,j}}
   \delta_{p_i^{r^{+}_i}} (   \tilde{m}_i )  .
\end{multline*}
%
Note that we will use the variable $\zeta_{\br,\fm}$ only in Lemma 18.

Let
\begin{equation*} 
\jmath_1(\bm): =\b1(\max_j |m_j| < (n^{4 \sh s}-1)/2) \qquad \ad \quad \jmath_2(\bm): =1-\jmath_1(\bm).
\end{equation*}
We have
\begin{multline*}  
   |E_s( \sD_{T_{\fs,\sh,1}} (Q,N)) | \ll n^{-s} +	\ddot{\cD}_{T_{\fs,1}}+ \ddot{\cD}_{T_{\fs,2}}, \;\;\;  \with  \\
 \ddot{\cD}_{T_{\fs},\nu}=
\sum_{\substack{\br_j \in \UU_{ T_{\fs}}  \\ j=1,...,\sh}} \;\;
\sum_{\substack{|\fm_{i,j}|  <n^{4 \sh s} \\ i=1,...,s,\; j=1,...,\sh}} \;\; \sum_{\substack{m_j \in I^{*}_{P_{\br_j}} \\  j=1,...,\sh}}\;
   \prod_{j=1}^{\sh} \frac{ \varpi_{\br,\bm,1} \zeta_{\br,\fm}  }{\bar{m}_j}
      \prod_{i =1 }^s \frac{ \jmath_{\nu}(\bm)}{\bar{\fm}_{i,j}}  \delta_{p_i^{r^{+}_i}} ( \tilde{m}_i ) , \; \nu=1,2.
\end{multline*}	
From \eqref{Le5-1}, we obtain that $ \ddot{\cD}_{T_{\fs},1} \leq \cD_{T_{\fs}} $, and the assertion of Lemma 4 follows from the estimate:
\begin{equation} \label{L5a-7}
     \ddot{\cD}_{T_{\fs},2} \ll 1.
\end{equation}
Now we consider $\ddot{\cD}_{T_{\fs,2}}$  (the case of $\jmath_{2}(\bm) =1 $ or $\max_j |m_j| \geq (n^{4 \sh s}-1)/2$).
Let  $i_0 \in T_{\fs}$. By  \eqref{Beg-22} and \eqref{Beg-26}, we get that $ r_{i_0,j} > V_1= [\log_2^3 n]$ $(j=1,...,\sh)$.

Let $|m_{j_0}| \geq (n^{4 \sh s}-1)/2 $  with some $j_0 \in [1,h]$.
 Let $j_0=1$. The case $j_0 \in [2,\sh] $ is similar.
It is easy to verify that \eqref{L5a-7} ensue from the following inequality
\begin{equation} \label{Le5-7}
W=W(m_2,...,m_{\sh},\br,\fm):  =  \sum_{p_0^n \geq |m_{1}| \geq (n^{4 \sh s}-1)/2} \frac{ \delta_{p_{i_0}^{r^{+}_{i_0}}}(\tilde{m}_{i_0} ) }{\bar{m}_{1}} \ll n^{-2\sh s}.
\end{equation}
Now we will prove \eqref{Le5-7}.
We fix $\fm_{i,j}\; (i \in [1,s], j \in [1, \sh])$, $\br_1,...,\br_{\sh}$ and $ m_2,...,m_{\sh}$.
Let
\begin{equation} \nonumber
K_1=K_1(m_2,...,m_{\sh},\br,\fm):=   -\fm_{{i_0},1}p_{i_0}^{ r^{+}_{i_0} -   r_{i_0,1 }}
-\sum_{j \in [2,\sh]}   ( \fm_{{i_0},j}
 + m_{j} M_{i, \br_j} )
  p_{i_0}^{ r^{+}_{i_0} -   r_{{i_0},j }}.
\end{equation}
We have
\begin{multline} \label{Le5-77}
W=W_1+W_2,   \quad \with  \quad W_{\nu}=  \sum_{p_0^n \geq |m_{1}| \geq (n^{4 \sh s}-1)/2} \frac{1}{\bar{m}_{1}}  \delta_{p_{i_0}^{r^{+}_{i_0}}}
  ( \tilde{m}_{{i_0}}  ) \tilde{\theta}_{\nu}(K_1), \quad \where   \\
  \tilde{\theta}_{1}(K_1)=\delta_{ p_{i_0}^{ r^{+}_{i_0} -   r_{{i_0},1}} }
  ( K_1 ),
    \;\;\;\; \tilde{\theta}_{2}(K_1) =1-\tilde{\theta}_{1}(K_1).
\end{multline}
Bearing in mind that $ \tilde{m}_{i_0} \equiv 0 \mod p_{i_0}^{ r^{+}_{i_0}} $, we get from   \eqref{Le5-50} that
\begin{equation} \nonumber
\sum_{j=1}^{\sh}   ( \fm_{{i_0},j} + m_{j} M_{i_0, \br_j} ) p_{i_0}^{ r^{+}_{i_0} -   r_{{i_0}, j}  }
\equiv 0 \mod p_{i_0}^{ r^{+}_{i_0}}.
\end{equation}
Hence
\begin{equation} \nonumber
 ( \fm_{{i_0},1 }
  + m_{1} M_{i_0, \br_{1}} ) p_{i_0}^{ r^{+}_{i_0} -   r_{i_0,1} }
  \equiv -\sum_{j \in [2,\sh]}    ( \fm_{{i_0},j} + m_{j}  M_{i_0, \br_{j }} )
 p_{i_0}^{ r^{+}_{i_0} -   r_{{i_0},j}  }
\mod p_{i_0}^{ r^{+}_{i_0}}
\end{equation}
and
\begin{equation*}
  m_{1}  M_{i_0, \br_{ i_0}, 1}  p_{i_0}^{ r^{+}_{i_0} -   r_{{i_0},1}  }
     \equiv
 -\fm_{{i_0},1 } p_{i_0}^{ r^{+}_{i_0} -   r_{{i_0},1 }  }  -\sum_{j \in [2,\sh]}   ( \fm_{{i_0},j } + m_{j}  M_{i_0, \br_{j}} )
 p_{i_0}^{ r^{+}_{i_0} -   r_{{i_0}, j}  }
\mod p_{i_0}^{ r^{+}_{i_0}}.
\end{equation*}
Therefore
\begin{equation} \label{Le5-76}
 m_{1} M_{i_0, \br_{ i_0 , 1 }} p_{i_0}^{ r^{+}_{i_0} -   r_{{i_0}, 1 }  }
       \equiv K_1 \mod p_{i_0}^{ r^{+}_{i_0}}.
\end{equation}
In view of \eqref{Le5-76} and \eqref{Le5-77}, we have
\begin{equation} \label{Le5-78}
 K_1 \equiv 0 \mod  p_{i_0}^{ r^{+}_{i_0}- r_{{i_0},1}} ,
 \quad \tilde{\theta}_{1}(K_1) =1 \qquad \ad \qquad W_2=0.
\end{equation}
Now let's consider $W_1$.  
 Let $K_2 = K_1    p_{i_0}^{-r^{+}_{i_0} +r_{{i_0}, 1} }$.
%
Then  $ m_{1} M_{i_0, \br_{1}} \equiv K_2 \mod p_{i_0}^{ r_{i_0,1}} $.
According to \eqref{Beg-3}, we get $ M_{i_0, \br_{1}} P_{\br_1}/  p_{i_0}^{r_{i_0,1}}
  \equiv 1 \mod p_{i_0}^{r_{i_0,1}}$.\\
Let $ K_3= K_2 P_{\br_1}/  p_{i_0}^{r_{i_0,1}}$.
Therefore
\begin{equation*} 
 m_{1} \equiv K_3    \mod p_{i_0}^{r_{i_0,1}} , \quad \with  \quad
 r_{i_0,1} \geq V_1 \geq \log_2^3 n .
\end{equation*}
By \eqref{Le5-77}, we have
\begin{multline*}  
W_1 \leq    \sum_{p_0^n \geq |m_{1}| \geq (n^{4 \sh s}-1)/2} \frac{1}{\bar{m}_{1}} \b1(  m_{1} \equiv K_3    \mod p_{i_0}^{r_{i_0,1}})\\
\ll \sum_{k \in \ZZ} \frac{1}{|K_3 +kp_{i_0}^{r_{i_0,1}}|}
 \b1 \big( p_0^{s n}  \geq |K_3 +kp_{i_0}^{r_{i_0,1}}|  \geq (n^{4 \sh s}-1)/2 \big)  \\
 \ll \sum_{k=0}^{p_0^{sn}} \frac{1}{n^{4 \sh s}+k n^{5\sh s}}
  \ll n^{-4 \sh s}+ n^{-5 \sh s} \log (p_0^{sn}) \ll  n^{-4\sh s}.
\end{multline*}
Taking into account \eqref{Le5-78} and \eqref{Le5-77}, we obtain that   \eqref{Le5-7}. So,  Lemma 4 is proved. \qed \\

We see that in order to prove \eqref{Le3-3} it is enough to verify that
\begin{equation}  \label{Le3-4}
\cD_{T_{\fs}}\ll  n^{\sh s/2 }.
\end{equation}\\ \\
%
%

\subsection{\textcolor{blue}{Linear forms in $\bbp$-adic logarithm}}
Let $\ddot{\alpha}_1,...,\ddot{\alpha}_{\ddot{n}}$ $(\ddot{n} \geq 1)$ be non-zero algebraic numbers and $K$ be a number
field containing $\ddot{\alpha}_1,...,\ddot{\alpha}_{\ddot{n}}$ with $d=[K: \QQ]$. Denote by $\fd$  a prime ideal
of the ring $\fO_K$ of integers in $K$, lying above the prime number $\bbp$, and by $\ff_{\fd}$  the residue class degree of $\fd$. For
$\gamma \in K$, $\gamma \neq 0$, write ${\rm ord}_{\fd}(\gamma)$ for the exponent to which $\fd$ divides the principal fractional ideal generated by $\gamma $ in $ K$.
Define
\begin{equation*}  
    h^{'}(\ddot{\alpha}_j) = \max( h_0(\ddot{\alpha}_j), \ff_{\fd} (\log \bbp)/d)
		\qquad  (1 \leq j \leq \ddot{n}),
\end{equation*}
where $h_0(\gamma)$ denotes the absolute logarithmic Weil height of an algebraic
number $\gamma$, i.e.,
\begin{equation*}  
    h_0(\gamma) = k^{-1} \Big(  \log a_0  + \sum_{i=1}^k \log \max(1, |\gamma^{(i)}|) \Big) ,
\end{equation*}
where the minimal polynomial for $\gamma$ is
\begin{equation*}  
  a_0x^k+   a_1 x^{k-1} + \cdots + a_k =a_0(x-\gamma^{(1)}) \cdots (x-\gamma^{(k)}), \qquad a_0 >0.
\end{equation*}

{\bf Theorem A.}  (\cite[Theorem 1]{Yu}, \cite[Theorem 2.9]{Bu}) {\it Let $\ddot{\Xi} = \ddot{\alpha}_1^{b_1} \cdots \ddot{\alpha}_{\ddot{n}}^{b_{\ddot{n}}} -1 \neq 0 $. Suppose that }
\begin{equation*}  
    {\rm ord}_{\fd}(\ddot{\alpha}_1) =0 \qquad (1 \leq j \leq \ddot{n}).
\end{equation*}
Then there exists a constant $C$, depending only on $\ddot{n}, d$ and $\fd$, such that
\begin{equation*}  
    {\rm ord}_{\fd}(\ddot{\Xi})  < C h^{'}(\ddot{\alpha}_1) \cdots h^{'}(\ddot{\alpha}_{\ddot{n}}) \log \ddot{B},
\end{equation*}
where
\begin{equation*}  
\ddot{B} = \max(|b_1|,...,|b_{\ddot{n}}| ,3).
\end{equation*}

We will use Theorem A with $\ddot{n}=s,\;k= d=1,\; \fd=p_i$,
$\{\ddot{\alpha}_1,...,\ddot{\alpha}_{s-1} \} =\{p_1,p_2,...,p_{i-1},p_{i+1},...,p_s\} $, $\ddot{\alpha}_s =k_1/k_2$, $k_1,k_2 \in \ZZ$.\\

{\bf Corollary 1.}  {\it Let $0<|k_j| \leq n^{{4 \sh s}},\; j=1,2$,  $  {\rm ord}_{p_i}(k_1/k_2)=0$
 ($i \in [1,s]$) and $\ddot{B} =n$. Then there exists a constant $C_1>0$ such that }
\begin{equation*}  
   {\rm ord}_{p_i}\Big( \ddot{\Xi} \Big)  < C_1 \log n \log \ddot{B} =C_1 \log^2 n, \;\;\;\; \with \;\;\;\; \ddot{\Xi}=(k_1 / k_2) \prod_{1 \leq j \leq s, j \neq i} p_j^{b_j} -1 \neq 0.
\end{equation*}\\ \\
%
\section{\textcolor{blue}{Applications of the $\bbp$-adic logarithm form theorem}}

\subsection{\textcolor{blue}{Notations}}
Basically, all we need to do is calculate the number of
 solutions of the
congruence $\tilde{m}_{i} \equiv 0 \; \mod  p_i^{r^{+}_i} $ (see Lemma 4).
By Theorem A, it is very easy to compute
this number for the case $\min |r_{i_1,j_1} - r_{i_1,j_2}| \geq \log^3_2 n$.
 The case $\min |r_{i_1,j_1} - r_{i_1,j_2}| < \log^3_2 n$ is much more complicated.
Therefore, we need to divide the interval $[1,n]$ into large and small
 intervals by the points $r_{i,j}$: \\

 First, we order sequences $(r_{i,j})_{j = 1}^{\sh}$ in the increasing order. We do this using a permutation of $\{ 1,2,...,\sh\}$.  Next, we select values $r_{i,j}, \; j \geq 1$ that are close to each other for $i \in [1,...,s]$.
  We introduce for this purpose sequences $b^{'}_{i, \br, 1}  $ and $\tilde{a}_{i,\br_{i}} $.

Let  $\Xi_{\sh}$ be the set of all permutations of  $\{ 1,2,...,\sh\}$ and
let $\sigma_{i,r}  \in \Xi_{\sh}$ with
$r_{i,\sigma_{i,r}(j+1) }  \geq r_{i,\sigma_{i,r}(j) } $ for $j =1,2,...,\sh-1$.
For a given nondecreasing sequence $(r_{i,\sigma_{i,r}(j)})_{1 \leq j \leq \sh}$, we define the following partition of the interval $[0,\sh]$: $0=b^{'}_{i, \br, 0} < b^{'}_{i, \br, 1} < \cdots < b^{'}_{i, \br, \tilde{a}_{i,\br_{i}}}=\sh $.
More precisely,  we define integer variables $\tilde{a}_{i,\br_{i}}$,
$b^{'}_{i, \br, k}$, $ b_{i, \br, k} \geq 1 $     from  the following conditions:
\begin{equation*}   
     b^{'}_{i, \br, k}   = b^{'}_{i, \br, k-1} +  b_{i, \br, k} \;\; \;\;  1 \leq k \leq \tilde{a}_{i,\br_{i}},\; b^{'}_{i, \br,0} =0, \; b^{'}_{i, \br,\tilde{a}_{i,\br_{i}}  } =\sh,\;\;\;
     b_{i, \br, 1} + \cdots + b_{i, \br, \tilde{a}_{i,\br_{i}} } =\sh,
\end{equation*}
with
\begin{multline}  \label{Par3-1}
 0 \leq	r_{i,\sigma_{i,r}(j_1) }  - r_{i,\sigma_{i,r}(j_2) }	\leq V_1 \;\;\;\;
	{\rm if} \;\;\; \;  j_1,j_2 \in (b^{'}_{i, \br, k-1},b^{'}_{i, \br, k}  ],
\;\;j_1>j_2,\;\; k  \in [1, \tilde{a}_{i,\br_{i}}], \\
r_{i,\sigma_{i,r}(b^{'}_{i, \br, k+1}) }  - r_{i,\sigma_{i,r}(b^{'}_{i, \br, k+1}-1)}
	>  V_1= [\log_2^3 n]\;\;\;
 \for \;\;\;
k  \in [1, \tilde{a}_{i,\br_{i}}-1], \;\; i=1,...,s.
\end{multline}
Now we partition $\cD_{T_{\fs}}$ (see \eqref{Le5-1a}, \eqref{Le5-50} and \eqref{Le5-0}) according to condition \eqref{Par3-1}:
\begin{multline}  \nonumber 
 \cD_{T_{\fs}}  =
\sum_{ m_j \in I^{*}_{n^{4 \sh s}},\; |\fm_{i,j}| \leq n^{4 \sh s},   i=1,...,s,j=1,...,\sh} \;\;\;
\sum_{\br_{j} \in \UU_{T_{\fs}},    j=1,...,\sh}  \; \prod_{i=1}^s
    \delta_{p_i^{r^{+}_i}} ( \tilde{m}_{i} )  \\
   \times
   \prod_{j=1}^{\sh} \frac{\varpi_{\br,\bm,1}  \zeta_{\br,\bm} }{\bar{m}_j} \;
    \prod_{i =1}^s  \frac{1}{\bar{\fm}_{i,j}} \sum_{a_1,...,a_s =1}^{\sh} \;
\; \;    \sum_{ 1 \leq \lambda_{i,k} \leq \sh, i=1,...,s ,\; \lambda_{i,1}+...+\lambda_{i,a_i}=\sh }\;\;\;
 \sum_{ \tau_i \in \Xi_{\sh},  i=1,...,s}\; 1\\
   \times \prod_{i =1}^s  \b1( \tilde{a}_{i,\br_{i}}=a_i ) \b1(b_{i,\br_{i} ,k}=\lambda_{i,k} )
  \b1( \sigma_{i,r}=\tau_i), \\
\where \quad \tilde{m}_{i} =
  \sum_{j=1}^{\sh}  (\fm_{i,j}  +  m_{j} M_{i, \br_j}) p_i^{ r^{+}_i -  r_{i,j}  }
   \quad \; \ad \quad \; \zeta_{\br,\fm} = \prod_{i=1}^s \b1( \max_{1 \leq j \leq \sh} |\fm_{i,j}| \leq  p_i^{r^{+}_i}).
\end{multline}
Changing the order of the summation, we obtain
\begin{equation} \label{Par3.2-3}
 \cD_{T_{\fs}}  =  \sum_{a_1,...,a_s =1}^{\sh} \;\;
    \sum_{1 \leq \lambda_{i,k} \leq \sh, i=1,...s , \; \lambda_{i,1}+...+\lambda_{i,a_i}=\sh }\;\;
 \sum_{ \tau_i \in \Xi_{\sh},  i=1,...,s}  \DD_{T_{\fs}, a, \lambda,\tau}
\end{equation}
with
\begin{equation} \label{Par3.2-4}
\DD_{T_{\fs}, a, \lambda,\tau}   = \sum_{ m^{'}_j \in I^{*}_{n^{4 \sh s}},\; |\fm^{'}_{i,j}| \leq n^{4 \sh s},  \; i=1,...,s,\; j=1,...,\sh} \;\;
 \prod_{j=1}^{\sh} \frac{1}{\bar{m}^{'}_j}
     \prod_{i =1}^s\frac{1}{\bar{\fm}^{'}_{i,j}}
      \SS_{m^{'},  \tau},
\end{equation}
where
\begin{multline} \nonumber
  \SS_{m^{'}, \tau} =\sum_{\br_{j} \in \UU_{T_{\fs}},    j=1,...,\sh} \;\;
    \prod_{i =1}^s \delta_{p_i^{r^{+}_i}} ( \tilde{m}_i ) \varpi_{\br,\bm^{'},1}  \zeta_{\br,\fm^{'}}   \\
    \times \prod_{i=1}^s  \b1( \tilde{a}_{i,\br_{i}}=a_i ) \b1(b_{i,\br_{i} ,k}=\lambda_{i,k} )
  \b1( \sigma_{i,r}=\tau_i)
\end{multline}
and
\begin{equation}  \label{Par3.2-40}
 \tilde{m}_{i} =
  \sum_{j=1}^{\sh}  (\fm^{'}_{i,j}  +  m^{'}_{j} M_{i, \br_j}) p_i^{ r^{+}_i -  r_{i,j}  }.
\end{equation}
By \eqref{Par3-1}, 
 we get
\begin{multline} \nonumber   
\SS_{m^{'}, \tau} \leq \sum_{ \br_{j} \in \UU_{T_{\fs}}, \;   j=1,...,\sh}  \;\;
    \prod_{i\in  T_{\fs}}  \delta_{p_i^{r^{+}_i}} (  \tilde{m}_i ) \varpi_{\br,\bm^{'},1}  \zeta_{\br,\fm^{'}}  \\
  \times   \b1\Big( 0 \leq 	r_{i,\tau_{i}(j_1) }  - 	r_{i,\tau_{i}(j_2) }	 \leq V_1 \;\;\;
	{\rm if} \;\;\;  j_1,j_2 \in (\Lambda_{i,k-1},\Lambda_{i,k}],\;\;j_1>j_2,\;\; \for \;\;  k  \in [1, a_{i}] \\
 \& \quad  r_{i,\tau_{i}(\Lambda_{i,k+1}) } - r_{i,\tau_{i}(\Lambda_{i,k+1} -1) }	> V_1 \;\;
 \for \;
k  \in [1, a_i-1], \; i\in  T_{\fs} \Big),
\end{multline}
where
\begin{equation}   \label{Par3-5}
   \Lambda_{i,  k}   = \Lambda_{i,  k-1} +  \lambda_{i,  k},\; \;\;1 \leq k \leq a_i,\;
   \Lambda_{i,  0} =0, \; \Lambda_{i,  a_i} =\sh, \;\; \sum_{k=1}^{a_i} \lambda_{i,  k} =\sh,
\end{equation}
By \eqref{Beg-22} and \eqref{Beg-26}, we get
\begin{equation}   \label{Par3-5a}
\SS_{m^{'}, \tau} \leq
   V_1^{\sh(s-\fs)} \; \prod_{ i\in  T_{\fs}} n^{a_i} V_1^{\sh}=  n^{a_1+ \cdots +a_{\fs}} V_1^{\sh s}.
\end{equation}
Let
\begin{equation}   \label{Par3-7}
 a_{i_0} = \max_{i \in T_{\fs}} a_i \quad \for \; {\rm some} \;\; i_0 \in T_{\fs}.
\end{equation}
In the following we fix $i_0 \in T_{\fs}$ and $\tau_{i_0} \in \Xi_{\sh}$.
Let's consider the case $i_0=1$. The proof for the case $i_0 \neq 1$ is similar.
Let
\begin{equation}  \label{Par3-8}
\rho_{i,j}= r_{i,\tau_{1}(j)}, \;\; \brho_j =\br_{\tau_{1}(j)}, \;\; m_j = m^{'}_{\tau_{1}(j)},
\;\; \fm_{i,j} = \fm^{'}_{i,\tau_{1}(j)}, \;\;
        \dot{\tau}_{i}(j) = \tau_{1}^{-1}(j)(\tau_{i}(j))
\end{equation}
$ i \in [1,s], \;j \in [1,\sh]$. By \eqref{Lem2-6}, we get
\begin{equation}   \nonumber
    \varpi_{\brho,\bm,j}=\varpi_{\br,\bm^{'},j} \quad  \with \quad \brho=(\brho_1,...,\brho_{\sh}), \;\; \br=(\br_1,...,\br_{\sh}), \;\; j=1,2.
\end{equation}

In view of  \eqref{Lem2-6}, \eqref{Le3-17}, \eqref{Par3.2-4}, \eqref{Par3.2-40}, \eqref{Par3-5} and \eqref{Par3-8}, we proved the following lemma: \\ \\
{\bf Lemma 5.} {\it We have}
\begin{equation} \label{Par3.2-45}
\DD_{T_{\fs}, a, \lambda,\tau}   = \sum_{m_j \in I^{*}_{n^{4 \sh s}},\; |\fm_{i,j}| \leq n^{4 \sh s},   i=1,...,s, j=1,...,\sh } \;\;
 \prod_{j=1}^{\sh} \frac{1}{\bar{m}_j} \;
     \prod_{i =1}^s\frac{1}{\bar{\fm}_{i,j}}
      \SS_{m,  \tau}
\end{equation}
with
\begin{equation}  \label{Le6-1a}
           \SS_{m, \tau} \leq       \sum_{\brho_{j} \in \UU_{T_{\fs}}, j\in[1,\sh]} \;   \varpi_{\brho,\bm,1} \; \zeta_{\brho,\fm} \;
    \chi_{\brho} \;   \delta_{p_1^{\rho^{+}_{1}}}(\tilde{m}_1),
\end{equation}
where
\begin{equation}  \label{Le6-2a}
\chi_{\brho}  =\prod_{i \in T_{\fs}}  \chi_{i,\brho}, \qquad \chi_{i,\brho} =\prod_{k=1}^{a_i} \chi_{i,k,\brho}, \qquad
  \chi_{i,k,\brho} =\prod_{ j \in [\Lambda_{i,k-1}+1,\Lambda_{i,k}] }  \chi_{i,k,j,\brho},
\end{equation}
\begin{equation} \nonumber
\chi_{i,k,j,\brho}=\b1\big( 0 \leq  \rho_{i,\dot{\tau}_{i}(j+1) }  - 	\rho_{i,\dot{\tau}_{i}(j)} 	\leq V_1\big) \;\;\;
	 \;\;\for \;  j, j+1 \in (\Lambda_{i,k-1},\Lambda_{i,k}],\;\; k \in [1,a_i],
\end{equation}
\begin{equation} \nonumber
\chi_{i,k,j,\brho}=\b1\big(      \rho_{i,\dot{\tau}_{i}(\Lambda_{i,k+1}) } - \rho_{i, \dot{\tau}_{i}(\Lambda_{i,k+1}-1) }	> V_1 \Big),\;\; \for \;\; j=\Lambda_{i,k}, \; k <a_i, \qquad \qquad \qquad \qquad \qquad \qquad \qquad
\end{equation}
\begin{equation} \nonumber
%
\chi_{i,a_i,\sh,\brho}=\b1\big(  0 \leq \rho_{i,\dot{\tau}_{i}(\sh) }  - 	\rho_{i,\dot{\tau}_{i}(\sh-1) }	\leq V_1\big) \;\;\;
	 \;\;\for \;  \sh-1 > \Lambda_{i,a_i-1}, \qquad \qquad \qquad \qquad \qquad \qquad\qquad
\end{equation}
\begin{equation} \nonumber
\chi_{i,a_i,\sh,\brho}=\b1\big(  \rho_{i,\dot{\tau}_{i}(\sh) }  - 	\rho_{i,\dot{\tau}_{i}(\sh-1) }	> V_1\big) \;\;\;
	 \;\;\for \;  \sh-1 = \Lambda_{i,a_i-1}, \qquad \qquad \qquad \qquad \qquad \qquad\qquad \qquad
\end{equation}
\begin{equation} \nonumber
\;\;\rho_{1,0}=0,\; \rho_{1,\sh} = \rho^{+}_{1}=\max_j \rho_{1,j}, \;     \rho_{1,j} \geq \rho_{1,j-1} \; (j \geq 1)   \qquad \qquad \qquad \qquad \qquad \qquad\qquad \qquad
\end{equation}
and
\begin{equation}  \label{Par3.2-41}
 \tilde{m}_{1} =
  \sum_{j=1}^{\sh}  (\fm_{1,j}  +  m_{j} M_{1, \brho_j}) p_1^{ \rho^{+}_1 -  \rho_{1,j}  }, \qquad
     \zeta_{\brho,\fm} = \prod_{i=1}^s \b1( \max_{1 \leq j \leq \sh} |\fm_{i,j}| \leq  p_i^{\rho^{+}_i}),
\end{equation}
\begin{multline}  \label{Lem2-6d}
\varpi_{\brho,\bm,1} = 1- \varpi_{\brho,\bm,2}, \quad \with \qquad
 \varpi_{\brho,\bm,2} = \b1 \big(\sh=2q, \; \fs=s,\; \exists \eta \in \Xi_{2q}  \; : \\
 m_{\eta(2k-1)}/P_{\brho_{\eta(2k-1)}}  =- m_{\eta(2k)}/P_{\brho_{\eta(2k)}},
  k \in [1,q]  \big) .
\end{multline} \\ \\
{\bf Corollary 2.}  {\it Let $  \ff_k: = \Lambda_{1,k}$, $i_1 \in [2,s]$,$k_0 \in \{1,...,a_{i_1}\}, \;\ff_{0}:=0, \; \ff_{a_{1}+1} :=\ff_{a_{1}}=\sh$,
$j  \in [\Lambda_{i,k-1}+1,\Lambda_{i,k}]  $. Then}

\begin{equation}  \label{Le6-1}
           \SS_{m, \tau} \leq       \sum_{\brho_{j} \in \UU_{T_{\fs}}, j\in[1,\sh]} \;   \varpi_{\brho,\bm,1} \;\zeta_{\brho,\fm} \;
    \chi_{i_1,k_0,j,\brho}\; \prod_{k=1}^{a_1} \chi_{1,k,\brho}    \; \delta_{p_1^{\rho^{+}_{1}}}(\tilde{m}_{1} ),
\end{equation}
where
\begin{align}  \label{Le6-2}
\chi_{1,k,j,\brho}&=\b1\big( 0 \leq \rho_{1,j+1 }  - 	 \rho_{1,j }	\leq V_1 \big) \quad
	\for \quad  j, j+1 \in (\ff_{k-1},\ff_{k}],\\
\chi_{1,k,j,\brho}&=\b1\big(\rho_{1,\ff_{k+1} } - \rho_{1, \ff_{k+1} -1 }	> V_1 \big), \quad
\for\; j = \ff_{k}, \;\; k <a_1, \nonumber \\
\chi_{1,a_1,\sh,\brho}&=\b1\big(  0 \leq \rho_{1,\sh }  - 	 \rho_{1,\sh-1 }	\leq V_1 \big) \quad
	\for \quad  \sh-1>\ff_{a_1-1},  \nonumber\\
\chi_{1,a_1,\sh,\brho}&=\b1\big(  \rho_{1,\sh }  - 	 \rho_{1,\sh-1 }	> V_1\big) \quad
	\for \quad  \sh-1=\ff_{a_1-1}  \nonumber.
\end{align}

In view of \eqref{Par3.2-3}, we get from \eqref{Par3.2-45} that
\begin{equation}  \label{Le3-5}
 \cD_{T_{\fs}}  \ll \max_{m_j \in I^{*}_{n^{4 \sh s}},\; |\fm_{i,j}| \leq n^{4 \sh s},   i=1,...,s, j=1,...,\sh } \;|\SS_{m, \tau}| \; (\log n)^{\sh(s+1)}.
\end{equation} \\

The main idea of the proof of this section is follows:\\
To calculate the number of solution of the system of congruences
 $ \delta_{p_i^{r^{+}_i}} ( \tilde{m}_{i} )=1$ ($i=1,...,\fs$),
   we fixed $\fs-$dimensional partitions
   $ \prod_{1 \leq i \leq \fs} [\Lambda_{i,k_i-1}+1,\Lambda_{i,k_i}]$,\\ $ k_i \geq 1 $. In this section, we show that using Theorem A, we can examine only $1-$dimensional partitions $ [\ff_{k-1}+1,\ff_{k}] , \; k \geq 1$.   This is the most fundamental point in the proof.\\

We apply Theorem A in Lemma 6, Lemma 8, Lemma 12 and Lemma 15. Lemma 7, Lemma 9 -- Lemma 11, Lemma 13 and Lemma 14 are auxiliary.
In Lemma 9, we prove that the calculation of  $\SS_{m, \tau} $ can be reduced to considering the number of solutions only on the interval  $[\ff_{k-1}+1,\ff_{k}]$.
 We will prove \eqref{Le3-4} or \eqref{Le3-5}    separately  for the cases $\fs=1$, $\fs  \geq 2$, $\min \lambda_{i,j}=1$,
  and $\min \lambda_{i,j} =\max \lambda_{i,j}=2$.

Theorem A allows us to move from congruences to equalities.
We described a method for applying Theorem A to Halton's sequences in \cite{Le5}. Here we use the same approach. But in this article, the description is much more technically complicated. \\ \\

\subsection{\textcolor{blue}{Case $\fs=1$, $s \geq 2$}}
%
{\bf Lemma 6.} {\it Let $\fs =1$, $s \geq 2$. Then }\\ \\
\begin{equation}  \label{Le3-40}
   E_s (\sD_{T_{\fs,\sh,1}}(Q,N))\ll   n^{\sh s/2 }.
\end{equation} \\
{\bf Proof.}  In view of Lemma 4 and \eqref{Par3.2-3}, in order to prove \eqref{Le3-40} it is enough to verify that
\begin{equation}  \label{Le3-401}
\DD_{T_{\fs}, a, \lambda,\tau}  \ll    n^{\sh s/2 }.
\end{equation}

Let  $ s =3$. We obtain $a_1 +\cdots +a_{\fs} =a_1 \leq \sh <\sh s/2$. Now \eqref{Le3-4} and  \eqref{Le3-5} follow from  \eqref{Par3-5a}. The similar is true for the case  $ s =2$ and  $   a_1 < \sh$. \\

 Now let  $ s =2$, $   a_1 = \sh$.
 Hence $\lambda_{1,k}=1$ and $\ff_k=k$ for all $k \in [0,\sh]$.\\
  We define $ m_{\sh+1} =\fm_{1,\sh+1} =\rho_{2,\sh+1} = 0$.
Let's consider the following equations:
\begin{equation}  \label{Le14-52}
 m_{j}+\fm_{1,j} p_2^{\rho_{2,j}} =0, \quad \for \quad j >k_0.
\end{equation}
We will prove \eqref{Le14-52} by induction. We see that \eqref{Le14-52} is true for $k_0=\sh$.    Suppose that \eqref{Le14-52} is true for $k_0= k$. We will prove that \eqref{Le14-52} is true for $k_0=k-1$, with  $k=\sh,\sh-1,...,1$ .
According to \eqref{Beg-3},
\begin{equation}  \nonumber
 M_{1,\brho_j} P_{\brho_j}/ p_1^{\rho_{1,j}}
	\equiv 1 \mod p_1^{\rho_{1,j}} \quad \ad \quad  M_{1,\brho_j}
	\equiv p_2^{-\rho_{2,j}} \mod p_1^{\rho_{1,j}} \quad \for \;\; s=2.
\end{equation}
By \eqref{Le14-52}, we have for $j > k_0$ that $m_{j}+\fm_{1,j} p_2^{\rho_{2,j}} =0$,
\begin{equation} \label{Le14-4a}
 \fm_{1,j}  +  m_{j} M_{1, \brho_j}
     \equiv 0 \mod p_1^{\rho_{1,j}} \quad \ad \quad
 (\fm_{1,j}  +  m_{j} M_{1, \brho_j}) p_1^{ \rho^{+}_1 -  \rho_{1,j}  }
     \equiv 0 \mod p_1^{\rho^{+}_{1}}.
\end{equation}
Bearing in mind that $   \tilde{m}_{1}\equiv 0 \mod p_1^{\rho^{+}_{1}}  $ (see \eqref{Le6-1}), we obtain from
\eqref{Par3.2-41} and \eqref{Le14-4a} that
\begin{equation}  \nonumber
  \sum_{j=1}^{k}  (\fm_{1,j}  +  m_{j} M_{1, \brho_j}) p_1^{ \rho^{+}_1 -  \rho_{1,j}  }
     \equiv 0 \mod p_1^{\rho^{+}_{1}} \quad \for \quad k=k_0.
\end{equation}
Taking into account that  $ \lambda_{1,j}=1$ and $ \rho_{1,j} - \rho_{1,j-1} \geq V_1 $ $(j=1, \ldots,\sh)$, we get
\begin{multline}  \label{Le14-4}
 (\fm_{1,k}  +  m_{k} M_{1, \brho_k}) p_1^{ \rho^{+}_1 -  \rho_{1,k}  } +v_0
 p_1^{ \rho^{+}_1 -  \rho_{1,k-1}  }
     \equiv 0 \mod p_1^{\rho^{+}_{1}} \quad \with \; {\rm some}\; v_0 \in \ZZ \; \quad \ad\\
(\fm_{1,k} + m_{k} M_{1, \brho_{k}}) p_1^{ \rho^{+}_1 -  \rho_{1,k}  }   \equiv 0 \mod
 p_1^{ \rho^{+}_1 -  \rho_{1,k-1}  },  \;\;  \fm_{1,k} + m_{k} M_{1, \brho_{k}}  \equiv 0 \mod
 p_1^{ \rho_{1,k} -  \rho_{1,k-1}  } \qquad \qquad\\
\fm_{1,k} + m_{k} M_{1, \brho_{k}}   \equiv 0 \mod p_1^{V_1} \quad   \ad \quad
      p_2^{\rho_{2,k}} \fm_{1,k} + m_{k}  \equiv 0 \mod p_1^{V_1}.
\end{multline}
We have $ \max( |\fm_{1,k} |,|m_{k}| ) \leq n^{4 \sh s},\; m_k \neq 0$. Let $\ord_{p_1}(m_k) =\beta_k$, \\ $m_k =\ddot{m}_k p_{1}^{\beta_k}$, $(p_{1},\ddot{m}_k)=1$.
Consequently $\beta_k \leq 4\sh s \log_2 n < 0.5 V_1 = [\log_2^3 n]/2$.\\
Therefore
 $\fm_{1,k} =\ddot{\fm}_{1,k}p_{1}^{\beta_k}$, $(p_{1},\ddot{\fm}_{1,k})=1$ and
\begin{equation*}
     -\ddot{\fm}_{1,k}  \ddot{m}^{-1}_{k}  p_2^{\rho_{2,k}} -1 \equiv 0 \mod p_1^{[V_1/2]}.
\end{equation*}
Applying Corollary 1,  we get that this congruence is equivalence.
Hence the assumption of the induction is true:
\begin{equation}  \label{Le14-6}
 m_{k} +\fm_{1,k}p_2^{\rho_{2,k}} = 0 , \quad k=1,...,\sh.
\end{equation}
From \eqref{Par3.2-45}, we derive for $\fs=1$ that
\begin{multline} \label{Le14-7}
\DD_{T_{\fs}, a, \lambda,\tau}   \leq       \sum_{\substack{\rho_{ 1, j+1} -\rho_{1, j} \geq V_1 \\
\rho_{1,j} \in [1,n],  j\in[1,\sh] \\ \rho_{2,j} \in [1,V_1]   }} \;
 \sum_{\substack{ m_j \in I^{*}_{n^{4 \sh s}}, |\fm_{i,j}| \leq n^{4 \sh s}\\  |\fm_{i,j}| \leq  p_i^{\rho^{+}_i}  \\i=1,2, \; j=1,...,\sh}}
 \prod_{j=1}^{\sh} \frac{1}{\bar{m}_j}
     \prod_{i =1}^2 \frac{ \b1(m_{j} +\fm_{1,j}p_2^{\rho_{2,j}} = 0)   }{\bar{\fm}_{i,j}} \\
\leq       \sum_{\substack{\rho_{1,j} \in [1,n],  j\in[1,\sh] \\ \rho_{2,j} \in [1,V_1]   }} \;
 \sum_{\substack{ |\fm_{i,j}| \leq n^{4 \sh s}\\  |\fm_{i,j}| \leq  p_i^{\rho^{+}_i}, \; i=1,2, \; j=1,...,\sh}}
 \prod_{j=1}^{\sh}  \frac{1}{(\bar{\fm}_{1,j})^2 p_2^{\rho_{2,j}}\bar{\fm}_{2,j}}  \\
%
\leq    n^{\sh}     \sum_{\substack{\rho_{2,j} \in [1,n]\\  j \in [1,\sh]}}    \;
 \sum_{\substack{ |\fm_{i,j}| \leq n^{4 \sh s}\\  |\fm_{2,j}| \leq  p_2^{\rho^{+}_i}, \; i=1,2, \; j=1,...,\sh}}
 \prod_{j=1}^{\sh}  \frac{1}{(\bar{\fm}_{1,j})^2 p_2^{\rho_{2,j}}\bar{\fm}_{2,j}} \\
%
\ll    n^{\sh}     \sum_{\substack{\rho_{2,j} \in [1,n]\\  j \in [1,\sh]}}    \;
 \prod_{j=1}^{\sh}  \frac{(\rho^{+}_{2})^{\sh}}{ p_2^{\rho_{2,j}}}
  \ll   n^{\sh}     \sum_{\rho^{+}_{2} \in [1,n] }    \;
   \frac{(\rho^{+}_{2})^{2\sh}}{ p_2^{\rho^{+}_{2}}} \ll n^{\sh} \quad \for \quad \rho^{+}_{2} = \max_{1 \leq j \leq \sh} \rho^{+}_{2,j}.
\end{multline}
By \eqref{Par3.2-3}, assertion \eqref{Le3-401} and Lemma 6 are proved. \qed \\
\subsection{\textcolor{blue}{Case $\fs\geq 2$. General estimates}}
 Let
\begin{multline} \label{Le7-1}
\d1_{k,\brho}:= \delta_{  p_{1}^{ \rho^{+}_{1} -   \rho_{1,\ff_{k-1}}  }   }
 \left( \sum_{j=\ff_{k-1}+1}^{\sh}   (      \fm_{1,j} + m_{j} M_{1, \brho_{j}} ) p_{1}^{ \rho^{+}_{1} -   \rho_{1,j}  }  \right)
   =
      \delta_{  p_{1}^{ \rho^{+}_{1} -   \rho_{1,\ff_{k-1}}  }   }    (L_{k-1,2} )\\
 =\delta_{  p_{1}^{ \rho^{+}_{1} -   \rho_{1,\ff_{k-1}}  }   }    (  L_{k,1}  p_{1}^{\rho^{+}_{1}-\rho_{1,\ff_{k}}}  + L_{k,2} )
 =  \delta_{p_i^{\rho_{1,\ff_{k}} -   \rho_{1,\ff_{k-1}}}} \big( L_{k,1} + L_{k,2} p_{1}^{-\rho^{+}_{1}+\rho_{1,\ff_{k}}} \big) \quad \with\\
L_{k,1} =   \sum_{j=\ff_{k-1}+1}^{\ff_{k} }   (      \fm_{1,j} + m_{j} M_{1, \brho_{j}} ) p_{1}^{\rho_{1,\ff_{k}}  -   \rho_{{1},j}  }, \;
L_{k,2}=\sum_{j=\ff_{k}+1}^{\sh}   (      \fm_{1,j} + m_{j} M_{1, \brho_{j}} ) p_{1}^{ \rho^{+}_{1} -   \rho_{1,j}  }
\end{multline}
and $ \ff_0=\rho_{1,0}=0, \; \ff_{a_1} = \sh   $.
Note that by \eqref{Beg-0a}, $\delta_M(a)=0$ for non-integer $a$.\\ \\
{\bf Lemma 7.} {\it Let  $k \in [1,a_{1}]$.
%
Then}
\begin{equation*}
    \delta_{p_i^{\rho^{+}_{1}}}(\tilde{m}_{1} ) = \prod_{k=1}^{a_{1}} \one_{k,\brho}, \quad \;\; \where \quad \;\; \tilde{m}_{1} =
     \sum_{j=1}^{\sh}   \Big(      \fm_{1,j} + m_{j} M_{1, \brho_{j}} \Big)
      p_i^{ \rho^{+}_{1} -   \rho_{1,j}} .
\end{equation*} \\
{\bf Proof.} 
Suppose that  $ \d1_{k,\brho} =1$. In view of \eqref{Le7-1}, we get
that
\begin{equation*}  
L_{k,1}  p_{1}^{\rho^{+}_{1}-\rho_{1,\ff_{k}}}  + L_{k,2} \equiv 0 \mod p_{1}^{\rho^{+}_{1}-\rho_{1,\ff_{k-1}}}, \quad
 L_{k,2}  \equiv 0 \mod p_{1}^{  \rho^{+}_{1} -   \rho_{1,\ff_{k}}} \quad \ad  \quad \d1_{k+1,\brho} =1
\end{equation*}
and
\begin{equation} \label{Le7-2a}
 \d1_{k,\brho}= \d1_{k,\brho} \d1_{k+1,\brho}, \qquad  \tilde{m}_{1} = L_{0,2}, \qquad  \delta_{p_i^{\rho^{+}_{1}}}(\tilde{m}_{1} ) = \delta_{p_i^{\rho^{+}_{1}}}(L_{0,2} ) = \d1_{1,\brho}.
\end{equation}
By \eqref{Le7-1}, we have
\begin{equation}  \label{Le7-3}
  L_{a_{1},2}=0  \quad \ad \quad  \d1_{a_{1}+1,\brho} =  \delta_{p_i^0} ( L_{a_{1},2} ) =   \delta_{1}(0)= 1.
\end{equation}
%
Using \eqref{Le7-2a} and \eqref{Le7-3}, we obtain
\begin{equation*}
   \delta_{p_i^{\rho^{+}_{1}}}(\tilde{m}_{1} ) = \d1_{1,\brho} =\d1_{1,\brho}\d1_{2,\brho} \times \cdots \times \d1_{a_{1},\brho}  \d1_{a_{1}+1,\brho} =
 \d1_{1,\brho} \times \cdots \times \d1_{a_{1},\brho}.
\end{equation*}
Therefore, Lemma 7 is proved. \qed \\ \\

Let $\rho^{(k)} =\{\rho_j  : j \in (\ff_{k-1}, \ff_k]\}$ and let  $m^{(k)} =\{(m_j, \fm_{i,j}) : j \in (\ff_{k-1}, \ff_k], \; i=1,...,s\}$ for $k =1,..., a_1$. We introduce a sequence of indicator functions $(\fa(\rho^{(k)}, m^{(k)},k ))_{k=1}^{a_1}$.

 By \eqref{Le6-1} the calculation of $\SS_{m, \tau}$ contain the summation of $  \varpi_{\brho,\bm,1} \;  \zeta_{\brho,\fm} \;
    \chi_{\brho} $ $\delta_{p_i^{\rho^{+}_{1}}}(\tilde{m}_{1} )  $ over $\rho_j$ with
    $j \in [1, \sh]  $. We want to reduce this calculation to the sequential summations of
  $ \fa(\rho^{(k)}, m^{(k)},k )\;
    \chi_{\brho} \; \delta_{p_i^{\rho^{+}_{1}}}(\tilde{m}_{1} )  $ over $\rho_j$ with
       $j \in (\ff_{k-1}, \ff_k]  $  for $k =1,..., a_1$.
 We will construct these $\fa(\rho^{(k)}, m^{(k)},k )$ in Section 3.5 by a decomposition of the function $ \varpi_{\brho,\bm,1}  $. In this section, we only need to know that $\fa(\rho^{(k)}, m^{(k)},k ) \in \{0,1\}$. \\ \\ \\
{\bf Lemma 8.} {\it Let  $k  \in [1,a_{1}]$, $\fa(\rho^{(k)}, m^{(k)},k ) \in \{0,1\}$,
$\fh_k \in  [\ff_{k-1}+1,\ff_{k}]$  and let
\begin{multline} \label{Le8-0}
 S_{k,\fa}: = \sum_{\substack{\rho_{j} \in \UU_{T_{\fs}}\\  j \in [\ff_{k-1}+1,\ff_{k}]}} \fa(\rho^{(k)}, m^{(k)},k )\;
  \chi_{1,k,\rho} \; \d1_{k,\rho},    \\
  \vartheta_{\rho, \fh_k}:= \sum_{\substack{\rho_{i, \fh_k} \in [1,n] \\ i \in [2,s] }} \fa(\rho^{(k)}, m^{(k)},k ) \;  \chi_{1,k,\rho} \d1_{k,\rho}
  \quad \ad \quad \cW_{\fh_k}  :=  \sum_{V_1 \leq \rho_{1,\fh_k} \leq n} \vartheta_{\rho, \fh_k}.
 \end{multline}
 Then}
\begin{equation}  \label{Le8-1}
S_{k,\fa} \leq \sum_{\substack{\brho_{j} \in \UU_{T_{\fs}}   \\  j \in [\ff_{k-1}+1,\ff_{k}] \setminus \{\fh_k\} }} \cW_{\fh_k} , \quad \;\;  \with  \;\;  \quad  \vartheta_{\rho_{1, \ff_k}} \in \{0,1\}, \quad \cW_{\fh_k} \leq 2V_1+1 \;\;\;\for \; \lambda_{1,k} \geq 2
\end{equation}
and
\begin{equation}   \label{Le8-1a}          
 S_{k,\fa} \ll n^{(\lambda_{1,k}-1) (\fs-1) +1} V_1^{\lambda_{1,k} s} \;\;\;\; \for \; \lambda_{1,k} \geq 1.
\end{equation} \\
{\bf Proof.} 
The first assertion is obvious (see also \eqref{Le6-2}). Now we examine
the second  assertion of statement \eqref{Le8-1}.
In view of \eqref{Le7-1}, we have
\begin{equation*}
   \one_{k,\rho}=
   \delta_{p_i^{\rho^{+}_{1} -   \rho_{1,\ff_{k-1}}}} \big(  L_{k-1,2} \big), \;\; L_{k,2}=\sum_{j=\ff_{k}+1}^{\sh}   (      \fm_{1,j} + m_{j} M_{1, \rho_{j}} ) p_{1}^{ \rho^{+}_{1} -   \rho_{1,j}  }.
\end{equation*}
Suppose that $\one_{k,\rho}= 1$.
Then
\begin{equation} \nonumber  
     \sum_{j=\ff_{k-1}+1}^{\sh}   ( \fm_{1,j} + m_{j} M_{1, \brho_{j}} )
      p_{1}^{ \rho^{+}_{1} -   \rho_{1,j}}    \equiv 0 \mod p_{1}^{\rho^{+}_{1} -  \rho_{1,\ff_{k-1}}  }.
\end{equation}
We fix $ \brho_j$ for $j \in \{1,...,\sh \} \setminus \{ \fh_k\}$ and we have
\begin{multline*} \nonumber  
   ( \fm_{1,\fh_k} + m_{j} M_{1, \brho_{\fh_k}} )
      p_{1}^{ \rho^{+}_{1} -   \rho_{1,\fh_k}}  \equiv \alpha_0 \mod  p_{1}^{\rho^{+}_{1} -  \rho_{1,\ff_{k-1}}  }, \quad \where \\
-\alpha_0 =
     \sum_{\ff_{k-1}+1 \leq j \leq \sh, \; j \neq \fh_k}   ( \fm_{1,j} + m_{j} M_{1, \brho_{j}} )
      p_{1}^{ \rho^{+}_{1} -   \rho_{1,j}}    \equiv 0 \mod p_{1}^{\rho^{+}_{1} -  \rho_{1,\fh_k} }.
\end{multline*}
We see that $ \alpha_0 \equiv 0  \mod  p_{1}^{\rho^{+}_{1} -  \rho_{1,\fh_{k}} } $.
Let $\alpha_1 = \alpha_0  p_{1}^{-\rho^{+}_{1} + \rho_{1,\fh_{k}} }$. Then
\begin{equation}  \nonumber      
   \fm_{1,\fh_k} + m_{{\fh_k}} M_{1, \brho_{\fh_k}}
  \equiv \alpha_1 \mod p_{1}^{\rho_{1,\fh_k}-  \rho_{1,\ff_{k-1}}  }.
\end{equation}
Hence
\begin{equation}  \nonumber     
 m_{{\fh_k}} M_{1, \brho_{\fh_k}}
  \equiv \alpha_2 \mod p_{1}^{\rho_{1,\fh_k}-  \rho_{1,\ff_{k-1}}  }, \;\; \with \;\;\alpha_2 = \alpha_1 - \fm_{1,\fh_k} .
\end{equation}
By \eqref{Le5-1a} and \eqref{Par3.2-4}, we get $0 <|m_{\fh_k}| \leq n^{4 \sh s}$.
Let $\beta ={\rm ord}_{p_{1}}(m_{\fh_k})$, $m_{\fh_k} =\ddot{m}_{\fh_k} p_{1}^{\beta}$, $( \ddot{m}_{\fh_k} ,p_{1})=1$.
Hence $\beta \ll \log n$.  According to \eqref{Le6-2} and \eqref{Le8-0}, we have
\begin{equation*}
      \rho_{1,\fh_k}-  \rho_{1,\ff_{k-1}} > V_1 = [\log_2^3 n] \quad \ad \quad \rho_{1,\fh_k}-  \rho_{1,\ff_{k-1}} -\beta> V_1-\beta \geq [V_1/2]  .
\end{equation*}
 Using  \eqref{Beg-3} and taking into account that $\rho_{1,\fh_k} \geq V_1 $, we obtain
\begin{equation*}
       \prod_{2 \leq i \leq s} p_i^{- \rho_{i,\fh_k}} \equiv M_{1, \brho_{\fh_k}} \equiv
\alpha_2/ \ddot{m}_{\fh_k} \mod p_{1}^{[V_1/2]}.
\end{equation*}
Suppose that $( \rho_{2,\fh_k}^{'},   ...,\rho_{s,\fh_k}^{'}) $ and
$(\rho_{2,\fh_k}^{''},  ...,\rho_{s,\fh_k}^{''}) $
 are  two different solutions of this congruence.
Then
\begin{equation*}  
 \prod_{2 \leq i \leq s } p_i^{\rho_{i,\fh_k}^{'}} - \prod_{2 \leq i \leq s }
  p_i^{\rho_{i,\fh_k}^{''}} =
 \prod_{2 \leq i \leq s  }  p_i^{\rho_{i,\fh_k}^{''}}
\Big( \prod_{2 \leq i \leq s } p_i^{\rho_{i,\fh_k}^{'}   - \rho_{i,\fh_k}^{''}}  - 1\Big)	
	  \equiv 0 \mod p_{1}^{[V_1/2]}.
\end{equation*}
Here $ V_1/2 = [\log_2^3 n]/2 \gg C_1 \log^2 n $. Therefore, we can apply Corollary 1. We get that this congruence is the equality, having only one solution. Namely, $\rho_{i,\fh_k}^{'} = \rho_{i,\fh_k}^{''}$ $(i=2,...,s)$.
By  \eqref{Le8-0}, we get $ \vartheta_{\rho_{1, \fh_k}} \leq 1$.
From conditions of the lemma, we have $   \fa(\rho^{(k)}, m^{(k)},k ) \; \chi_{1,k, \fh_k,\rho}\; \d1_{k,\brho}  \in \{0,1\} $. Hence
$\vartheta_{\rho_{1, \fh_k}} \in \{0,1\}$ .\\

Let's consider
the third  assertion of statement \eqref{Le8-1}:\\
In view of \eqref{Le6-2}, we have
\begin{equation} \label{Le9-100}
   |\rho_{1,j_1}-  \rho_{1,j_2} | \leq V_1 \quad   \for \;\; j_1,j_2 \in [\ff_{k-1}+1,\ff_k] \;\; \ad \;\;  \chi_{1,k,\rho}=1.
\end{equation}
Bearing in mind that $\lambda_{1,k} \geq 2 $, for fixed $ (\brho_j \; | \; j \in [\ff_{k_1}+1,\ff_k] \setminus \fh_k) $, we have from \eqref{Le9-100} and
\eqref{Le8-0}, that
\begin{equation}   \nonumber        
\cW_{\fh_k}  =  \sum_{V_1 \leq \rho_{1,\fh_k} \leq n}  \vartheta_{\rho, \fh_k}
 \leq 2V_1 +1.
\end{equation}

Let's consider
 assertion  \eqref{Le8-1a}:\\
Let $\lambda_{1,k} \geq 1 $. Taking into account \eqref{Par3-5}, \eqref{Le6-2} and that $\vartheta_{\rho_{1, \fh_k}} \in \{0,1\}  $, we get from Lemma 7:
\begin{equation}   \nonumber        
 S_{k,\fa} \ll  \sum_{\substack{\brho_{j} \in \UU_{T_{\fs}}   \\  j \in [\ff_{k-1}+1,\ff_{k}]
 \setminus \{\fh_k\}}} \sum_{V_1 \leq \rho_{1,\fh_k} \leq n} \vartheta_{\rho_{1, \fh_k}}
  \ll
 n^{(\lambda_{1,k}-1) (\fs-1) +1} V_1^{\lambda_{1,k} s }.
\end{equation}
Therefore, Lemma 8 is proved. \qed  \\
%
%
%

In the following Lemma 9, we show that the calculation of the sum:
\begin{equation}  \label{Lem11-0}
           \SS_{m, \tau,\fa}: =       \sum_{\brho_{j} \in \UU_{T_{\fs}}, j\in[1,\sh]} \;
           \prod_{k=1}^{a_1}\fa(\rho^{(k)}, m^{(k)},k )  \;
    \chi_{\brho}  \; \delta_{p_1^{\rho^{+}_{1}}}(\tilde{m}_{1} )
\end{equation}
can be reduced to the product of simple sums over $ \rho_j \in \UU_{T_{\fs}}$ for $j \in [\ff_{k-1}+1, \ff_k]  $. If $\fa(\rho^{(k)}, m^{(k)},k )=1$ for all
  $\rho^{(k)}, m^{(k)},k $, then we will write $1$ instead of $\fa$. \\ \\ \\
{\bf Lemma 9.} {\it Let
\begin{equation*} \nonumber 
 \dddot{\cS}_{k_1,k_2}:=    \max_{\brho_1,...,\brho_{\sh}}  \hat{\cS}_{k_1,k_2} , \;\; \;\; \quad
 \hat{\cS}_{k_1,k_2}=        \sum_{\substack{\brho_{j} \in \UU_{T_{\fs}} , \quad
    \\ j\in[\ff_{k_1}+1,\sh]}}  \prod_{k=k_1+1}^{k_2}    \fa(\rho^{(k)}, m^{(k)},k ) \;
   \chi_{1,k,\brho} \; \one_{k,\brho}.
\end{equation*}
Then}
\begin{equation} \label{Lem11-1}
    \SS_{m, \tau, \fa} =  \dddot{\cS}_{0,a_1}=\hat{\cS}_{0,a_1}  \leq
   \prod_{k=1}^{a_{1}} \ddot{S}_{k,\fa} \;\; \ad \;\;
     \dddot{\cS}_{k_1,k_2} \leq \prod_{k=k_1+1}^{k_2} \ddot{S}_{k,\fa}
 \;\; \with \;\;   \ddot{S}_{k,\fa}=    \max_{\rho_1,...,\rho_{\sh}}   S_{k,\fa},
 \end{equation}
\begin{equation} \nonumber
    S_{k,\fa} = \sum_{\substack{\rho_{j} \in \UU_{T_{\fs}}   \\  j \in [\ff_{k-1}+1,\ff_{k}]}}
\fa(\rho^{(k)}, m^{(k)},k ) \;  \chi_{1,k,\rho} \; \one_{k,\rho}
 \end{equation}
and
\begin{equation}   \label{Lem11-1a}
   \SS_{m,\tau,\fa}  \ll n^{\sh (\fs-1) -(\fs-2)a_{1}} \log^{3\sh s} n, \quad
     \dddot{\cS}_{k_1,k_2}  \ll \prod_{k=k_1+1}^{k_2} n^{(\lambda_{1,k}-1) (\fs-1) +1}  \log^{3 \lambda_{1,k} s} n.
\end{equation} \\ \\
{\bf Proof.} Let's consider statement \eqref{Lem11-1}:\\
By  \eqref{Lem11-0}, \eqref{Le7-1}, \eqref{Le6-2} and Lemma 7, we get
\begin{equation} \label{Le10-01}
           \SS_{m, \tau,\fa} \leq       \sum_{\rho_{j} \in \UU_{T_{\fs}}, j\in[1,\sh]} \;\prod_{k=1}^{a_{1}}   \fa(\rho^{(k)}, m^{(k)},k ) \;
    \chi_{1,k,\rho} \; \one_{k,\rho} ,
\end{equation}
with
\begin{equation*}
 \one_{k,\rho}= \delta_{  p_{1}^{ \rho^{+}_{1} -   \rho_{1,\ff_{k-1}}  }   }    (L_{k-1,2} ), \qquad \qquad
L_{k,2}=\sum_{j=\ff_{k}+1}^{\sh}   (      \fm_{1,j} + m_{j} M_{1, \brho_{j}} ) p_{1}^{ \rho^{+}_{1} -   \rho_{1,j}  }
\end{equation*}
and
\begin{equation*}
  \chi_{1,k,\brho} =\prod_{ j \in [\ff_{k-1}+1,\ff_{k}] }  \chi_{1,k,j,\brho}, \qquad \qquad
\chi_{1,k,j,\brho}=\b1\big(  0 \leq \rho_{1,j+1 }  - 	 \rho_{1,j }	\leq V_1 \big).
\end{equation*}
for $j, j+1 \in (\ff_{k-1},\ff_{k}]$.

Suppose that $a_{1}=1$.  According to \eqref{Le6-2}, we have
  $\ff_0=0, \; \ff_{a_{1}}=\sh $.\\
Hence  $[\ff_{0}+1,\ff_{1}] =[1,\sh] $, $ \SS_{m, \tau,\fa}=   S_{1,\fa} $ and \eqref{Lem11-1} follows.\\

Let's consider the case $a_{1} \geq 2$. 
%
%
%
Bearing in mind that $ \fa(\rho^{(k)}, m^{(k)},k )$ $\;\chi_{1,k,\brho} \; \one_{k,\brho}  $ does not depend  on  $\brho_{j} $ for $j \leq \ff_{k-1}$, we have
\begin{align*}
 \hat{\cS}_{k_1,k_2}&=     \sum_{\substack{\brho_{j} \in \UU_{T_{\fs}} , \quad
    \\ j\in[\ff_{k_1}+1,\ff_{k_2}]}}  \prod_{k=k_1+1}^{k_2}    \fa(\rho^{(k)}, m^{(k)},k ) \;
   \chi_{1,k,\brho} \; \one_{k,\brho} \\
    &=  \sum_{\substack{\brho_{j} \in \UU_{T_{\fs}} , \quad
    \\ j\in[\ff_{k_1+1}+1,\ff_{k_2}]}}  \prod_{k=k_1+2}^{k_2}
   \fa(\rho^{(k)}, m^{(k)},k ) \;   \chi_{1,k,\brho} \; \one_{k,\brho} \\
   & \qquad\qquad\qquad \times   \sum_{\substack{\brho_{j} \in \UU_{T_{\fs}} , \quad
    \\ j\in[\ff_{k_1}+1, \ff_{k_1+1}  ]}} \fa(\rho^{(k_1+1)}, m^{(k_1+1)},k_1+1 )
      \chi_{1,k_1+1,\brho}  \one_{k_1+1,\brho} \\
&=
 \sum_{\substack{\brho_{j} \in \UU_{T_{\fs}} , \\
    \\ j\in[\ff_{k_1+1}+1,\ff_{k_2}]}}  \prod_{k=k_1+2}^{k_2}
  \fa(\rho^{(k)}, m^{(k)},k ) \;  \chi_{1,k,\brho} \; \one_{k,\brho}  \;  S_{k_1+1,\fa} \\
      &\leq
 \sum_{\substack{\brho_{j} \in \UU_{T_{\fs}} , \\
    \\ j\in[\ff_{k_1+1}+1,\ff_{k_2}]}}  \prod_{k=k_1+2}^{k_2}
   \fa(\rho^{(k)}, m^{(k)},k ) \;  \chi_{1,k,\brho} \; \one_{k,\brho}  \;  \ddot{S}_{k_1+1,\fa}=
  \hat{\cS}_{k_1+1,k_2} \; \ddot{S}_{k_1+1,\fa}
\end{align*}
and $\dddot{\cS}_{k_1,k_2}  \leq \dddot{\cS}_{k_1+1,k_2} \; \ddot{S}_{k_1+1,\fa} $
for $k_1 \in [0,a_1-1] $.\\
We have $ \hat{\cS}_{a_{1}-1,a_1}= S_{a_{1}, \fa} \leq \ddot{S}_{a_{1}, \fa} =\dddot{\cS}_{a_{1}-1,a_1}$. Therefore
\begin{equation*}  
   \dddot{\cS}_{k_1,k_2}   \leq \ddot{S}_{k_1+1,\fa} \times \cdots\times  \ddot{S}_{k_2,\fa}  \quad \ad \quad   \SS_{m, \tau,\fa}=\dddot{\cS}_{0,a_1}   \leq \ddot{S}_{1,\fa} \times \cdots\times  \ddot{S}_{a_{1},\fa}
\end{equation*}
and \eqref{Lem11-1} follows.\\

Let's consider statement \eqref{Lem11-1a}:\\
Using \eqref{Beg-22}, \eqref{Par3-5} and Lemma 8, we obtain
\begin{multline} \label{Le10-10}
 \dddot{\cS}_{k_1,k_2}  \ll  \prod_{k=k_1+1}^{k_2}    n^{(\lambda_{1,k}-1) (\fs-1) +1}V_1^{\lambda_{1,k}s} \quad \ad \quad \SS_{m,\tau,\fa} \ll
  n^v V_1^{\sh s} \quad \with \\
   v=\sum_{k=1}^{a_{1}} (\lambda_{1,k} (\fs-1) -\fs+2 )
  =\sh (\fs-1) -(\fs-2)a_{1}, \;\; V_1 = [\log_2^3 n].
\end{multline}
Hence, Lemma 9 is proved. \qed  \\  \\
%
%
%
%
%
{\bf Lemma 10.} {\it Let $s \geq 3$, $s > \fs \geq 1$ or $s=\fs$  and
  $ a_{1}  \neq \sh/2 $.  Then}
\begin{equation}            \label{Le12-1}
   \SS_{m,\tau,\fa} \ll  n^{\sh s/2-1/3}.
\end{equation}
{\bf Proof.} %
In view of \eqref{Par3-7}, we have
 $a_{1} = \max_{i \in T_{\fs}} a_i $. Taking into account \eqref{Beg-26} and \eqref{Lem11-0}, we get that estimate \eqref{Par3-5a} for $\SS_{m,\tau}$ remains correct for $\SS_{m,\tau,\fa}$:
\begin{equation}   \label{Par3-5b}
 \SS_{m,\tau,\fa}  \ll
  n^{a_1+ \cdots +a_{\fs}} V_1^{\sh s} \ll n^{\fs a_1} \log_2^{3 s \sh} n.
\end{equation}
Let's consider the case $ a_{1} \geq (\sh+1)/2$.
By Lemma 9, we get
\begin{multline}         \nonumber 
    \SS_{m,\tau,\fa} \ll  n^{v_0} \log^{3 \sh s} n, \quad v_0= \sh(\fs-1) -(\fs-2)a_{1} \leq
     \sh(\fs-1) -(\fs-2) (\sh+1)/2 \\
                =\sh \fs/2 - (\fs-2)/2.
\end{multline}
Suppose that $\fs \geq 3$. Then $v_0 \leq \sh s/2 -1/2$. If $\fs =2$, then $v_0 = \sh \fs/2 \leq \sh s/2 -1/2$.
If $\fs =1$, then $v_0 = a_1 \leq \sh   \leq \sh s/2 -1/2$ for $s \geq 3$.
 Hence \eqref{Le12-1} is true for  $a_{1} \geq (\sh+1)/2$. \\

Let's consider the case $ a_{1} \leq \sh/2$.
Suppose that $s >\fs$.
From \eqref{Par3-5b}, we get
\begin{equation}         \nonumber  
    \SS_{m,\tau,\fa} \ll  n^{\sh \fs/2} \log^{3\sh s} n \ll  n^{\sh s/2 -1/3}.
\end{equation}
Now suppose that $s=\fs$. According to the condition of the lemma, $ a_{1} \neq \sh/2$.
Therefore,   $ a_{1} < \sh/2$ and
\eqref{Le12-1} follows from \eqref{Par3-5b}.

Hence, Lemma 10 is proved. \qed \\
\subsection{\textcolor{blue}{Case $\fs \geq 2$ and there exists $k_0$ with $\lambda_{1,k_0} =1$}}

%
%
%
{\bf Lemma 11.} {\it Let $k  \in [1,a_{1}]$, $\fa(\rho^{(k)}, m^{(k)},k ) \in \{0,1\}$. Then}
\begin{multline}  \label{Le9-0}
    S_{k,\fa} \leq \tilde{S}_{k,\fa}   +O(n^{(\lambda_{1,k}-1) (\fs-1)}V_1^{ \sh s}) , \\
     \where \quad
 \tilde{S}_{k,\fa}= \sum_{\substack{\brho_{j} \in \UU_{T_{\fs}} \\  j \in [\ff_{k-1} +1,\ff_{k}]}}
  \fa(\rho^{(k)}, m^{(k)},k )  \chi_{1,k,\brho} \tilde{\d1}_{k, \brho},   \\
\tilde{\d1}_{k, \brho}:=\delta_{p_{1}^{V_1}}( L_{k,1}),   \;\; \;\; V_1 = [\log_2^3 n].
\end{multline} \\
{\bf Proof.}  By Lemma 8, we obtain
\begin{multline}  \label{Le9-1}
S_{k,\fa} \leq \sum_{\substack{\brho_{j} \in \UU_{T_{\fs}}   \\  j \in [\ff_{k-1}+1,\ff_{k})}}
 \sum_{V_1 \leq \rho_{1,\ff_k} \leq n} \vartheta_{\rho_{1, \ff_k}}, \\
 \vartheta_{\rho_{1, \ff_k}}= \sum_{\substack{\rho_{i, \ff_k} \in [1,n] \\ i \in [2,s] }} \fa(\rho^{(k)}, m^{(k)},k )\;  \chi_{1,k, \rho} \;  \d1_{k,\brho}, \quad  \where \quad
              \vartheta_{\rho_{1, \ff_k}} \in \{0,1\}.
\end{multline}
Similarly to \eqref{Le8-0}, we have
\begin{equation} \nonumber
   \tilde{S}_{k,\fa} \leq \sum_{\substack{\brho_{j} \in \UU_{T_{\fs}}   \\  j \in [\ff_{k-1}+1,\ff_{k})}}
 \sum_{V_1 \leq \rho_{1,\ff_k} \leq n}  \ddot{\vartheta}_{\rho_{1, \ff_k}},
\end{equation}
where
\begin{equation} \label{Le9-2}
   \ddot{\vartheta}_{\rho_{1, \ff_k}}= \sum_{\substack{\rho_{i, \ff_k} \in [1,n] \\ i \in [2,s]}} \fa(\rho^{(k)}, m^{(k)},k ) \; \chi_{1,k, \rho} \; \tilde{\d1}_{k, \brho}.
\end{equation}
We see that $  \ddot{\vartheta}_{\rho_{1, \ff_k}}$ is a  non-negative  integer. Let
\begin{equation} \label{Le9-4}
    \beta_{\rho} =  \card \{  V_1 \leq  \rho_{{1},\ff_k} \leq n \; : \;
   \vartheta_{\rho_{1, \ff_k}} =1   \;\; \ad \;\;  \ddot{\vartheta}_{\rho_{1, \ff_k}}=0   \} .
\end{equation}
Using  \eqref{Beg-26} and \eqref{Par3-5}, we get
\begin{multline*}
S_{k,\fa} \leq \tilde{S}_{k,\fa}  + \tilde{\beta}, \; \with \;
    \tilde{\beta}:=  \sum_{\substack{\brho_{j} \in \UU_{T_{\fs}}   \\  j \in [\ff_{k-1}+1,\ff_{k})
}}
      \sum_{V_1 \leq \rho_{1,\ff_k} \leq n} \vartheta_{\rho_{1, \ff_k}} \\
     \times  \b1( \ddot{\vartheta}_{\rho_{1, \ff_k}}   <  \vartheta_{\rho_{1, \ff_k}} )
 \leq \sum_{\substack{\brho_{j} \in \UU_{T_{\fs}}   \\  j \in [\ff_{k-1}+1,\ff_{k})}}
 \beta_{\rho} \\
 \leq n^{(\lambda_{1,k}-1) (\fs-1)}V_1^{ \sh(s-\fs) } \max_{\brho_{2},...,\brho_{\sh}} \beta_{\rho}.
\end{multline*}
Hence  \eqref{Le9-0} may be derived from the inequality
\begin{equation} \label{Le9-5}
        \max_{\brho_{2},...,\brho_{\sh}} \beta_{\rho} \leq (\sh +1)V_1.
\end{equation}
Now we will prove \eqref{Le9-5}:

We fix $ \brho_j$ for $j =2,...,\sh $. Let's consider $\beta_{\rho} $. Suppose that  $\vartheta_{\rho_{1, \ff_k}}=1$  and $ \ddot{\vartheta}_{\rho_{1, \ff_k}}=0  $ for some $\rho_{1, \ff_k} \in [V_1,n] $. We will calculate the number of such $\rho_{1, \ff_k}$.

Applying \eqref{Le7-1}, we obtain
\begin{multline}  \label{Le9-6}
    \one_{k,\brho} =  \delta_{  p_{1}^{ \rho^{+}_{1} -   \rho_{1,\ff_{k-1}}  }   }    (  L_{k,1}  p_{1}^{\rho^{+}_{1}-\rho_{1,\ff_{k}}}  + L_{k,2} ), \quad \with \\
L_{k,1} =   \sum_{j=\ff_{k-1}+1}^{\ff_{k} }   (      \fm_{1,j} + m_{j} M_{1, \brho_{j}} ) p_{1}^{\rho_{1,\ff_{k}}  -   \rho_{{1},j}  }, \;
L_{k,2}=\sum_{j=\ff_{k}+1}^{\sh}   (      \fm_{1,j} + m_{j} M_{1, \brho_{j}} ) p_{1}^{ \rho^{+}_{1} -   \rho_{1,j}  }.
\end{multline}
Suppose that  $\d1_{k,\brho} =1  , \;\; \ad \;\;  \tilde{\d1}_{k, \brho}=0  $.

By \eqref{Le6-2}, $\rho_{1,\ff_{k}} -   \rho_{1,\ff_{k-1}} >V_1$.
If $L_{k,2} =0$ then
\begin{equation*}
L_{k,1} \equiv 0 \mod p_{1}^{\rho_{1,\ff_{k}} -   \rho_{1,\ff_{k-1}}} , \quad \tilde{\d1}_{k, \brho} = \delta_{p_{1}^{V_1}}(L_{k,1}) =1 \quad \ad \quad
\d1_{k, \brho}= \tilde{\d1}_{k, \brho}.
\end{equation*}
We have a contradiction.

Now let $L_{k,2} \neq 0$ and let  $\ord_{p_{1}}(L_{k,2}) =\xi$, $L_{k,2} = L_{k,3}   p_{1}^{\xi} $, $(L_{k,3},p_{1}) =1$.\\
From  \eqref{Le9-6},  we obtain that $ L_{k,2} $ does not depend on $ \rho_{1,\ff_k}$  and
\begin{equation}  \label{Le9-7}
  \one_{k, \brho}=
\delta_{  p_{1}^{ \rho^{+}_{1} -   \rho_{1,\ff_{k-1}}  }   }    (  L_{k,1}  p_{1}^{\rho^{+}_{1}-\rho_{1,\ff_{k}}}  + L_{k,3}   p_{1}^{\xi} )
  =1.
\end{equation}
It is easy to see that this congruence is false for $\xi <   \rho^{+}_{1}- \rho_{i,\ff_k}$.
By  \eqref{Le9-4}, in order to prove \eqref{Le9-5}, it is enough to verify that there is no solution of
 \eqref{Le9-7} with $  \tilde{\d1}_{k, \brho}=0 $ for $\xi  \geq   \rho^{+}_{1} - \rho_{i,\ff_k} +(\sh +1)V_1 $.

Indeed. Let $\xi=  \rho^{+}_{1} - \rho_{i,\ff_k} +(\sh +1)V_1 +\xi_1$ with $\xi_1 \geq 0 $.
For $\rho_{i,\ff_k} \geq  \rho^{+}_{1}$ $ +(\lambda_{1,k}  +1)V_1  -\xi $, we derive from \eqref{Le9-7}  that
\begin{multline*}
1=\one_{k, \brho} =\delta_{  p_{1}^{ \rho^{+}_{1} -   \rho_{1,\ff_{k-1}}  }   }    (  L_{k,1}  p_{1}^{\rho^{+}_{1}-\rho_{1,\ff_{k}}}  + L_{k,3}   p_{1}^{\rho^{+}_{1} - \rho_{i,\ff_k} +(\sh +1)V_1 +\xi_1} ) \\
= \delta_{  p_{1}^{ \rho_{i,\ff_k} -   \rho_{1,\ff_{k-1}}  }   }    (  L_{k,1}  + L_{k,3}   p_{1}^{(\sh +1)V_1 +\xi_1} ) \\
=\delta_{  p_{1}^{ V_1 }   }    (  L_{k,1}  + L_{k,3}   p_{1}^{(\sh +1)V_1 +\xi_1} ) \leq \delta_{p_{1}^{V_1}}( L_{k,1}) \leq \tilde{\d1}_{k, \brho}.
\end{multline*}
We have a contradiction.
Hence, assertion \eqref{Le9-5} and Lemma 11 are proved. \qed  \\ \\
{\bf Lemma 12.} {\it Let $\fs \geq 2$, and let $\lambda_{1,k_0}=1$ for some  $k_0 \in [1,a_{1}]$. Then}
\begin{equation}  \nonumber 
\SS_{m,\tau,\fa}   \ll  n^{\sh s/2-1/4}.
\end{equation}\\
{\bf Proof.}  By Lemma 11, we get
\begin{equation}  \nonumber
    S_{k,\fa} \leq  \sum_{\substack{\brho_{j} \in \UU_{T_{\fs}} \\  j \in [\ff_{k-1} +1,\ff_{k}]}}
  \fa(\rho^{(k)}, m^{(k)},k )  \chi_{1,k,\brho} \delta_{p_{1}^{V_1}}( L_{k,1})
+O(n^{(\lambda_{1,k}-1) (\fs-1)}V_1^{ \sh s}).
\end{equation}
Hence
\begin{multline*}
S_{k_0,\fa} \ll  \dot{S}_{k_0,\fa} + \log^{3 \sh s} n, \quad \dot{S}_{ k_0,\fa} =\sum_{\brho_{\ff_{k_0}} \in \UU_{T_{\fs}} }
 \delta_{p_{1}^{V_1}} \big(  L_{k_0,1} \big), \quad \with \\
L_{k_0,1} =   m_{\ff_{k_0}}
  M_{1, \brho_{\ff_{k_0}} }  + \fm_{{1},\ff_{k_0}}.
\end{multline*}
Suppose that $\dot{S}_{k_0,\fa} \geq 1$. Then
\begin{equation} \nonumber 
  m_{\ff_{k_0}}
  M_{1, \brho_{\ff_{k_0}} }  + \fm_{{1},\ff_{k_0}} \equiv 0 \mod p_{1}^{V_1}.
\end{equation}
Taking into account that
$ M_{1,\brho_j} P_{\brho_j}/ p_1^{\rho_{1,j}}
	\equiv 1 \mod p_1^{\rho_{1,j}} $ (see \eqref{Beg-3}), we derive
\begin{equation}  \nonumber %
 m_{\ff_{k_0}} + \fm_{1,\ff_{k_0}} \prod_{i \in [2,s] }p_i^{\rho_{i,}} \equiv 0 \mod p_{1}^{V_1}.
\end{equation}
Bearing in mind that $\fs \geq 2 $, we obtain from \eqref{Beg-22} and \eqref{Beg-26} that $\max_{ i \neq 1} \rho_{i,\ff_{k_0}}    \geq   \\  V_1= [\log_2^3 n] $. Taking into account Corollary 1 and that $\max( | \fm_{i,\ff_{k_0}} |, |  m_{\ff_{k_0}}  |) \leq n^{4 \sh s}$, we have a contradiction. Thus
\begin{equation} \label{Le15-4}
 S_{k_0,\fa} = O(\log^{3 \sh s} n).
\end{equation}
By Lemma 8 and  Lemma 9, we have
\begin{equation*}
    \SS_{m, \tau, \fa} \leq
   \prod_{k=1}^{a_{1}}  \max_{\brho_1,...,\brho_{\sh}}   S_{k,\fa} , \qquad  S_{k,\fa} \ll
    n^{(\lambda_{1,k}-1) (\fs-1) +1} V_1^{\sh s}.
 \end{equation*}
In view of  \eqref{Le15-4} and \eqref{Le10-10}, we get
\begin{multline} \label{Le15-4a}
  \SS_{m,\tau,\fa} \ll \log^{3 \sh s} n  \; \prod_{1 \leq k \leq a_1, \; k \neq k_0}   n^{(\lambda_{1,k}-1) (\fs-1) +1}\log^{{3 \sh s}} n \\
 \ll  n^{-1} \prod_{k=1}^{a_{1}}    n^{(\lambda_{1,k}-1) (\fs-1) +1}\log^{{3 \sh s}} n
 \ll n^{\sh (\fs-1) -(\fs-2)a_{1} -1} \log ^{3\sh^2 s} n .
\end{multline}
If $\fs=2$, then $  \SS_{m,\tau, \fa} \ll n^{\sh -3/4}$.
Let's consider the case  $\fs \geq 3$.   By Lemma 10, we get
$$
\SS_{m,\tau, \fa} \ll n^{\sh s/2 -1/4} \quad \for \quad s > \fs
 \quad  {\rm or} \quad  s=\fs \quad  \ad \quad    a_{1}  \neq \sh/2.
$$
Now let  $s = \fs \geq 3$ and $a_{1} =\sh/2$. By \eqref{Le15-4a}, we have
\begin{equation} \nonumber
  \SS_{m,\tau,\fa} \ll      n^{\sh (s-1)/2  +\sh/2 -1} \log^{3\sh^2 s } n  \ll  n^{\sh s/2  -1} \log^{3\sh^2 s }
  \ll  n^{\sh  s/2-1/2}.
\end{equation}
Therefore, Lemma 12 is proved. \qed \\
\subsection{\textcolor{blue}{Case  $s=\fs  \geq 2$ and $\lambda_{1,k} =2$ for all $k$}}

The purpose of this subsection is to exclude and consider the case $\lambda_{i,k} =2$ for all $k$ and all $i$ separately from other possible cases.

 Taking into account that $\lambda_{1,k} =2$ for all $k$, we get $\ff_{k}=\Lambda_{1,k}=2k$ for all  $k$ and $\sh=2q$. \\ \\
{\bf Lemma 13.} {\it Let $s=\fs \geq 2$, $\lambda_{1,k} = 2$ for all $k$,
 and let
 there exist $i_1\in [2,s] $, \;  $ k_1 \neq k_2 \leq q$,
  $  l \in [1,a_{i_1}]$, $j_{k_1} \in [2k_1-1,2k_1]$, $j_{k_2} \in [2k_2-1,2k_2]$,
      $l_1,l_2  \in  [\Lambda_{i_1,l-1}+1,\Lambda_{i_1,l}]    $,
                    with
             $\dot{\tau}_{i_1}(l_1) = j_{k_1} $,
               $\dot{\tau}_{i_1}(l_2) =j_{k_2}$.
          Then}
 \begin{equation}  \label{Le16-1}
   \SS_{m, \tau, \fa}  \ll n^{q s-3/4}.
\end{equation} \\
{\bf Proof.} Let's consider the case $k_1 <k_2$. The proof for the case $k_2 >k_1$ is similar.

In view of  \eqref{Lem11-0},   \eqref{Le6-2a} and  \eqref{Le6-2}, we obtain
\begin{equation} \nonumber 
           \SS_{m, \tau,\fa} \leq       \sum_{\rho_{j} \in \UU_{T_{\fs}}, j\in[1,2q]} \; \chi_{i_1,l,\rho} \; \prod_{k=1}^{q}   \fa(\rho^{(k)}, m^{(k)},k ) \;
    \chi_{1,k,\rho} \; \one_{k,\rho}.
\end{equation}
By   \eqref{Le6-2a} and  \eqref{Le7-1}, we have that $ \chi_{1,k,\rho}\one_{k,\rho}$
 does not depend on $\rho_j$ for $j \leq  2k-2$ and  $ \chi_{i_1,l,\rho}$
 does not depend on $\rho_j$ for $j \leq  2k_1-2$. Similarly to the proof of Lemma 9, we have
\begin{multline}  \nonumber  
           \SS_{m, \tau,\fa} \leq       \sum_{\substack{\rho_{j} \in \UU_{T_{\fs}} \\ j\in[2k_1-1,2q]}} \; \max_{\rho_{j}, j\in[2k_1-1,2q]}
            \Bigg(  \sum_{\substack{\rho_{j} \in \UU_{T_{\fs}}\\  j\in[1, 2k_1-2]}} \;
            \prod_{k=1}^{k_1-1}   \fa(\rho^{(k)}, m^{(k)},k ) \;
    \chi_{1,k,\rho} \; \one_{k,\rho} \Bigg)  \\
     \times       \chi_{i_1,l,\rho} \; \prod_{k=k_1}^{q}   \fa(\rho^{(k)}, m^{(k)},k ) \;
    \chi_{1,k,\rho} \; \one_{k,\rho}  \leq \dddot{\cS}_{0,k_1-1} \; A_0
 \quad \with \\
  A_0:=\sum_{\substack{\rho_{j} \in \UU_{T_{\fs}} \\ j\in[2k_1-1,2q]}} \;
                        \chi_{i_1,l,\rho} \; \prod_{k=k_1}^{q}   \fa(\rho^{(k)}, m^{(k)},k ) \;
    \chi_{1,k,\rho} \; \one_{k,\rho}.
\end{multline}
Similarly, we get
\begin{multline} \nonumber 
           A_0 \leq       \sum_{\substack{\rho_{j} \in \UU_{T_{\fs}} \\ j\in[2k_1+1,2q]}} \; \max_{\rho_{j}, j\in[2k_1+1,2q]}
            \Bigg(  \sum_{\substack{\rho_{j} \in \UU_{T_{\fs}} \\ j\in[2k_1-1, 2k_1]}} \;
              \chi_{i_1,l,\rho} \;
             \fa(\rho^{(k_1)}, m^{(k_1)},k_1 ) \;
    \chi_{1,k_1,\rho} \; \one_{k_1,\rho} \Bigg)\\
           \times \prod_{k=k_1+1}^{q}   \fa(\rho^{(k)}, m^{(k)},k )
    \chi_{1,k,\rho} \; \one_{k,\rho}
  \leq B_0 \dddot{\cS}_{k_1,q}  \quad \with \quad  B_0 :=  \max_{\rho_{j}, j\in[2k_1+1,2q]} B_1, \quad \\
   B_1 :=  \sum_{\rho_{j} \in \UU_{T_{\fs}}, j\in[2k_1-1, 2k_1]} \;
              \chi_{i_1,l,\rho} \;
             \fa(\rho^{(k_1)}, m^{(k_1)},k_1 ) \;
    \chi_{1,k_1,\rho} \; \one_{k_1,\rho}.
\end{multline}
Let consider $B_1$:\\
Applying Lemma 8, we get
\begin{equation}   \nonumber 
B_1 \leq \sum_{1 \leq \rho_{i,j_{k_1}} \leq n, i=1,\ldots,s} \; \chi_{i_1,l,\rho}\;  \cW_{\fh_{k_1}} \leq
 (2V_1+1) \sum_{1 \leq \rho_{i,j_{k_1}} \leq n, i=1,\ldots,s} \; \chi_{i_1,l,\rho} ,
\end{equation}
with $ \{j_{k_1},\fh_{k_1}   \} =\{2k_1-1,2k_1\}$.\\
In view of \eqref{Le6-2a}, we have
\begin{equation}  \nonumber
  \chi_{i_1,l,\brho} =\prod_{ j \in [\Lambda_{i_1,l-1}+1,\Lambda_{i_1,l}] }  \chi_{i_1,l,j,\brho},\quad
  \chi_{i_1,k,j,\brho}=\b1\big( 0 \leq  \rho_{i_1,\dot{\tau}_{i_1}(j+1) }  - 	\rho_{i_1,\dot{\tau}_{i_1}(j)} 	\leq V_1\big)
\end{equation}
for $ j, j+1 \in (\Lambda_{i_1,l-1},\Lambda_{i_1,l}],\;\; l \in [1,a_{i_1}]$.

Taking into account that $\dot{\tau}_{i_1}(l_1) = j_{k_1} $,
               $\dot{\tau}_{i_1}(l_2) =j_{k_2}$, $l_1,l_2  \in  [\Lambda_{i_1,l-1}+1,\Lambda_{i_1,l}]  $ and $|l_2- l_1| < 2q$, we get
\begin{equation}   \nonumber 
B_1 \leq (2V_1+1) \sum_{1 \leq \rho_{i,j_{k_1}} \leq n, i=1,\ldots,s} \;
              \b1( |\rho_{i_1,j_{k_1}} -\rho_{i_1,j_{k_2}}  | <2q V_1) \leq n^{s-1}(2V_1+1) (4qV_1+1).
\end{equation}
By Lemma 9, we obtain
\begin{multline}   \nonumber  
   \SS_{m, \tau,\fa} \leq \dddot{\cS}_{0,k_1-1} B_0 \dddot{\cS}_{k_1,q} \\
               \ll n^{-1} \log^{3  s}  \prod_{k=1}^{q} n^{(\lambda_{1,k}-1) (\fs-1) +1}  \log^{3 \lambda_{1,k} s} n \ll n^{-1+ 2q(s-1) -q(s-1)+q}
               \log^{9q s} n \ll n^{qs-3/4} .
\end{multline}
Therefore, Lemma 13 is proved. \qed \\ \\
{\bf Lemma 14.} {\it Let $s=\fs \geq 2$, $\lambda_{1,k}=2$ for all $k \in [1,q]$ and let  $\lambda_{i_1,l} \neq 2$
 for some $i_1 \in [2,s]$ and $ l \in [1,a_{i_1}]$.
 Then}
\begin{equation}  \label{Le14-1a}
 \SS_{m, \tau }  \ll  n^{\sh s/2-1/4}.
\end{equation}\\
{\bf Proof.} Suppose that there exist $i_1 \in [2,s]$ and $ l \in [1,a_{i_1}]$ with  $\lambda_{i_1,l} \geq 3$.

 By \eqref{Par3-5},  we obtain that $\#\{ \dot{\tau}_{i_1}(j) \; : \; j  \in  [\Lambda_{i_1,l-1}+1,\Lambda_{i_1,l}]  \} \geq 3$.\\
Takigng into account that $\lambda_{1,k}=2$ for all $k \in [1,q]$, we get that  there exist
                     $1 \leq  k_1 \neq k_2 \leq q$, with
             $\dot{\tau}_{i_1}(l_1) \in  [2k_1-1, 2k_1] $,
               $\dot{\tau}_{i_1}(l_2) \in [2k_2-1, 2k_2]$.
  for some $l_1, l_2 \in [\Lambda_{i_1,l-1}+1,\Lambda_{i_1,l}]$.\\
From Lemma 13, we get the assertion of the lemma.

Now suppose that $\lambda_{i,k}\leq 2$ for all $i,k$ and there exist $i_1 \in [2,s]$ and $ l \in [1,a_{i_1}]$ with  $\lambda_{i_1,l} =1$.\\
By \eqref{Par3-5} and \eqref{Par3-7},   we get that
$\lambda_{i_1,1}+ \cdots + \lambda_{i_1,a_{i_1}}=2q$ and
 $q=a_{1} \geq a_{i_1}$. Hence
$$
\sum_{j=1}^{a_{i_1}}  \lambda_{i_1,j}= 1+ \sum_{1 \leq j \leq a_{i_1}, j \neq l} \leq  1+ 2(a_{i_1}-1) <2q.
$$
 We have a contradiction.\\
Therefore, Lemma 14 is proved. \qed \\ \\

We have $\sh=2q$. To control $\varpi_{\brho,\bm,1}$, we define  variables $ \fb_{\rho,m,k,0}$
  and $ \fb_{\rho,m,k,1}$ as follows
\begin{equation}\label{Le18-0a}
 \fb_{\rho,m,k,0} =\b1(m_{2k-1}/P_{\brho_{2k-1}}  =- m_{2k}/P_{\brho_{2k}}) \quad \ad \quad  \fb_{\rho,m,k,1}=1 - \fb_{\rho,m,k,0}.
\end{equation}
It is easy to see that
\begin{equation} \nonumber
1 = \prod_{k=1}^q (\fb_{\rho,m,k,0}+\fb_{\rho,m,k,1}) = \sum_{\substack{\fj_i \in \{0,1\} \\ i \in \{1,s\}}} \ddot{\fb}_{\rho,m,\fj} \;\; \with \;\; \ddot{\fb}_{\rho,m,\fj}=
  \prod_{k=1}^q \fb_{\rho,m,k,\fj_k}, \; \fj=(\fj_1,...,\fj_q).
\end{equation}
By \eqref{Lem2-6}, we get that if $\ddot{\fb}_{\rho,m,\fj}=1$  for $\fj=(0,...,0)$ then $ \varpi_{\brho,\bm,2} =1$. Hence  $\varpi_{\brho,\bm,2}  \geq  \ddot{\fb}_{\rho,m,0}$,
  $\varpi_{\brho,\bm,1} = 1- \varpi_{\brho,\bm,2}  \leq  1-\ddot{\fb}_{\rho,m,0}$ and
\begin{equation}\label{Le18-0b}
 \varpi_{\brho,\bm,1} \leq   \sum_{\fj_1,...,\fj_q \in \{0,1\}, \fj_1+\cdots \fj_q \geq 1} \ddot{\fb}_{\rho,m,\fj}.
\end{equation}
Let
\begin{equation}  \label{Le6-100}
           \SS_{m, \tau,\fb_{\fj} } =     \sum_{\brho_{j} \in \UU_{T_{\fs}}, j\in[1,\sh]} \;  \ddot{\fb}_{\rho,m,\fj}  \;
    \chi_{\brho} \; \delta_{p_i^{\rho^{+}_{1}}}(\tilde{m}_{1} ) =
    \sum_{\brho_{j} \in \UU_{T_{\fs}}, j\in[1,\sh]} \;    \prod_{k=1}^q \fb_{\rho,m,k,\fj_k}  \;
    \chi_{\brho} \; \delta_{p_i^{\rho^{+}_{1}}}(\tilde{m}_{1} ).
\end{equation}
From Lemma 5, we obtain
\begin{equation}  \label{Le6-101}
         \SS_{m, \tau} \leq       \sum_{\brho_{j} \in \UU_{T_{\fs}}, j\in[1,\sh]} \;   \varpi_{\brho,\bm,1} \;
    \chi_{\brho} \;   \zeta_{\brho,\fm} \; \delta_{p_1^{\rho^{+}_{1}}}(\tilde{m}_{1} ) \leq  \sum_{\fj_1,...,\fj_q \in \{0,1\}, \fj_1+\cdots +\fj_q \geq 1}
            \SS_{m, \tau,\fb_{\fj} } .
\end{equation} \\ \\
{\bf Lemma 15.} {\it Let $s=\fs \geq 2$,  $\lambda_{i,k}=2$ for all $i \in [1,s]$ and all $k \in [1,q] $. Then}
\begin{equation}\label{Le18-01}
 \SS_{m, \tau }  \ll  n^{\sh s/2-1/4}.
\end{equation} \\
{\bf Proof.} Taking into account that $\lambda_{i,k} =2$ for all $i,k$, we get $\Lambda_{i,k}= \ff_{k}=2k$ for  $i=1,...,s$, $k=1,...,q$ and $\sh=2q$.

Suppose that there exist $i_1 \in [2,s]$, $k_1, l_1,l_2 \in [1,q]$, $l_1 \neq l_2$ and
  $\ell_0,\ell_1   \in \{ 0,1\}$, such that
\begin{equation}\label{Le18-01a}
  2k_1-1= \dot{\tau}_{i_1}(2l_1 -l_0) \quad \ad \quad 2k_1=\dot{\tau}_{i_1}(2l_2 -l_1) .
\end{equation}
Let $\ell_2 =1-\ell_1$ ($\ell_2 \in \{0,1\}$) and $ j_{k_2}:= \dot{\tau}_{i_1}(2l_2 -l_2)  \in [2k_2-1, 2k_2] $ with some $k_2 \in [1,q]$.

Suppose that $k_1=k_2$. Then $ \dot{\tau}_{i_1}(2l_2 -l_2)=2k_1-1 $. Bearing in mind that  $2k_1-1  =\dot{\tau}_{i_1}(2l_1 -l_0)  = \dot{\tau}_{i_1}(2l_2 -l_2) $, we get $\dot{\tau}_{i_1}(2l_1 -l_0)  = \dot{\tau}_{i_1}(2l_2 -l_2) $ and $l_1=l_2$.
 We have a contradiction. Hence $k_1 \neq k_2$.

Taking $j_{k_1}:= 2k_1=\dot{\tau}_{i_1}(2l_2 -l_1)$ and applying Lemma 13 with $l=l_2$, we get assertion \eqref{Le18-01}.

Now let's consider the case when assertion \eqref{Le18-01a} does not true. Indeed, for all $i_1 \in [2,s]$ and all $k_1 \in [1,q]$ there exists $ l  \in [1,q]$ , such that
\begin{equation}\label{Le18-01b}
  \{2k_1-1, 2k_1 \}= \{ \dot{\tau}_{i_1}(2l-1), \dot{\tau}_{i_1}(2l)\}.
\end{equation}
Suppose that $\chi_{\rho}=1$. From \eqref{Par3-5} and \eqref{Le18-01b}, we obtain for $i_1 =1, ...,s$ and $k_1 =1, ...,q$, that
\begin{equation}  \nonumber 
 |\rho_{i_1, 2k_1 } - \rho_{i_1, 2k_1-1}  | =
    |\rho_{i, \dot{\tau}_i(2l) } - \rho_{i, \dot{\tau}_i(2l-1)  } |  <  V_1,
       \;\;\; \; V_1=[\log^3 n].
\end{equation}
Hence, for fixed $(\rho_{1,2k_1-1},\rho_{2,2k_1-1},...,\rho_{s,2k_1-1})$, we have
\begin{equation}  \label{Le18-4}
   \# \{ \brho_{2k_1} -  \brho_{2k_1-1}  \} \leq (2V_1+1)^s   \ll \log^{3s} n
   \quad \for \quad
      \chi_{1,k_1,\brho}=1.
\end{equation}
 In view of \eqref{Le6-101}, we have that in order to prove \eqref{Le18-01}, it is enough to verify that
\begin{equation}  \label{Le6-102}
       \SS_{m, \tau,\fb_{\fj} } \ll n^{ q s-1/4} \quad \forall (\fj_1,...,\fj_q ) \in \{0,1\}^q \quad
       \with \quad \fj_1+\cdots +\fj_q \geq 1.
\end{equation}
We fix $(\fj_1,...,\fj_q ) \in \{0,1\}^q$. We have that there exists $k_1 \in [1,q]$ with   $\fj_{k_1} =1$.\\
Using  \eqref{Le18-0a}, we get
\begin{equation}  \label{Le18-6}
\{0,1 \} \ni \fb_{\rho,m,k_1,1}=1   \Longleftrightarrow                m_{2k_1-1}/P_{\brho_{2k_1-1}} \neq - m_{2k_1}/P_{\brho_{2k_1}} .
\end{equation}
Applying Lemma 9 for $ \fa_{\rho,m,k} =\fb_{\rho,m,k,\fj_k}$, we obtain from \eqref{Le6-100} and \eqref{Le6-2} that
\begin{equation}  \nonumber
    \SS_{m, \tau,\fb_{\fj}} \leq
   \prod_{k=1}^{q} \ddot{S}_{k,\fb_{\fj}},   \quad \with \quad  \ddot{S}_{k,\fb_{\fj}}=    \max_{\rho_1,...,\rho_{\sh}}   S_{k,\fb_{\fj}}
 \end{equation}
and
\begin{equation} \label{Le18-1}
  S_{k,\fb_{\fj}} = \sum_{\substack{\rho_{j} \in \UU_{T_{\fs}}   \\  j \in [2k-1,2k]}} \fb_{\rho,m,k,\fj_k} \; \chi_{1,k,\rho} \; \d1_{k,\rho}.
\end{equation}
According to  Lemma 8, we get
\begin{equation}  
S_{k, \fb_{\fj}} \ll n^{s}V_1^{\sh s}, \quad k=1,...,q.
\end{equation}
We have that in order to prove \eqref{Le6-102} it is enough to verify that
\begin{equation}  \label{Le18-6a1}
  S_{k, \fb_{\fj}} \ll n^{s-3/4} \quad \for \quad k=k_1.
\end{equation}
In view of Lemma 11 and \eqref{Le7-1}, we get
\begin{multline} \label{Le18-2}
  S_{k,\fb_{\fj}} = \tilde{S}_{k,\fb_{\fj}} +O(n^{\fs-1} V_1^{\sh s}), \; \;
   \tilde{S}_{k,\fb_{\fj}}=
   \sum_{\substack{\brho_{j} \in \UU_{T_{\fs}}   \\  j \in [2k-1,2k]}}
\fb_{\rho,m,k,\fj_k} \; \chi_{1,k,\brho}
\; \delta_{p_{1}^{\rho_{1,\ff_{k}} -   \rho_{1,\ff_{k-1}}}}( L_{k,1}) ,  \\
 \where \quad L_{k,1} =   \sum_{j=\ff_{k-1}+1}^{\ff_{k} }   (      \fm_{1,j} + m_{j} M_{1, \brho_{j}} ) p_{1}^{\rho_{1,\ff_{k}}  -   \rho_{{1},j}  } = \\
 (      \fm_{1,\ff_{k}} + m_{j} M_{1, \brho_{\ff_{k}}} )
 + (      \fm_{1,\ff_{k-1}+1} + m_{\ff_{k-1}+1} M_{1, \brho_{\ff_{k-1}+1}} ) p_{1}^{\rho_{1,\ff_{k}}  -   \rho_{{1},\ff_{k-1}+1}  }  .
 \end{multline}
Hence, to prove \eqref{Le18-6a1} it is enough to verify that
\begin{equation}  \label{Le18-6a}
\tilde{S}_{k_1, \fb_{\fj}} \ll n^{s-3/4}.
\end{equation}
Now we will prove \eqref{Le18-6a}:\\
By \eqref{Le18-2}, we need to consider the case $\fb_{\rho,m,k,\fj_k}  \neq 0$. From \eqref{Le18-6}, we have
\begin{equation}  \label{Le18-6c}
                m_{2k_1-1}/P_{\brho_{2k_1-1}} \neq - m_{2k_1}/P_{\brho_{2k_1}}.
\end{equation}
Suppose that $ \chi_{1,k,\brho}=1$ and $\delta_{p_{1}^{\rho_{1,\ff_{k}} -   \rho_{1,\ff_{k-1}}}}( L_{k,1}) =1$.
By \eqref{Le18-2}, we obtain for $f_{k_1} =2k_1$ that
\begin{equation}  \nonumber
\sum_{j=2{k_1}-1}^{2{k_1}} \Big(      \fm_{1,j} + m_{j} M_{{1}, \brho_{j}} \Big)  p_{1}^{\rho_{1,2k_1} - \rho_{{1},j} }
  \equiv 0 \mod p_{1}^{\rho_{1,2k_1} -\rho_{1,2k_1-2} }.
\end{equation}
Let $v_{k_1} = \rho_{{1},2{k_1}} -\rho_{{1},2{k_1}-1}$.
From  \eqref{Le6-2}, we get
\begin{equation}  \nonumber
\rho_{1,2k_1} -\rho_{1,2k_1-2} -  v_{k_1} = \rho_{{1},2{k_1}-1} -\rho_{{1},2{k_1}-2} > V_1
 \quad \ad \quad  0 \leq v_{k_1} \leq V_1.
\end{equation}
Therefore
\begin{equation}  \label{Le18-8}
\sum_{j=2{k_1}-1}^{2{k_1}} \Big(      \fm_{1,j} + m_{j} M_{{1}, \brho_{j}} \Big)  p_{1}^{\rho_{{1},2{k_1}} - \rho_{{1},j} }
  \equiv 0 \mod p_{1}^{v_{k_1}+V_1}, \; v_{k_1} = \rho_{{1},2{k_1}} -\rho_{{1},2{k_1}-1}\le V_1.
\end{equation}

Let's consider the case $v_{k_1} \geq V_1/4$. By  \eqref{Le18-8}, we get $ \fm_{1,2{k_1}} +m_{2{k_1}} M_{{1}, \brho_{2{k_1}}}  \equiv 0 \mod p_{1}^{V_1/4}$, with $ \max( |\fm_{1,2{k_1}} |,|m_{2{k_1}}| ) \leq n^{4 \sh s},\; m_{2{k_1}}\neq 0$. Let $\ord_{p_{1}}(m_{2{k_1}}) =\beta_{k_1}$, $m_{2{k_1}} =\ddot{m}_{2{k_1}} p_{1}^{\beta_{k_1}}$, $(p_{1},\ddot{m}_{2{k_1}})=1$. Hence $\beta_{k_1} \leq 4 \sh s\log_2 n < V_1/8 $. We get from Corollary 1 that the congruence
\begin{equation*}
  ( -\fm_{1,2{k_1}} / p_{1}^{\beta_{k_1}} ) \ddot{m}^{-1}_{2{k_1}}
   \prod_{1 \leq i \leq s, i \neq 1} p_i^{\rho_{i,2{k_1}}} -1 \equiv 0 \mod p_{1}^{[V_1/8]}
\end{equation*}
 is equality. But this is impossible because $\fs=s \geq 2$ and $\rho_{i_2,j} \geq V_1 =[\log_2^3 n]$ for all $j$ (see \eqref{Beg-22} and \eqref{Beg-26}). \\

Let's consider the case $0 \leq v_{k_1} < V_1/4$ and $g_{k_1}:=\fm_{1,2{k_1}} + \fm_{1,2{k_1}-1} p_{1}^{v_{k_1}} =0$.\\
 By  \eqref{Le18-8} and \eqref{Beg-3}, we get
\begin{equation*}
m_{2{k_1}} \prod_{1 \leq i \leq s, i \neq 1} p_i^{\rho_{i,2{k_1}-1 }-\rho_{i,2{k_1}}}\equiv  -p_{1}^{v_{k_1} } m_{2{k_1}-1} \mod p_{1}^{v_{k_1}+V_1}
\end{equation*}
with $ 0< |m_{2k_1 -1}|, |m_{2k_1}|  \leq  n^{4 \sh s}$.
 Let $\ord_{p_{1}}(m_j) =\beta_j$, $m_j =\ddot{m}_j p_{1}^{\beta_j}$, $(p_{1},\ddot{m}_j)=1$, $ j=2{k_1}-1,2{k_1}$.
By Corollary 1, we get that the congruence
\begin{equation*}
  \prod_{1 \leq i \leq s, i \neq 1} p_i^{\rho_{i,2{k_1}-1}-\rho_{i,2{k_1}}} \equiv  - p_{1}^{v_{k_1}+\beta_{2{k_1}-1}-\beta_{2{k_1}}}  \ddot{m}_{2{k_1}-1}/ \ddot{m}_{2{k_1}}  \mod p_{1}^{v_{k_1}+V_1 -\beta_{2k_1}}
\end{equation*}
 is equality (here $V_1-  \beta_{2k_1} \geq V_1/2$,  $\;\beta_{2k_1}, \beta_{2k_1 -1} \ll \log n  $ and  $v_{k_1}+ \beta_{2k_1 -1} - \beta_{2k_1}=0  $). Hence $m_{2{k_1}-1}/P_{\brho_{2{k_1}-1}} = -m_{2{k_1}}/P_{\brho_{2{k_1}}}$.
But according to   \eqref{Le18-6c} , it  is impossible.  \\

Let's consider the case $0 \leq v_{k_1} < V_1/4$ and $g_{k_1} =  \fm_{1,2{k_1}} + \fm_{1,2{k_1}-1} p_{1}^{v_{k_1}}  \neq 0$.
Let
\begin{equation} \nonumber 
  \xi_{k_1} = m_{2{k_1}}  \prod_{i=2}^s p_i^{n + \rho_{i,2k_1 -1} - \rho_{i,2k_1}}   +p_{1}^{v_{k_1}}  m_{2{k_1}-1} \pp_n \quad \with \quad \pp_n = \prod_{i=2}^s p_i^{n} .
\end{equation}
From  \eqref{Le18-4}, we get
\begin{equation} \label{Le18-9}
  \# \{ \xi_{k_1} \; : \; \brho_{2k_1-1}, \brho_{2k_1} \in \UU_{T_{\fs}}, \; \chi_{1,k_1,\brho}=1 \} \leq
  (2V_1+1)^s V_1/4  \ll \log^{3(s+1)} n.
\end{equation}

Taking into account that $\rho_{i,j} \leq n$ for all $i$ and $j$, we obtain that $\xi_{k_1}$ is integer. By \eqref{Beg-3}, we have $M_{i, \br} \equiv \prod_{i=2}^s p_i^{-r_{i}} \mod p_1^{r_1} $. Hence
\begin{equation*}
   M_{{1}, \brho_{2{k_1}}} M^{-1}_{{1}, \brho_{2{k_1}-1}} \equiv
   \prod_{i=2}^s p_i^{\rho_{i,2k_1 -1} - \rho_{i,2k_1}}
   \mod p_{1}^{\rho_0} \;\; \with \;\; \rho_0:=\min(\rho_{1,2{k_1}-1}, \rho_{1,2{k_1}} )
\end{equation*}
and
\begin{equation*}
  \xi_{k_1} \equiv  (m_{2{k_1}}    M_{{1}, \brho_{2{k_1}}} M^{-1}_{{1}, \brho_{2{k_1}-1}} +p_{1}^{v_{k_1}}m_{2{k_1}-1}) \pp_n \mod p_{1}^{\rho_0} .
\end{equation*}
Bearing in mind \eqref{Le18-8} and that $\rho_{0} \geq V_1$, we get
\begin{equation} \label{Le18-10}
 g_{k_1}\pp_n +   \xi_{k_1}   M_{{1}, \brho_{2{k_1}-1}}    \equiv 0 \mod p_{1}^{V_1}.
\end{equation}
We fix $m_{2{k_1}-1}, m_{2{k_1}}, \fm_{1,2{k_1}-1} \fm_{1,2{k_1}}$ and $\rho_{i,2{k_1}} -\rho_{i,2{k_1}-1}$ for $i=1,...,s$.

Let $\beta_1 =\ord_{p_{1}}(g_{k_1})$ and let $g_{k_1}=\ddot{g}_{k_1} p_{1}^{\beta_1}$, $ (p_{1},\ddot{g}_{k_1} )  =1$. Let $\beta_2 =\ord_{p_{1}}(\xi_{k_1})$ and let $\xi_{k_1}=\ddot{\xi}_{k_1} p_{1}^{\beta_2}$,   $ (p_{1},\ddot{\xi}_{k_1} )  =1$.
We see that
$$
\beta_1 \leq v_{k_1}+ \log_{p_1} (2n^{4\sh s}) \leq v_{k_1}+ V_1/2 < 3V_1/4.
$$
By  \eqref{Le18-10}, we get $\beta_1 = \beta_2$  and
\begin{equation}  \label{Le18-12}
    \prod_{1 \leq i \leq s, i \neq 1} p_i^{\rho_{i,2{k_1}-1}} \equiv - \ddot{\xi}_{k_1} /(\ddot{g}_{k_1}\pp_n) \mod p_{1}^{V_1/4}.
\end{equation}
By Corollary 1, we get for fixed $\rho_{1,2k_{1}-1}$,  that the number of solutions in the variables $(\rho_{2,2k_1-1},\rho_{3,2k_1-1},...,\rho_{s,2k_1-1})$ of this congruence is no more than one.
Therefore, the number of vectors $(\rho_{1,2k_1-1},\rho_{2,2k_1-1},...,\rho_{s,2k_1-1})$ satisfying  \eqref{Le18-12} is less than $n+1$.
According to  \eqref{Le18-9}, there are
 only $O(\log^{3(s+1)} n )$ opportunities to  choose $\xi_{k_1}$.
Applying \eqref{Le18-2}, we obtain  $\tilde{S}_{1,{k_1,\fa}} \ll n \log^{3(s+1)} n \ll  n^{\fs-3/4}$  $(\fs \geq 2)$, and
   \eqref{Le18-6a} follows.
Thus,  Lemma 15 is proved. \qed \\ \\
\subsection{\textcolor{blue}{Final result }}
{\bf Lemma 16.} {\it We have  }
\begin{equation}\label{Le13-01}
E_s (\sD_{T_{\fs,\sh,1}}(Q,N))\ll n^{\sh s/2  -1/5 } \quad \for \quad s \geq 3 \;\; {\rm or}\;\; s=\fs=2
\end{equation}
and
\begin{equation}\label{Le13-10}
E_s (\sD_{T_{\fs,\sh,1}}(Q,N))\ll n^{\sh s/2  } \quad \for \quad s =2 \;\; \ad \;\; \fs=1.
\end{equation}
\\
{\bf Proof.}   In view of Lemma 4 and \eqref{Par3.2-3}, in order to prove \eqref{Le13-01} it is enough to verify that
\begin{equation}  \nonumber 
\DD_{T_{\fs}, a, \lambda,\tau}  \ll    n^{\sh s/2 -1/5 }.
\end{equation}
By \eqref{Le3-5}, it is enough to verify that
\begin{equation}  \label{Cor3}
\SS_{m, \tau}    \ll n^{\sh s/2-1/4}.
\end{equation}
Taking $\fa(\rho^{(k)}, m^{(k)},k ) \equiv 1$, we obtain from \eqref{Le6-1} and \eqref{Lem11-0}, that
\begin{equation} \nonumber  
\SS_{m, \tau}    \ll \SS_{m, \tau,1}.
\end{equation}
In the following, we will use Lemma 10, Lemma 12  and Lemma 14 with $\fa(\rho^{(k)}, m^{(k)},k ) \equiv$ 1.\\
Let
\begin{equation} \nonumber
\bfc = \begin{cases}
     1 , & {\rm if} \; \fs=1,  \; s=2 ,  \\
     2 , & {\rm if}  \; s >\fs, \; s \geq 3,   \\
     3, & {\rm if} \;  s=\fs \geq 2, \;  \exists k \; \with \;\lambda_{1,k}=1, \\
     4 , & {\rm if} \;  s=\fs \geq 2,  \; \forall k \; \lambda_{1,k}\geq 2, \; \exists k \;
                      \with \; \lambda_{1,k}\geq 3, \\
     5 , & {\rm if} \; s=\fs \geq 2,  \; \forall k \; \lambda_{1,k} =2, \; \exists i_1,k \;
                      \with \; \lambda_{i_1,k}\neq 2, \\
     6 , & {\rm if} \; s=\fs \geq 2,  \; \forall i,k \; \lambda_{i,k}= 2.
   \end{cases}
\end{equation}

Let's consider the case $\bfc =1$. By Lemma 6, we get \eqref{Le13-10}.

Let's consider the case $\bfc =2$. By Lemma 10, we get  \eqref{Cor3}.

Let's consider the case $\bfc =3$. By Lemma 12, we get \eqref{Cor3}.

Let's consider the case $\bfc =4$. In view of \eqref{Par3-5}, $a_1 <\sh/2$.

  $\;\;\;  \qquad \qquad \qquad \qquad \qquad \qquad $ By Lemma 10, we get  \eqref{Cor3}.

Let's consider the case $\bfc =5$. By Lemma 14, we get \eqref{Cor3}.

Let's consider the case $\bfc =6$. By Lemma 15, we get \eqref{Cor3}.\\
Hence, Lemma 16  is proved. \qed \\ \\
\section{\textcolor{blue}{Completion  of proofs of Theorems}}

\subsection{\textcolor{blue}{Proof of Theorem 1}}
%

 Bearing in mind the monotony of the $L_p$ norm, we get that it is enough to consider only the case of $p=\sh=2q$.

Applying  \eqref{Beg-28}  and   Minkowski's inequality, we derive :
\begin{equation} \label{Le13-0}
\left\| D(\bx, (H(k))_{k=Q}^{Q+N-1}  )   \right\|_{s,2q}  \ll \sum_{\fs=1}^{s}  \sum_{T_{\fs} \subseteq\{1,...,s\} }
    \left\| \fD_{T_{\fs}}(Q,N) \right\|_{s,2q} + \log^s n .
\end{equation}
In view of \eqref{Lem2-7} and \eqref{Le3-2}, we have
\begin{equation}  \nonumber
  \fD_{T_{\fs}}^{2q}(Q,N)  = \sD_{T_{\fs,2q,1}}(Q,N) + \sD_{T_{\fs,2q,2}}(Q,N) \quad \ad \quad
    \sD_{T_{\fs,2q,2}}(Q,N) \ll n^{qs}.
\end{equation}
Using  Lemma 16, we get
\begin{equation}  \nonumber
      E_s(\sD_{T_{\fs,\sh,1}}(Q,N)) \ll n^{qs}.
\end{equation}
Hence
\begin{equation}  \label{Le13-2}
   E_s(\fD_{T_{\fs}}^{2q}(Q,N))    \ll n^{qs} \quad \ad \quad  \left\| \fD_{T_{\fs}}(Q,N) \right\|_{s,2q} \ll n^{s/2} .
\end{equation}
From \eqref{Le13-0}, we obtain
\begin{equation}  \nonumber
   \left\| D(\bx, (H(k))_{k=Q}^{Q+N-1}  )   \right\|_{s,2q}  \ll    n^{ s/2   }.
\end{equation}
By \eqref{In2}, Theorem 1 is proved. \qed \\ \\
\subsection{\textcolor{blue}{Proof of Theorem 2}}

  The assertion of Theorem 2 follows essentially from Lemma 19. To prove Lemma 19,
we need firstly to compute the main value of the product of functions $\varphi_{ \br,0,N,m}$ (see \eqref{Le2-1}):
\begin{equation} \label{Th2-00}
 \ddot{\gamma}^{(q)}_{\br,\bm} :=  \int_0^1 \prod_{j=1}^{q}
   |\varphi_{ \br_j,0,[Nx_{s+1}],m_j}|^2  d x_{s+1}, \quad \quad\varphi_{ \br,0,N,m}  =\frac{e(m  N/P_{\br})-1}{P_{\br}(e(m/P_{\br})-1)}.
\end{equation}
 Taking into account \eqref{Beg-2}, \eqref{Le2-1} and that $ 2|x|/\pi \leq |\sin(x)| \leq |x|$, we obtain
\begin{equation} \nonumber
|\varphi_{ \br,0,N,m}|^2  \leq \left| \frac{\sin(\pi m  N/P_{\br})}{P_{\br} \sin(\pi m/P_{\br})}\right|^2 \leq  \frac{\sin^2(\pi m  N/P_{\br})}{ \bar{m}^2}
\end{equation}
and
\begin{equation}\label{Th2-01}
  |\ddot{\gamma}^{(q)}_{\br,\bm}|
  \leq \max_{0 \leq x_{s+1} \leq 1} \prod_{j=1}^{q} \frac{\sin^2\big( \frac{\pi [N x_{s+1}] |m_j|}{P_{\br_j}}\big) }{\bar{m}^2_{j}}
  \leq \prod_{j=1}^{q} \frac{ \min \big(1, \big(\frac{ \pi N |m_j|}{P_{\br_j}} \big)^2) }{\bar{m}^2_{j}}.
\end{equation}\\
In the following Lemma 17, we find a simple expression that approximates $ \ddot{\gamma}^{(q)}_{\br,\bm}$.
To prove Lemma 17, we use  inequality \eqref{Beg-2a}.\\ \\
{\bf Lemma 17.} {\it Let $0 < |m_j| \leq n^{4 \sh s}, \; j=1,...,q$,  and let
\begin{equation}  \label{Th2-002}
    \fU: = \{\br\in [V_1,n]^s \;|\; P_{\br} \leq 2^{n +\log_2^3 n}\}. 
\end{equation}
 Then}
\begin{multline}   \nonumber
\imath:  = \sum_{\substack{\br_j \in \fU  \\ j \in [1,q]}} |   \ddot{\gamma}^{(q)}_{\br,\bm}-
       \theta^{(q)}_{\br,\bm}|  \ll  \frac{ n^{qs-1} \log^3 n}{\bar{m}_1^2 ... \bar{m}_{q}^2}   , \;\;\;\; \with \\
        \theta^{(q)}_{\br,\bm}:=
    \prod_{j=1}^{q}   \frac{2}{|P_{\br_j}(1-e(m_j/P_{\br_j}))|^2}, \;s \geq 3,
\end{multline}
\begin{multline}  \label{Th2-02}
\sum_{\substack{\br_j \in \fU  \\  j \in [1,q]}} |   \dddot{\gamma}^{(q)}_{\br,\bm}-
       \theta^{(q)}_{\br,\bm}|  \ll  \frac{ n^{qs-1} \log^3 n}{\bar{m}_1^2 ... \bar{m}_{q}^2}   , \quad  \with \\
 \dddot{\gamma}^{(q)}_{\br,\bm} =  \int_0^1 \prod_{j=1}^{q}
   |\varphi_{ \br_j,-[Nx_{3}],2[Nx_{3}],m_j}|^2  d x_{3},\; s=2.
\end{multline} \\
{\bf Proof.} We will prove the first statement. The proof of the second statement is similar.
By  \eqref{Th2-00} and \eqref{Th2-02},  we obtain
\begin{equation}   \nonumber
  \ddot{\gamma}^{(q)}_{\br,\bm} =   \frac{1}{N} \sum_{k=0}^{N-1} \prod_{j=1}^{q}
   |\varphi_{ \br,0,k,m_j}|^2 =  \frac{1}{N} \sum_{k=0}^{N-1} \prod_{j=1}^{q}
 \frac{2-2 \cos(2 \pi m_jk/P_{\br_j} )}{|P_{\br_j}(1-e(m_j/P_{\br_j}))|^2} = \theta^{(q)}_{\br,\bm}
  (1+\epsilon g_{\br,\bm}),
\end{equation}
\begin{equation}  \label{Th2-50}
 \with \quad  g_{\br,\bm}:= \sum_{J \subseteq \{1,...,q \}, J \neq \emptyset} \Big|\frac{1}{N} \sum_{k=0}^{N-1}  \prod_{j \in J}  \cos \Big(2 \pi k  m_j/P_{\br_j}  \Big)\Big|, \;\;\; |\epsilon| \leq 2 .
\end{equation}
It is easy to see that
\begin{equation}  \nonumber
 g_{\br,\bm} \leq  \sum_{J \subseteq \{1,...,q \}, J \neq \emptyset} \;
  \sum_{ \nu_j \in \{ -1,1\}, j \in J} |\cX_{\br}|, \quad \cX_{\br}=   \frac{1}{N} \sum_{k=0}^{N-1}  e \Big( k  \sum_{j \in J} \nu_j m_j/P_{\br_j} \Big).
\end{equation}
Using   \eqref{Beg-2a},  we have
\begin{equation}   \label{Th2-1005a}    
  |\cX_{\br}| \leq
   \min \Big(1, \frac{1}{2 N  \llangle  Y^{'}_{\br}  \rrangle }    \Big), \quad \with \quad
     Y^{'}_{\br} = \sum_{j \in J} \nu_j m_j/P_{\br_j}.
\end{equation}
Similarly to \eqref{Beg-2} and  \eqref{Le2-1}, we obtain
\begin{equation}   \label{Th2-1005}    
    \theta^{(q)}_{\br,\bm}    \leq     \frac{ 1}{\bar{m}_1^2 ... \bar{m}_{q}^2}.
\end{equation}
Let $ \sL = \card(J)$, $J =(\fj_1,...,\fj_{\sL})$, $\fj_i < \fj_{i+1}\; (i=1,2,...)$,  $\cJ= \fj_{\sL}$.\\
By  \eqref{Th2-02} - \eqref{Th2-1005}, we get
\begin{multline} \label{Th2-70}
 \imath \leq \sum_{\substack{J \subseteq \{1,...,q \}\\ J \neq \emptyset}} \;
  \sum_{\substack{\nu_i \{ -1,1\} \\ i \in [1,\sL]  }} \;\;  \sum_{\br_j \in \fU,\;
  j \in \{ 1, \ldots, q \}\setminus \{\cJ  \}}  \; \frac{ \fF_{\cJ}}{\bar{m}_1^2 ... \bar{m}_{q}^2}, \quad \with\\
  \fF_{\cJ}: =   \sum_{ \br_{\cJ} \in \fU}   \min \Big(1, \frac{1}{2 N  \llangle  Y^{'}_{\br}  \rrangle }    \Big), \quad Y^{'}_{\br} = \sum_{j \in J} \nu_j m_j/P_{\br_j}.
\end{multline}
From \eqref{Th2-002}, we have that in order to prove \eqref{Th2-02}, it is enough to verify that
\begin{equation}\label{Th2-71}
      \fF_{\cJ}   \ll n^{s-1} \log^3 n .
\end{equation}
Now we will prove \eqref{Th2-71}:

We fix $m_j, \br_j$, for all $j \in J \setminus \{\cJ \} $. Then  $  Y^{'}_{\br}  =f' + \nu_{\cJ}  m_{\cJ}/P_{\br_{\cJ}}$  
for some $f'$. Let $f \equiv f' \mod 1$ and $f \in (-1/2, 1/2]$. Taking into account that the function $\llangle  z  \rrangle $ has period one, we get  $\llangle  Y_{\br}  \rrangle =\llangle  Y^{'}_{\br}  \rrangle $ with  $  Y_{\br}   =f +  \nu_{\cJ} m_{\cJ}/P_{\br_{\cJ}}$.
We see that
\begin{multline}\label{Th2-72}
   \fF_{\cJ} = \fF_{\cJ,1} + \cdots +  \fF_{\cJ,5}, \quad \where \\
       \fF_{\cJ,i}  = \sum_{\br_{\cJ} \in [1,n]^s, P_{\br_{\cJ}} \leq 2^{n +\log_2^3 n}}  \min \left(1, \frac{1}{ 2 N\llangle Y_{\br} \rrangle} \right)    \b1(b_{\br} =i),\;\;   m_{\cJ} \neq 0,
\end{multline}
with
\begin{equation}\label{Th2-73}
 b_{\br} = \begin{cases}
     1 , & {\rm if} \; \llangle Y_{\br} \rrangle  \geq n^s/N,  \\
     2 , & {\rm if} \;  |m_{\cJ}|/P_{\br_{\cJ}} < 4n^s/N,  \\
  3 , & {\rm if} \;  |m_{\cJ}|/P_{\br_{\cJ}} \geq 1/4, \\
    4 , & {\rm if} \; \llangle Y_{\br} \rrangle  < n^s/N, \;1/4 > |m_{\cJ}|/P_{\br_{\cJ}} \geq 4n^s/N,\; |f| >  2n^s/N,  \\
    5 , & {\rm if} \; \llangle Y_{\br} \rrangle  < n^s/N, \;1/4 >  |m_{\cJ}|/P_{\br_{\cJ}} \geq 4n^s/N ,\; |f| \leq  2n^s/N .
  \end{cases}
\end{equation}

Let's consider the case $ b_{\br}=1$.
By \eqref{Th2-72} and \eqref{Th2-73},  we obtain
\begin{equation}   \nonumber  
    \fF_{\cJ,1}  \leq \sum_{\br_{\cJ} \in [V_1,n]^s } \frac{1}{N} \cdot \frac{N}{n^s} \leq 1.
\end{equation}

Let's consider the case $ b_{\br}=2$.
By \eqref{Th2-72} and \eqref{Th2-73},  we derive
\begin{equation}    \nonumber  
   1 \leq |m_{\cJ}|, \;\; |m_{\cJ}|/P_{\br_{\cJ}} < 4n^s/N, \quad 0.25N /n^s  \leq 0.25N|m_{\cJ}| /n^s  \leq P_{\br_{\cJ}} \leq   2^{n +\log_2^3 n},
\end{equation}
\begin{equation}    \nonumber  
  n-3 -s\log_2 n \leq \sum_{i=1}^s r_{\cJ,i} \log_2 p_i  \leq  n + \log_2^3 n \quad (n =[\log_2 N] +1, \; P_{\br}=p_1^{r_1} \cdots p_s^{r_s}).
\end{equation}
It is easy to verify that the number of solutions of this inequality is equal to $O(n^{s-1} \log_2^3 n)$ and $\fF_{\cJ,2} =O(n^{s-1} \log_2^3 n)$.\\

Let's consider the case $ b_{\br}=3$.
By \eqref{Th2-73},  we get
\begin{equation*}
 p_1^{r_{1,\cJ}} \cdots   p_{\sh}^{r_{\sh,\cJ}} =   P_{\br_{\cJ}} \leq  4 | m_{\cJ}| \leq 4 n^{4\sh s}.
\end{equation*}
We see that the number of solutions of this inequality is equal to $ O(\log_2^s n)$ and $\fF_{\cJ,3} =O( \log_2^s n)$.\\

Let's consider the case $ b_{\br}=4$. We have
\begin{equation*}
  Y_{\br} = \nu_{\cJ} m_{\cJ}/P_{\br_{\cJ}}  +f \;\; \with \;\;   2n^s/N < |f| \leq 1/2, \quad \llangle Y_{\br} \rrangle  < n^s/N \;\;  \ad \;\; |m_{\cJ}|/P_{\br_{\cJ}} <1/4.
\end{equation*}
Hence
\begin{equation*}
 |Y_{\br} | \leq  | m_{\cJ}/P_{\br_{\cJ}}| +|f| \leq 3/4.
\end{equation*}
Bearing in mind that
\begin{equation*}
 \llangle Y_{\br} \rrangle  < n^s/N , \quad n^s/N \to 0  \quad \for \quad N \to \infty,\quad {\rm we\; have}  \quad |Y_{\br} | = \llangle Y_{\br} \rrangle .
\end{equation*}
Taking into account that $ \nu_{\cJ} m_{\cJ}/P_{\br_{\cJ}} = -f +Y_{\br}$, we derive
\begin{equation*}
  \nu_{\cJ} m_{\cJ}/P_{\br_{\cJ}} \in [-f - n^s/N , -f + n^s/N] \;\;
 \ad \;\; |m_{\cJ}|/P_{\br_{\cJ}} \in [|f|/2 , 2 |f|] \; \for \;
 |f|>2n^s/N.
\end{equation*}
 Therefore
\begin{equation}    \label{Th2-60}
  \log_2 |m_{\cJ}| - \sum_{i=1}^s r_{\cJ,i} \log_2 p_i \in [\log_2 |f| -1, \log_2 |f| +1].
\end{equation}
Bearing in mind that $ | \log_2 |f| | \leq 2 \log_2 N  \leq 2n $, we get that  the number of solutions of  \eqref{Th2-60}
is equal to $O(n^{s-1})$. By \eqref{Th2-72}, we have $  \fF_{\cJ,4}   =O(n^{s-1})$.  \\

Let's consider the case $ b_{\br}=5$. Taking into account that
\begin{equation*}
  |f| \leq 2n^s/N, \;\;\; \; n^s/N \to 0  \quad \ad \quad |m_{\cJ}|/P_{\br_{\cJ}} \leq 1/4,
\end{equation*}
 we get that
$|f +  \nu_{\cJ} m_{\cJ}/P_{\br_{\cJ}}|=|Y_{\br} | \leq 3/8 $.  Hence  $|Y_{\br} | =\llangle Y_{\br} \rrangle < n^s/N$.\\
By \eqref{Th2-73}, we obtain
\begin{equation*}
    4n^s/N \leq |m_{\cJ}|/P_{\br_{\cJ}} \leq  |Y_{\br} |+|f| < 3n^s/N  .
\end{equation*}
We have a contradiction and $ \fF_{\cJ,5}=0$.
Hence
\begin{equation}    \nonumber
    \fF_{\cJ}= \fF_{\cJ,1} + \cdots +  \fF_{\cJ,5} =O(n^{s-1} \log^3 n).
\end{equation}
 Therefore, \eqref{Th2-71} and  Lemma 17 are proved. \qed \\  \\

In the following Lemma 18, we consider the expectation for the case $\fs=1$ and $s=2$. The estimate of the $s$-dimensional expectation in \eqref{Le13-10} is not enough for our purpose. Thereby, we will consider $s+1$-dimensional expectation to prove Theorem 2.  The main tool is  inequality \eqref{Beg-2a}. \\ \\
{\bf Lemma 18.} {\it  Let $\fs=1,\; s=2$. Then}
\begin{equation}\label{Lem8-1}
   E_{s+1}(\sD_{T_{\fs},\sh,1}(-[Nx_{s+1}],2[Nx_{s+1}] ))  = O( n^{ \sh  s /2 -1/2 }).
\end{equation} \\
{\bf Proof.}
 Let's consider the case $T_{\fs} =\{ 1\}$. The proof for the case $T_{\fs} =\{ 2\}$ is similar.
By Lemma 4, we have
\begin{multline} \label{Lem8-1a}
  |E_{s+1} (\sD_{T_{\fs,\sh,1}}(-[Nx_{s+1}], 2[Nx_{s+1}])) | \ll   1 + \tilde{\cD}_{T_{\fs}}, \;\with   \\
   \tilde{\cD}_{T_{\fs}}=
\sum_{\substack{m_j \in I^{*}_{n^{4 \sh s}}, \; |\fm_{i,j}| \leq n^{4 \sh s}\\  i=1,...,s,\; j=1,...,\sh}}\;
 \sum_{\substack{\br_j \in \UU_{ T_{\fs}}\\ j=1,...,\sh}}\;
\frac{\varpi_{\br,\bm,1} \zeta_{\br,\fm}\; |\hat{\gamma}^{(\sh)}_{\br,\bm}| }{\bar{m}_1 \cdots \bar{m}_{\sh}}  \;
\;
\prod_{i=1}^s  \prod_{j=1}^{\sh} \frac{1}{\bar{\fm}_{i, j}}
       \delta_{p_i^{r^{+}_i}} (\tilde{m}_{i} ) .
\end{multline}
Similarly to \eqref{Par3.2-3} and \eqref{Le6-1a}, we get
\begin{equation} \label{Lem8-1b}
\tilde{\cD}_{T_{\fs}}  =  \sum_{a_1,...,a_s =1}^{\sh} \;\;
    \sum_{1 \leq \lambda_{i,k} \leq \sh, i=1,...s , \; \lambda_{i,1}+...+\lambda_{i,a_i}=\sh }\;\;
 \sum_{ \tau_i \in \Xi_{\sh},  i=1,...,s}  \tilde{\DD}_{T_{\fs}, a, \lambda,\tau}
\end{equation}
with
\begin{equation} \label{Lem8-2}
\tilde{\DD}_{T_{\fs}, a, \lambda,\tau}   = \sum_{m_j \in I^{*}_{n^{4 \sh s}},\; |\fm_{i,j}| \leq n^{4 \sh s},   i=1,...,s, j=1,...,\sh} \;\;\;
 \prod_{j=1}^{\sh} \frac{1}{\bar{m}_j}
     \prod_{i =1}^s\frac{1}{\bar{\fm}_{i,j}}
      \SS_{m,  \tau},
\end{equation}
where
\begin{equation}  \nonumber
           \SS_{m, \tau} \leq       \sum_{\brho_{j} \in \UU_{T_{\fs}}, j\in[1,\sh]} \;   \varpi_{\brho,\bm,1} \;
    \chi_{\brho} \;   \zeta_{\brho,\fm} \; |\hat{\gamma}^{(\sh)}_{\br,\bm}|\; \delta_{p_1^{\rho^{+}_{1}}}(\tilde{m}_{1} ).
\end{equation}
Bearing in mind \eqref{Lem8-1a} and \eqref{Lem8-1b}, we see that in order to prove   \eqref{Lem8-1}, it is enough verify that
\begin{equation}  \label{Lem8-3}
      \DD_{T_{\fs}, a, \lambda,\tau} \ll  n^{ \sh  s /2 -1/2 }.
\end{equation}

 Let's consider the case $a_1 = \sh$. From  \eqref{Par3-5}, we get $\lambda_{1,j}=1$ for $j=1,...,\sh$.\\
By \eqref{Le14-52},  we have that
\begin{equation}  \label{Lem8-4}
 m_{j}+\fm_{1,j} p_2^{\rho_{2,j}} =0, \quad \for \quad j=1,...,\sh.
\end{equation}
Similarly to  \eqref{Le14-7}, we derive
\begin{multline} \label{Lem8-4a}
\DD_{T_{1}, a, \lambda,\tau}
\leq       \sum_{\substack{\rho_{1,j} \in [1,n],  j\in[1,\sh] \\ \rho_{2,j} \in [1,V_1]   }} \;
 \sum_{\substack{ |\fm_{i,j}| \leq n^{4 \sh s}\\  |\fm_{i,j}| \leq  p_i^{\rho^{+}_i}, \; i=1,2, \; j=1,...,\sh}}
 \prod_{j=1}^{\sh}  \frac{|\hat{\gamma}^{(\sh)}_{\br,\bm}|}{(\bar{\fm}_{1,j})^2 p_2^{\rho_{2,j}}\bar{\fm}_{2,j}}  \\
%
\leq      \sum_{\substack{\rho_{2,j} \in [1,n]\\  j \in [1,\sh]}}    \;
 \sum_{\substack{ |\fm_{i,j}| \leq n^{4 \sh s}\\  |\fm_{2,j}| \leq  p_2^{\rho^{+}_2}, \; j=1,...,\sh}}
 \prod_{j=1}^{\sh}  \frac{1}{(\bar{\fm}_{1,j})^2 p_2^{\rho_{2,j}}\bar{\fm}_{2,j}} \;  \Gamma^{\prime}   \quad \with \\
%
\Gamma^{\prime} =\max_{\fm_{i,j},\rho_{i,j} } \Gamma, \quad \where \quad \Gamma =  \sum_{\substack{\rho_{1,j} \in [1,n]\\  j \in [1,\sh]}}    \; |\hat{\gamma}^{(\sh)}_{\br,\bm}| .
\end{multline}
Hence
\begin{equation} \nonumber
\DD_{T_{1}, a, \lambda,\tau}
\ll     \Gamma^{\prime}   \sum_{\substack{\rho_{2,j} \in [1,n]\\  j \in [1,\sh]}}    \;
 \sum_{\substack{  |\fm_{2,j}| \leq  p_2^{\rho^{+}_2}\\ j=1,...,\sh}}
 \prod_{j=1}^{\sh}  \frac{1}{ p_2^{\rho_{2,j}}\bar{\fm}_{2,j}}
\ll     \Gamma^{\prime}   \sum_{\rho^{+}_2 \in [1,n]}    \;
   \frac{ (\rho^{+}_2)^{\sh} }{ p_2^{\rho^{+}_2}} \ll \Gamma^{\prime}.
\end{equation}
We have that in order to prove   \eqref{Lem8-3}, it is enough verify that
\begin{equation}  \label{Lem8-5}
    \Gamma^{\prime}  \ll  n^{ \sh  s /2 -1/2 }.
\end{equation}

Now we will prove \eqref{Lem8-5}. Let's consider $ \hat{\gamma}^{(\sh)}_{\br,\bm} $.\\
In view of Lemma 1 and Lemma 2, we have
\begin{equation*} 
  \hat{\gamma}^{(\sh)}_{\brho,\bm}= \bar{m}_1 \cdots \bar{m}_{\sh} \int_0^1 \prod_{j=1}^{\sh} \varphi_{ \brho_j,-[Nx_{s+1}], 2[Nx_{s+1}],m_j} d x_{s+1},
\end{equation*}
with
\begin{equation}  \nonumber
  \varphi_{ \brho_j,-[Nx_{3}],2[Nx_{3}],m_j}  =\frac{e(m_j[Nx_{3}])/P_{\brho_j}) - e(-m_j[Nx_{3}]/P_{\brho_j} )}{P_{\brho_j}(e(m_j/P_{\brho_j})-1)}.
\end{equation}
Applying Lemma 1,  we get
\begin{equation}  \label{Lem8-5a}
    |\hat{\gamma}^{(\sh)}_{\br,\bm}|  \leq 1.
\end{equation}
We have
\begin{equation*}
 \hat{\gamma}^{(\sh)}_{\brho,\bm} =   \frac{1}{N} \sum_{k=0}^{N-1} \prod_{j=1}^{\sh} \bar{m}_j
    \frac{e(m_j k/P_{\brho_j}) - e(-m_j k/P_{\brho_j} )}{P_{\brho_j}(e(m_j/P_{\brho_j})-1)}
\end{equation*}
\begin{equation*}
 =   \frac{1}{N} \sum_{k=0}^{N-1} \prod_{j=1}^{\sh} \bar{m}_j
    \frac{2\sqrt{-1} \sin(2\pi m_j k/P_{\brho_j})}{P_{\brho_j}(e(m_j/P_{\brho_j})-1)}
 =   \sum_{\nu_1,...,\nu_{2q} \in \{-1,1\} } \nu_1 \nu_2 \cdots \nu_{\sh} \AA_{\nu},
\end{equation*}
with
\begin{equation*}
  \AA_{\nu} = \XX_{\rho} \; \prod_{j=1}^{\sh}
    \frac{ \bar{m}_j }{P_{\brho_j}(e(m_j/P_{\brho_j})-1)} , \quad
     \XX_{\rho}=  \frac{1}{N} \sum_{k=0}^{N-1}  e \Big( k  \sum_{j \in J} \nu_j m_j/P_{\brho_j} \Big)
\end{equation*}
and
\begin{equation}   \label{Lem8-7}   
  \hat{\gamma}^{(\sh)}_{\brho,\bm} =    \XX^{'}_{\rho} \; \prod_{j=1}^{\sh}
    \frac{ \bar{m}_j 2\sqrt{-1}  }{P_{\brho_j}(e(m_j/P_{\brho_j})-1)} , \quad  \XX^{'}_{\rho}=
    \frac{1}{N} \sum_{k=0}^{N-1} \prod_{j=1}^{\sh} \sin \big(2 \pi k   m_j/P_{\brho_j} \big).
\end{equation}
By \eqref{Beg-2a},  we have
\begin{equation*}     
      \frac{ 2\bar{m}_j }{P_{\brho_j}|e(m_j/P_{\brho_j})-1|}  \leq 1.
\end{equation*}
Using  Lemma 1, we get
\begin{equation}   \label{Lem8-8}    
  | \AA_{\nu}| \leq 2^{\sh}|\XX_{\rho}|, \quad  |\XX_{\rho}| \leq
    \min \Big(1, \frac{1}{2 N  \llangle  \YY_{\brho}  \rrangle }    \Big), \quad \with \quad
     \YY_{\brho} = \sum_{j \in [1,\sh]} \nu_j m_j/P_{\brho_j}.
\end{equation}
In view of \eqref{Lem8-7},  we derive
\begin{equation}   \label{Lem8-10}    
 | \hat{\gamma}^{(\sh)}_{\brho,\bm}| \leq \max_{\nu_1,...,\nu_{2q} \in \{-1,1\} } 2^{\sh}\; \min \Big(1, \frac{1}{2 N  \llangle  \YY_{\rho}  \rrangle }    \Big),\quad
 \ad   \quad  | \hat{\gamma}^{(\sh)}_{\brho,\bm}| \leq  |\XX^{'}_{\rho}|.
\end{equation} \\
We see that
\begin{multline} \nonumber 
   \Gamma = \Gamma_{1} + \cdots +  \Gamma_{5}, \quad \where \\
       \Gamma_{i}  = \sum_{\br_{\cJ} \in [1,n]^s, P_{\br_{\cJ}} \leq 2^{n +\log_2^3 n}}  \min \left(1, \frac{1}{ 2 N\llangle Y_{\br} \rrangle} \right)    \b1(\tilde{b}_{\br} =i),\;\;   m_{\cJ} \neq 0,
\end{multline}
with
\begin{equation}  \nonumber 
 \tilde{b}_{\rho} = \begin{cases}
     1 , & {\rm if} \; a_1 = \sh,\; \rho_1^{+} \geq n\log_{p_1}2 + 100s \sh \log_{p_1} n
     ,   \\
     2 , & {\rm if} \; a_1 = \sh,\; \rho_1^{+} \in [ n\log_{p_1}2 - 100s \sh \log_{p_1} n ,n\log_{p_1}2 + 100s \sh \log_{p_1} n  ],  \\
  3 , & {\rm if} \; a_1 = \sh,\;  \rho_1^{+} \leq   n\log_{p_1}2 - 100s \sh \log_{p_1} n,\; \; \YY_{\br} \neq 0, \\
    4 , & {\rm if} \; \; a_1 = \sh, \;  \rho_1^{+} \leq   n\log_{p_1}2 - 100s \sh \log_{p_1} n,\; \; \YY_{\br} = 0,     \\
     5 , & {\rm if} \; a_1 < \sh.
  \end{cases}
\end{equation}
Let's consider the case $ \tilde{b}_{\rho}=1$:\\
By \eqref{Lem8-7} and  \eqref{Lem8-10},  we get
\begin{equation} \nonumber
 | \hat{\gamma}^{(\sh)}_{\brho,\bm}| \leq  |\XX^{'}_{\rho}| \leq \max_{1 \leq k \leq N}|\sin (2 \pi   k\fm_{1,\sh}/ p_1^{\rho_1^{+}} ) | \leq 2 \pi   N |\fm_{1,\sh}|/ p_1^{\rho_1^{+}} \ll n^{-s\sh}
\end{equation}
From \eqref{Lem8-4a}, we get \eqref{Lem8-5}.\\ \\
Let's consider the case $ \tilde{b}_{\rho}=2$:\\
  Taking into account \eqref{Lem8-4a} and \eqref{Lem8-5a}, we obtain
\begin{multline}   \nonumber 
\Gamma^{\prime} =\max_{\fm_{i,j},\rho_{i,j} }  \sum_{\substack{\rho_{1,j} \in [1,n]\\  j \in [1,\sh]}}    \; |\hat{\gamma}^{(\sh)}_{\br,\bm}|
  \b1(\rho_1^{+} \in [ n\log_{p_1}2 - 100s \sh \log_{p_1} n ,n\log_{p_1}2 + 100s \sh \log_{p_1} n  ]) \\
   \ll n^{\sh s -1} \log_2 n
\end{multline}
and \eqref{Lem8-5} follows. \\  \\
Let's consider the case $ \tilde{b}_{\rho}=3$:\\
In view of \eqref{Lem8-8}, \eqref{Lem8-10} and \eqref{Lem8-4a}, we have
\begin{equation}  \nonumber 
|\XX_{\brho}| \leq
   \frac{p_1^{\rho_1^{+}}}{2 N  }   \ll n^{-10}, \quad  |\hat{\gamma}^{(\sh)}_{\br,\bm}| \ll n^{-10},
   \;\; \ad \;\;  \Gamma_1  \ll n^{s\sh/2 -10}
\end{equation}
and \eqref{Lem8-5} follows. \\  \\
Let's consider the case $ \tilde{b}_{\rho}=4$:
Applying  \eqref{Lem8-8}, we obtain
\begin{equation}   \label{Lem8-9}   
       \nu_1 \fm_{1,1} p_1^{\rho_1^{+}-   \rho_{1,1}} = -\sum_{j \in [2,\sh]} \nu_j \fm_{1,j} p_1^{\rho_1^{+}-   \rho_{1,j}}.
\end{equation}
By \eqref{Lem8-2} and \eqref{Lem8-4}, $\fm_{1,1}  \neq 0$. We see, that for fixed $\rho_{1,\sh}, \rho_{1,\sh-1},...,\rho_{1,2}$ there is at most one solution in the variable $\rho_{1,1}$ of equation
\eqref{Lem8-9}.
%
According to \eqref{Lem8-4a} and \eqref{Lem8-10}, statement \eqref{Lem8-5}    is proved.\\ \\
Let's consider the case $ \tilde{b}_{\rho}=5$:\\
In this case inequality \eqref{Lem8-5a} remains true.
Using   \eqref{Lem8-2},  we have
\begin{equation}  \nonumber
           \SS_{m, \tau} \leq       \sum_{\substack{1 \leq \rho_{1,j} \leq n,\; 1 \leq \rho_{1,j}  \\ j\in[1,\sh]}} \;
    \chi_{\brho}  \leq
            V_1^{\sh}  \max_{\rho_2}  \sum_{\substack{1 \leq \rho_{1,j} \leq n  \\ j\in[1,\sh]}} \;
    \chi_{\brho}  \leq   V_1^{2\sh} n^{a_1} \leq n^{\sh-1/3}.
\end{equation}
Now \eqref{Lem8-3} follows from \eqref{Lem8-1b} and\eqref{Lem8-2}.

Hence Lemma 18 is proved. \qed \\ \\

The next lemma is the main lemma in this section. Its proof is based on inequality \eqref{Beg-2a}, Lemma 17 ,  rearrangements of domains of summations and multiple changes  of  orders of  summations.

The main point of the proof is the transition from the condition \\ $  m_{\eta(2k-1)}/P_{\br_{\eta(2k-1)}}  =- m_{\eta(2k)}/P_{\br_{\eta(2k)}} $ to the condition $  \br_{\eta(2k-1)}  = \br_{\eta(2k)} $. The next point of the proof is the transition from the set $\Xi_{\sh}$  (the set of all permutations of the set $ \{1, ..., \sh \} $) to the set $  \tilde{\Xi}_q $, where
\begin{multline*}
  \tilde{\Xi}_q = \{ (\sigma_1,\sigma_2),  \sigma_i: \{1,...,q\} \to \{1,...,2q \}, \;i=1,2\;\; |\;\;
  \{\sigma_1(1),\sigma_2(1),...,\\ \sigma_1(q),  \sigma_2(q)\}
   = \{1,...,2q\},\;  \;\sigma_1(k)<\sigma_2(k) \; \forall k  \}, \quad \with \\
    \card(\tilde{\Xi}_q )   =  \binom{2q}{2}  \binom{2q-2}{2}  \cdots  \binom{2}{2}
    = \frac{(2q)!}{ 2^{q}}.
\end{multline*}
From Lemma 19, we will easily obtain the assertions of Lemma A and Theorem 2 in what follows.
\\  \\
{\bf Lemma 19.}   {\it With notations as above}
\begin{multline} \label{Th2-30}
  \varsigma_{2q} =   E_{s+1} (  \sD_{T_{s,2q,2}}(0,[Nx_{s+1}])  )
  = \frac{(2q)!}{ 2^{q} q!}\Big(E_{s+1} ( \sD_{T_{s,2,2}}( 0,[Nx_{s+1}])  ) \Big)^q \\+ O(n^{qs-9/10}),
       \quad s \geq 2,
\end{multline}
and
\begin{multline} \label{Th2-31}
E_{s+1} (\sD_{T_{s,2q,2}}(-[Nx_{s+1}], 2[Nx_{s+1}]))\\ =
  \frac{(2q)!}{ 2^{q} q!}\Big(E_{s+1} (\sD_{T_{s,2,2}}(-[Nx_{s+1}], 2[Nx_{s+1}])) \Big)^q+ O(n^{qs-9/10}), \; \;\; s=2.
\end{multline}  \\
{\bf Proof.} We will prove the first statement. The proof of the second statement is similar.
Using \eqref{Th2-30}, \eqref{Beg-26}, \eqref{Le2-1} and \eqref{Lem2-6},   we obtain
\begin{multline}  \label{Th2-30a1}
\varsigma_{2q}	=\sum_{\substack{r_{i,j} \in [V_1,n], i \in[1,s]  \\ j \in [1,2q]}} \;
\sum_{\substack{m_j \in I^{*}_{P_{\br_j}} \\ j \in [1,2q]}}  \varpi_{\br,\bm,2}  \gamma^{(2q)}_{\br,\bm}  \beta^{(2q)}_{\br,\bm},
 \; \gamma^{(2q)}_{\br,\bm} =  \int_0^1 \prod_{j=1}^{2q}
   \varphi_{ \br_j,0,[Nx_{s+1}],m_j}  d x_{s+1},  \\
\varphi_{ \br,0,N,m}  =\frac{e(m N/P_{\br})-1}{P_{\br}(e(m/P_{\br})-1)}, \quad
    \beta^{(2q)}_{\br,\bm} = E_{s} \Big( e(\wp_{\br,\bm}) \prod_{j=1}^{2q}
 \psi_{ \br_j}(m_j,\bx)  \Big),
\end{multline}
where
\begin{multline} \label{Th2-30a}
 \varpi_{\br,\bm,2} = \b1 \Big(\sh=2q, \; \fs=s,\; \exists \eta \in \Xi_{2q} \; : \; m_{\eta(2k-1)}/P_{\br_{\eta(2k-1)}}=-m_{\eta(2k)}/P_{\br_{\eta(2k)}},\; \\
 k \in [1,q]    \Big),\;\;
 \wp_{\br,\bm}=\sum_{j=1}^{2q}  \frac{-m_j}{P_{\br_j}} \cX_{\br_j}, \quad
           | \varphi_{ \br,0,N,m} |  \leq  \frac{1}{\bar{m}}, \quad N=1,2,...,\\
        \psi_{ \br}(m,\bx)  = \prod_{i=1}^s \ddot{\psi}(i,\{-m M_{i,\br}/p_i \}p_i, x_{i, r_i}) ,\quad \quad
\ddot{\psi}(i,0,x_{i, r_i})=x_{i, r_i}, \\
				\ddot{\psi}(i,m', x_{i, r_i}) = \frac{1- e(-m' x_{i, r_i}/p_i)}{e(m'/p_i)-1} \;\;\;
	 \for \; m' \neq 0,\;\;\;\; |\ddot{\psi}(i,m', x_{i, r_i})| \leq p_i.
\end{multline}
 By \eqref{Beg-1a}, we have  $ I_M^{*}=[-[(M-1)/2],[M/2]] \cap \ZZ \setminus \{ 0 \}.$  Let
\begin{equation*}
  m_j= \ddot{m}_j P_{\balpha_j}\quad  \with \quad (\ddot{m}_j ,p_0)=1, \quad P_{\balpha_j} = p_1^{\alpha_{1,j}} \cdots p_s^{\alpha_{s,j}}
\end{equation*}
and let
\begin{equation}   \label{Th2-30b}
 I^{**}_{\br, \balpha }:=\{ k \in
 [-[(P_{\br }-1)/2]P^{-1}_{\balpha },[P_{\br }/2]P^{-1}_{\balpha }] \cap \ZZ, \;\;\; (k,p_0)=1 \}.
\end{equation}\\

{\bf Step 1.}  We want summarize  by $ \balpha_1,..., \balpha_{2q} $ at the beginning.
For this purpose, we introduce the affirmation $A_{\bm,\ddot{\bm}, \balpha} $ as follows :
\begin{equation*}
    A_{\bm,\ddot{\bm}, \balpha}:= (m_j= \ddot{m}_j P_{\balpha_j},\quad (\ddot{m}_j ,p_0)=1, \quad {\rm for\; all} \; j \in \{1, \ldots,2q\}).
\end{equation*}
It is easy to see
\begin{equation*}
  \sum_{\alpha_{i,j} \in [0, p_0 s  n], i \in[1,s], \; \ddot{m}_j \in I^{**}_{\br_j, \balpha_j },j \in [1,2q]}
\b1( A_{\bm,\ddot{\bm}, \balpha} )=1.
\end{equation*}
Changing the order of the summation, we get from \eqref{Th2-30a1}  and \eqref{Th2-30b}:
\begin{equation} \nonumber
\varsigma_{2q}	=\sum_{\substack{r_{i,j} \in [V_1,n], i \in[1,s]  \\ j \in [1,2q]}}
\sum_{\substack{m_j \in I^{*}_{P_{\br_j}} \\ j \in [1,2q]}}
  \sum_{\substack{\alpha_{i,j} \in [0, p_0 s  n], i \in[1,s]  \\ j \in [1,2q]}}
\sum_{\substack{\ddot{m}_j \in  I^{**}_{\br_j, \balpha_j } \\ j \in [1,2q]}}
\b1( A_{\bm,\ddot{\bm}, \balpha} )  \varpi_{\br,\bm,2} \; \gamma^{(2q)}_{\br,\bm}  \beta^{(2q)}_{\br,\bm}
\end{equation}
\begin{equation} \nonumber
=    \sum_{\substack{\alpha_{i,j} \in [0, p_0 s  n], i \in [1,s]  \\ j \in [1,2q]}} \;
\sum_{\substack{r_{i,j} \in [V_1,n], i \in[1,s]  \\ j \in [1,2q]}} \;
\sum_{\substack{\ddot{m}_j \in I^{**}_{\br_j, \balpha_j } \\ j \in [1,2q]}}
\sum_{\substack{m_j \in I^{*}_{P_{\br_j}} \\ j \in [1,2q]}}
\b1( A_{\bm,\ddot{\bm}, \balpha} )  \varpi_{\br,\bm,2} \; \gamma^{(2q)}_{\br,\bm}  \beta^{(2q)}_{\br,\bm}
\end{equation}
\begin{multline}  \label{Th2-30ad}
 =\sum_{\substack{\alpha_{i,j} \in [0,p_0 s  n], i \in[1,s]  \\ j \in [1,2q]}} \;
\sum_{\substack{r_{i,j} \in [V_1,n], i \in[1,s]  \\ j \in [1,2q]}} \;
\sum_{\substack{\ddot{m}_j \in  I^{**}_{\br_j, \balpha_j } \\ j \in [1,2q]}}
 \ddot{\varpi}_{\br,\ddot{\bm}P_{\balpha}} \; \gamma^{(2q)}_{\br,\ddot{\bm}P_{\balpha}}  \beta^{(2q)}_{\br,\ddot{\bm}P_{\balpha}}\quad \with  \\
 \ddot{\varpi}_{\br,\ddot{\bm}P_{\balpha}} := \sum_{\substack{m_j \in I^{*}_{P_{\br_j}} \\ j \in [1,2q]}}
\b1( A_{\bm,\ddot{\bm}, \balpha} )  \varpi_{\br,\bm,2} .
\end{multline}
%
We get from \eqref{Th2-30a} and \eqref{Th2-30ad}    that
\begin{multline} \label{Th2-34a}
 \ddot{\varpi}_{\br,\ddot{\bm}P_{\balpha}}  =
 \b1 \Big( \exists \eta \in \Xi_{2q} \;\; | \;\; \ddot{m}_{\eta(2k-1)} = -\ddot{m}_{\eta(2k)}, \\
 \br_{\eta(2k-1)} - \balpha_{\eta(2k-1)} =   \br_{\eta(2k)} - \balpha_{\eta(2k)}  ,\;\;
 k \in [1,q]   \Big).
\end{multline}
Our goal is to make the summations over  variables $\br, \bm$ and $\balpha$ independent of one another.
To do this, in \eqref{Th2-34} we will move from $\varsigma_{2q}$ to $ \ddot{\varsigma}_{2q}$.
%

By  \eqref{Th2-30a1} and \eqref{Th2-30a}, we obtain for $p_0=p_1 \cdots p_s$, that
\begin{multline}  \label{Th2-31a}
| \beta^{(2q)}_{\br,\ddot{\bm}P_{\balpha}}| \leq  p_0^{2q}, \quad | \varphi_{ \br_j,0,[Nx_{s+1}],m_j} | \leq
    \frac{1}{\bar{m}_j}= \frac{1}{\bar{\ddot{m}}_j P_{\balpha_j}}, \\
    | \varphi_{ \br_j,0,[Nx_{s+1}],m_j} |  \leq \max_{x_{s+1} \in [0,1]} \frac{|\sin(\pi m_j [Nx_{s+1}]/P_{\br_j})|}{\bar{m}_j} \leq \min( \frac{1}{\bar{m}_j},\pi N/P_{\br_j})\\
 \ad \qquad  |\gamma^{(2q)}_{\br,\ddot{\bm}P_{\balpha}} \; \beta^{(2q)}_{\br,\ddot{\bm}P_{\balpha}}|
  \leq \prod_{j=1}^{2q} p_0 \pi
  \min\big( \bar{\ddot{m}}_j^{-1}  P^{-1}_{\balpha_j} , N/P_{\br_j}\big).
\end{multline} \\

{\bf Step 2.} In order to single out the cases  $\b1(A_1)=\cdots =\b1(A_4)=0$, we introduce affirmations $A_1,...,A_4$  as follows:
\begin{multline} \label{Th2-31a5}
A_1:= (\exists \eta \in \Xi_{2q}, \quad  \exists k \in [1,q] \;\; : \;\;  |\ddot{m}_{\eta(2k-1)}|=|\ddot{m}_{\eta(2k)}| > n ),
      \qquad \qquad \qquad \qquad \qquad \qquad  \qquad \qquad \qquad \qquad \qquad \qquad \qquad \qquad \qquad \qquad \qquad \qquad  \qquad \qquad \qquad \\
\;\;\;A_2:= ( \exists (i,j) \;: \; \alpha_{i,j} >  \log_2 n ),  \qquad \qquad \qquad \qquad \qquad \qquad \qquad \qquad \qquad \qquad \qquad \qquad \qquad \qquad \qquad \\
A_3:= \big(\exists \eta \in \Xi_{2q},\;
 \exists (i,j_1 ,j_2) \; : \; j_1 \neq j_2, \;\;
 |r_{i,\eta(2j_1)}-\alpha_{i,\eta(2j_1)} - r_{i,\eta(2j_2)} + \alpha_{i,\eta(2j_2)}| <V_1 \\
\ad \quad  \br_{\eta(2j-1)} - \balpha_{\eta(2j-1)} =   \br_{\eta(2j)} - \balpha_{\eta(2j)} \quad {\rm for \; all} \; j \in \{1, \ldots,2q\}  \big),  \\
    A_4:=( \exists j \; : \;  P_{\br_j -\balpha_j} > 2^{n+\log_2^3 n} ),
    \qquad \qquad \qquad \qquad \qquad \qquad  \qquad \qquad  \qquad\qquad
\end{multline}
with $  V_1 = [\log_2^3 n],\; n = [\log_2 N]+1$. Let
\begin{equation}  \label{Th2-31a2}
\iota_{\nu}:= \sum_{\substack{\alpha_{i,j} \in [0, p_0 s n]\\ i \in[1,s],   j \in [1,2q]}}\;
\sum_{\substack{r_{i,j} \in [V_1,n]\\  i \in[1,s],  j \in [1,2q]}} \;
\sum_{\substack{\ddot{m}_j \in I^{**}_{P_{\br_j -\balpha_j}} \\ j \in [1,2q]}}
 \ddot{\varpi}_{\br,\ddot{\bm}P_{\balpha}} \; \gamma^{(2q)}_{\br,\ddot{\bm}P_{\balpha}}  \beta^{(2q)}_{\br,\ddot{\bm}P_{\balpha}}    \b1( A_{\nu}).
 \end{equation}
Let's consider the case $\b1(A_1)$=1:

Taking into account that
\begin{equation*}
   \sum_{|m| > n} |m|^{-2} \leq 4/n,
\end{equation*}
we have from \eqref{Th2-34a} and \eqref{Th2-31a}  that
\begin{multline}  \nonumber
\iota_{1} \ll  \sum_{\substack{\alpha_{i,j} \in [0, p_0 s n]\\ i \in[1,s],   j \in [1,2q]}}\;
\sum_{\substack{r_{i,j} \in [V_1,n]\\  i \in[1,s],  j \in [1,2q]}} \;
\sum_{\substack{\ddot{m}_j \in I^{**}_{\br_j, \balpha_j } \\ j \in [1,q]}} \;
\;(P_{\balpha_1} \cdots P_{\balpha_{2q}})^{-1}  \; \frac{1}{\bar{\ddot{m}}^2_1 \cdots \bar{\ddot{m}}^2_q}   \\
 \times \b1( \br_{j}  =   \br_{j+q} - \balpha_{j+q} +\balpha_j \quad {\rm for \; all} \; j \in \{1, \ldots,2q\}) \b1( A_1) \\
\end{multline}
\begin{multline}   \label{Th2-31a1}
\ll n^{-1}  \sum_{\substack{\alpha_{i,j} \in [0, p_0 s n]\\ i \in[1,s],   j \in [1,2q]}}\;
\sum_{\substack{r_{i,j} \in [V_1,n]\\  i \in[1,s],  j \in [1,2q]}} \;
\;  \frac{\b1( \br_{j}  =   \br_{j+q} - \balpha_{j+q} +\balpha_j \quad {\rm for\; all} \; j \in \{1, \ldots,2q\})}{P_{\balpha_1} \cdots P_{\balpha_{2q}}}\\
\ll n^{-1}  \sum_{\substack{\alpha_{i,j} \in [0, p_0 s n]\\ i \in[1,s],   j \in [1,2q]}}\;
(P_{\balpha_1} \cdots P_{\balpha_{2q}})^{-1} n^{qs}  \ll n^{qs-1}.
\end{multline}
Let's consider the case $\b1(A_2)$=1:

Bearing in mind \eqref{Th2-34a} and \eqref{Th2-31a}, we get,   similarly to \eqref{Th2-31a1}, that
\begin{multline}  \nonumber 
\iota_{2} \ll   \sum_{\substack{i_0 \in [1, s]\\ j_0 \in[1,2q]}}\;
 \sum_{\substack{\alpha_{i,j} \in [0, p_0 s n]\\ i \in[1,s],   j \in [1,2q] \\ \alpha_{i_0,j_0} \geq \log_2 n }}\;
\sum_{\substack{r_{i,j} \in [V_1,n]\\  i \in[1,s],  j \in [1,2q]}} \;
\sum_{\substack{\ddot{m}_j \in I^{**}_{\br_j, \balpha_j } \\ j \in [1,q]}} \;
\;  \; \frac{(P_{\balpha_1} \cdots P_{\balpha_{2q}})^{-1}}{\bar{\ddot{m}}^2_1 \cdots \bar{\ddot{m}}^2_q}   \\
 \times \b1( \br_{j}  =   \br_{j+q} - \balpha_{j+q} +\balpha_j \quad {\rm for \; all} \; j \in \{1, \ldots,2q\})
\end{multline}
\begin{multline}  \nonumber
\ll
 \sum_{\substack{i_0 \in [1, s]\\ j_0 \in[1,2q]}}\;
 \sum_{\substack{\alpha_{i,j} \in [0, p_0 s n]\\ i \in[1,s],   j \in [1,2q] \\ \alpha_{i_0,j_0} \geq \log_2 n }}\;
\sum_{\substack{r_{i,j} \in [V_1,n]\\  i \in[1,s],  j \in [1,2q]}}
\; \frac{1}{P_{\balpha_1} \cdots P_{\balpha_{2q}}}  \\
 \times \b1( \br_{j}  =   \br_{j+q} - \balpha_{j+q} +\balpha_j \quad {\rm for \; all} \; j \in \{1, \ldots,2q\}) \\
\ll  n^{qs}
 \sum_{\substack{i_0 \in [1, s]\\ j_0 \in[1,2q]}}\;
 \sum_{\substack{\alpha_{i,j} \in [0, p_0 s n]\\ i \in[1,s],   j \in [1,2q] \\ \alpha_{i_0,j_0} \geq \log_2 n }}\;
 \frac{1}{P_{\balpha_1} \cdots P_{\balpha_{2q}}}   \ll n^{qs-1}.
\end{multline}
Let's consider the case $\b1(A_3)$=1:

Bearing in mind \eqref{Th2-34a} and \eqref{Th2-31a}, we get,  similarly to \eqref{Th2-31a1}, that

\begin{multline} \nonumber 
\iota_{3} \ll \sum_{\substack{i_0 \in [1, s]\\ j_0 \in[1,2q]}}\;
 \sum_{\substack{\alpha_{i,j} \in [0, p_0 s n]\\ i \in[1,s],   j \in [1,2q]}}\;
  \sum_{\substack{r_{i,j} \in [V_1,n]\\  i \in[1,s],  j \in [1,2q]}} \;
\sum_{\substack{\ddot{m}_j \in I^{**}_{\br_j, \balpha_j } \\ j \in [1,q]}} \;
\;  \; \frac{(P_{\balpha_1} \cdots P_{\balpha_{2q}})^{-1}}{\bar{\ddot{m}}^2_1 \cdots \bar{\ddot{m}}^2_q}   \\
 \times \sum_{ \eta \in \Xi_q}
   \sum_{\substack{i_0 \in [1, s]\\ j_1, j_2 \in[1,q], j_1 \neq j_2}}\;
 \b1\Bigg( |r_{i,\eta(2j_1)}-\alpha_{i,\eta(2j_1)} - r_{i,\eta(2j_2)} + \alpha_{i,\eta(2j_2)}| <V_1 \\
\ad \quad  \br_{\eta(2j-1)} - \balpha_{\eta(2j-1)} =   \br_{\eta(2j)} - \balpha_{\eta(2j)} \quad {\rm for \; all} \; j \in \{1, \ldots,2q\}  \Bigg)
\end{multline}
\begin{multline} \nonumber
\ll n^{-1} V_1
 \sum_{\substack{\alpha_{i,j} \in [0, p_0 s n]\\ i \in[1,s],   j \in [1,2q] }}\;
\sum_{\substack{r_{i,j} \in [V_1,n]\\  i \in[1,s],  j \in [1,2q]}} \;
\sum_{\substack{\ddot{m}_j \in I^{**}_{\br_j, \balpha_j } \\ j \in [1,q]}} \;
\;  \; \frac{(P_{\balpha_1} \cdots P_{\balpha_{2q}})^{-1}}{\bar{\ddot{m}}^2_1 \cdots \bar{\ddot{m}}^2_q}   \\
 \times \b1( \br_{j}  =   \br_{j+q} - \balpha_{j+q} +\balpha_j \quad {\rm for \; all} \; j \in \{1, \ldots,2q\}) .
\end{multline}
Now from  \eqref{Th2-31a1}, we get  $\iota_3 \ll n^{qs-1} (\log_2 n)^3$.\\ \\
Let's consider the case $\b1(A_4)$=1:
By \eqref{Th2-31a}, we obtain
\begin{equation*}
| \gamma^{(2q)}_{\br,\ddot{\bm}P_{\balpha}}  \; \beta^{(2q)}_{\br,\ddot{\bm}P_{\balpha}}|   \;  \b1( A_4)
\ll
\prod_{j=1}^{2q} \min\big( P^{-1}_{\balpha_j} /|\bar{\ddot{m}}_j|, N/P_{\br_j}\big)
\ll  N 2^{-n- \log_2^3 n} \leq 2^{- \log^3_2 n}.
\end{equation*}
In view of   \eqref{Th2-31a2}, we get
\begin{equation}  \nonumber
\iota_4 \ll  \sum_{\substack{\alpha_{i,j} \in [0, p_0 2s n]\\ i \in[1,s],   j \in [1,2q]}}\;
\sum_{\substack{r_{i,j} \in [V_1,n]\\  i \in[1,s],  j \in [1,2q]}} \;
\sum_{\substack{\ddot{m}_j \in I^{**}_{P_{\br_j -\balpha_j}} \\ j \in [1,2q]}}
   2^{- \log^3_2 n} \ll n^{10 s q}  2^{-\log^3_2 n}
 \ll 1/n.
 \end{equation}

%
Combining the above cases, we get from  \eqref{Th2-30ad}, that
\begin{multline*}
\varsigma_{2q} =\sum_{\substack{\alpha_{i,j} \in [0,p_0 s  n], i \in[1,s]  \\ j \in [1,2q]}} \;
\sum_{\substack{r_{i,j} \in [V_1,n], i \in[1,s]  \\ j \in [1,2q]}} \;
\sum_{\substack{\ddot{m}_j \in  I^{**}_{\br_j, \balpha_j } \\ j \in [1,2q]}}
 \ddot{\varpi}_{\br,\ddot{\bm}P_{\balpha}} \; \gamma^{(2q)}_{\br,\ddot{\bm}P_{\balpha}}  \beta^{(2q)}_{\br,\ddot{\bm}P_{\balpha}} \\
 = \sum_{\substack{\alpha_{i,j} \in [0,p_0 s  n], i \in[1,s]  \\ j \in [1,2q]}} \;
\sum_{\substack{r_{i,j} \in [V_1,n], i \in[1,s]  \\ j \in [1,2q]}} \;
\sum_{\substack{\ddot{m}_j \in  I^{**}_{\br_j, \balpha_j } \\ j \in [1,2q]}}
 \ddot{\varpi}_{\br,\ddot{\bm}P_{\balpha}} \; \gamma^{(2q)}_{\br,\ddot{\bm}P_{\balpha}}  \beta^{(2q)}_{\br,\ddot{\bm}P_{\balpha}} \prod_{j=1}^4 (1- \b1( A_j)) \\
 +  O(n^{qs-9/10}).
\end{multline*}
Taking into account that $\b1(A_1)=\b1(A_2)=0$ (i.e.,
$0 \leq \alpha_{i,j} \leq \log_2 n $  and $ \ddot{m}_j \in  I^{\prime}_n $
for all $1 \leq i \leq s$, $1 \leq j \leq 2q $ with $I^{\prime}_{n} :=\{ k \in [-n,n] \;  | \;k \neq 0, (k,p_0)=1\}$),
we have
\begin{equation}  \label{Th2-34}
\varsigma_{2q} = \ddot{\varsigma}_{2q} + O(n^{qs-9/10}),
\end{equation}
where
\begin{equation} \nonumber
     \ddot{\varsigma}_{2q}
	=  \sum_{\substack{\alpha_{i,j} \in [0, \log_2 n]\\ i \in[1,s], j \in [1,2q]}} \;\;
\sum_{\substack{r_{i,j} \in [V_1,n]\\ i \in[1,s], j \in [1,2q]}} \;\;
\sum_{\substack{\ddot{m}_j \in I^{\prime}_{n}  \\ j \in [1,2q]}} \ddot{\varpi}_{\br,\ddot{\bm}P_{\balpha}} \; \gamma^{(2q)}_{\br,\ddot{\bm}P_{\balpha}}   \beta^{(2q)}_{\br,\ddot{\bm}P_{\balpha}}   \prod_{j=3}^4 (1- \b1( A_j)),
\end{equation}
Thus,  we selected the case $\b1(A_1)=\cdots =\b1(A_4)=0$.
\\ \\

{\bf Step 3.}  We want to collect  vectors $\br_1,...,\br_{2q} $ with the same permutations $\eta$ in \eqref{Th2-34a}.
For this purpose, we introduce the set $\tilde{\Xi}_q $ and the sequences
$ \ell_{\varrho_1,\varrho_2,\br,\balpha} $
as follows :
\begin{multline*}
  \tilde{\Xi}_q = \Big\{ (\varrho_1,\varrho_2),  \varrho_i: \{1,...,q\} \to \{1,...,2q \}, \;i=1,2\;\; :\;\;
  \big\{\varrho_1(1),\varrho_2(1),...,\\ \varrho_1(q),  \varrho_2(q) \big\}
   = \big\{1,...,2q \big\},\;  \;\varrho_1(k)<\varrho_2(k) \; \forall k  \Big\}
\end{multline*}
and let
\begin{equation}  \label{Th2-33ab}
  \ell_{\varrho_1,\varrho_2,\br,\balpha} =\b1(  \br_{\varrho_1(k)} - \balpha_{\varrho_1(k)}=\br_{\varrho_2(k)} -\balpha_{\varrho_2(k)} \;\; {\rm for \; all} \;\; k=1,...,q )   .
\end{equation}
It is easy to verify that
\begin{equation}   \label{Th2-33a}
  \card(\tilde{\Xi}_q )   =  \binom{2q}{2}  \binom{2q-2}{2}  \cdots  \binom{2}{2}
    = \frac{(2q)!}{ 2^{q}}.
\end{equation}
Suppose that $ \ddot{\varpi}_{\br,\ddot{\bm}P_{\balpha}} =1$ and $ \b1( A_3)=0$ (see \eqref{Th2-34a}, \eqref{Th2-31a5}).
Then there exists $(\varrho_1,\varrho_2) \in \tilde{\Xi}_q$ with $\ell_{\varrho_1,\varrho_2,\br,\balpha}=1$ and for a fixed $\varrho_1$ such $\varrho_2$ is unique.

From \eqref{Th2-34a} and \eqref{Th2-31a5}, we have
\begin{equation} \nonumber
\sum_{(\varrho_1,\varrho_2) \in \tilde{\Xi}_q}  \ell_{\varrho_1,\varrho_2,\br,\balpha} \;
 \ddot{\varpi}_{\br,\ddot{\bm}P_{\balpha}} \;    (1- \b1( A_3)) =
q! \;\ddot{\varpi}_{\br,\ddot{\bm}P_{\balpha}} \;    (1- \b1( A_3)).
\end{equation}
Let
\begin{multline}  \label{Th2-35}
\partial^{(\bfr,\bmu,\varrho)}_{\br,\ddot{\bm}, \balpha} =\b1 \big(
  \bfr_{k}=\br_{\varrho_1(k)} - \balpha_{\varrho_1(k)}=\br_{\varrho_2(k)} -\balpha_{\varrho_2(k)},\\
 \mu_{k}= \ddot{m}_{\varrho_1(k)}=-\ddot{m}_{\varrho_2(k)} \;\;\;\; {\rm for \; all} \;\; k \in \{1,...,q \}   \big).
\end{multline}
We see that
$$
\fr_{k,i} \geq V_1 -\log_2 n >V_1/2, \quad i=1,\ldots,s, \; k=1, \ldots, q.
$$
It is easy to verify
\begin{equation*}
 \sum_{\substack{\fr_{i,j} \in [V_1/2,n]\\ i \in[1,s], j \in [1,2q]}} \;\;
\sum_{\substack{ \mu_j \in  I^{\prime}_{n}      \\ j \in [1,q]}}
\partial^{(\bfr,\bmu,\varrho)}_{\br,\ddot{\bm}, \balpha} \;  \ell_{\varrho_1,\varrho_2,\br,\balpha} \; \ddot{\varpi}_{\br,\ddot{\bm}P_{\balpha}} (1- \b1( A_4)) =   \ell_{\varrho_1,\varrho_2,\br,\balpha}\;
 \ddot{\varpi}_{\br,\ddot{\bm}P_{\balpha}}\; (1- \b1( A_4)).
\end{equation*}
In view of \eqref{Th2-34}, we have that
\begin{multline}    \nonumber 
 \ddot{\varsigma}_{2q}
	=  \frac{1}{q!} \sum_{(\varrho_1,\varrho_2) \in \tilde{\Xi}_q}
 \sum_{\substack{\alpha_{i,j} \in [0, \log_2 n]\\ i \in[1,s], j \in [1,2q]}} \;\;
\sum_{\substack{r_{i,j} \in [V_1,n]\\ i \in[1,s], j \in [1,2q]}} \;\;
\sum_{\substack{\ddot{m}_j \in I^{\prime}_{n}  \\ j \in [1,2q]}}
\ell_{\varrho_1,\varrho_2,\br,\balpha}  \;
\ddot{\varpi}_{\br,\ddot{\bm}P_{\balpha}} \; \gamma^{(2q)}_{\br,\ddot{\bm}P_{\balpha}}\\
   \times   \beta^{(2q)}_{\br,\ddot{\bm}P_{\balpha}}   \prod_{j=3}^4 (1- \b1( A_j))
\end{multline}
and
\begin{multline*}
 \ddot{\varsigma}_{2q} 	=   \frac{1}{q!} \sum_{(\varrho_1,\varrho_2) \in \tilde{\Xi}_q}
 \sum_{\substack{\alpha_{i,j} \in [0, \log_2 n]\\ i \in[1,s], j \in [1,2q]}} \;\;
\sum_{\substack{r_{i,j} \in [V_1,n]\\ i \in[1,s], j \in [1,2q]}} \;\;
\sum_{\substack{\ddot{m}_j \in I^{\prime}_{n}   \\ j \in [1,2q]}}
  \sum_{\substack{\fr_{i,j} \in [V_1/2,n]\\ i \in[1,s], j \in [1,q]}} \;\;
\sum_{\substack{\mu_j \in I^{\prime}_{n}  \\ j \in [1,q]}} 1
\\
  \times   \partial^{(\bfr,\bmu,\varrho)}_{\br,\ddot{\bm}, \balpha} \;
  \ell_{\varrho_1,\varrho_2,\br,\balpha} \;
  \ddot{\varpi}_{\br,\ddot{\bm}P_{\balpha}} \; \gamma^{(2q)}_{\br,\ddot{\bm}P_{\balpha}}  \beta^{(2q)}_{\br,\ddot{\bm}P_{\balpha}}   \prod_{j=3}^4 (1- \b1( A_j)).
\end{multline*}
Changing the order of the summation, we get
\begin{equation}  \label{Th2-33}
   \ddot{\varsigma}_{2q} 	=   \frac{1}{q!} \sum_{(\varrho_1,\varrho_2) \in \tilde{\Xi}_q}
    \sum_{\substack{\alpha_{i,j} \in [0,  \log_2 n]\\ i \in[1,s], j \in [1,2q]}} \;\; \sum_{\substack{\fr_{i,j} \in [ V_1/2,n]\\ i \in[1,s], j \in [1,q]}} \;\;
\sum_{\substack{  \mu_j  \in I^{\prime}_{n}  \\ j \in [1,q]}} Z_0,
\end{equation}
where
\begin{equation} \nonumber
 Z_0=\sum_{\substack{r_{i,j} \in [V_1,n]\\ i \in[1,s], j \in [1,2q]}} \;\;
\sum_{\substack{\ddot{m}_j \in I^{\prime}_{n}   \\ j \in [1,2q]}}
\partial^{(\bfr,\bmu,\varrho)}_{\br,\ddot{\bm}, \balpha} \;
 \ell_{\varrho_1,\varrho_2,\br,\balpha} \;
  \ddot{\varpi}_{\br,\ddot{\bm}P_{\balpha}} \; \gamma^{(2q)}_{\br,\ddot{\bm}P_{\balpha}}  \beta^{(2q)}_{\br,\ddot{\bm}P_{\balpha}}   \prod_{j=3}^4 (1- \b1( A_j)).
\end{equation}
Let
\begin{multline} \label{Th2-46a}
 \tilde{\kappa}_{\fr,\balpha}:
 = \b1 \big(
  \min_{i,j} \fr_{i,j}+\alpha_{i,j} \geq V_1  , \;\; \max_{i,j}   \fr_{i,j}+\alpha_{i,j} \leq n, \;\;   \max_j P_{\fr_j} \leq 2^{n+\log_2^3 n}, \\
\min_{i,j_1,j_2, j_1\neq j_2} |\fr_{i,j_1}- \fr_{i,j_2}| \geq V_1 \big).
\end{multline}
In view of \eqref{Th2-34a}, \eqref{Th2-31a5}, \eqref{Th2-33ab} and \eqref{Th2-35},  we have
\begin{equation}  \label{Th2-46b}
 Z_0 =  \sum_{\substack{\ddot{m}_j \in I^{\prime}_{n}   \\ j \in [1,2q]}}
 \b1( \mu_{k}= \ddot{m}_{\varrho_1(k)}=-\ddot{m}_{\varrho_2(k)}, \;\;  k \in [1,q] )
 \tilde{\kappa}_{\bfr,\balpha}\; \gamma^{(2q)}_{\bfr+\balpha,\ddot{\bm}P_{\balpha}}  \beta^{(2q)}_{\bfr+\balpha,\ddot{\bm}P_{\balpha}}.
\end{equation}
Thus, we have done in \eqref{Th2-33} summation over permutations at the beginning.

Taking
into account  \eqref{Th2-33a}, this will allow us at Step 6 to obtain the coefficient $\frac{(2q)!}{ 2^{q} q!} $ in \eqref{Th2-30}. \\ \\

{\bf Step 4.} In this step, we will rewrite the parameters $\gamma^{(2q)}_{\br,\ddot{\bm}P_{\balpha}}  $ and
 $\beta^{(2q)}_{\br,\ddot{\bm}P_{\balpha}}$  using  variables  $\ddot{\gamma}^{(q)}_{\br,\bm}$ and $\tilde{\beta}^{(q,\varrho)}_{\bfr,\bmu,\balpha}$  (see \eqref{Th2-50} and \eqref{Th2-35b}):

According to \eqref{Th2-50}, we obtain
\begin{equation*}
  \ddot{\gamma}^{(q)}_{\br,\bm}  =  \frac{1}{N} \sum_{k=0}^{N-1} \prod_{j=1}^{q}
 \frac{2-2 \cos(2 \pi m_jk/P_{\br_j} )}{|P_{\br_j}(1-e(m_j/P_{\br_j}))|^2} =
  \frac{1}{N} \sum_{k=0}^{N-1} \prod_{j=1}^{q}
 \frac{|e( m_jk/P_{\br_j}) -1|^2 }{|P_{\br_j}(e(m_j/P_{\br_j})-1)|^2}.
\end{equation*}
By \eqref{Th2-30a1},  \eqref{Th2-35} and \eqref{Th2-00} , we derive
\begin{multline} \label{Th2-46c}
\gamma^{(2q)}_{\bfr +\balpha,\ddot{\bm}P_{\balpha}} = \gamma^{(2q)}_{\br,\ddot{\bm}P_{\balpha}} = \frac{1}{N} \sum_{k=0}^{N-1} \prod_{j=1}^{2q}
 \frac{e(\ddot{m}_jk/P_{\br_j-\balpha_j})-1}{P_{\br_j}(e(\ddot{m}_j/P_{\br_j-\balpha_j}) -1)} \\
  =     \frac{1}{N P_{\balpha_1}\cdots P_{\balpha_{2q}}} \sum_{k=0}^{N-1} \prod_{j=1}^{q}
 \frac{|e(\mu_j k /P_{\bfr_j})-1|^2}{|P_{\bfr_j}(e(\mu_j/P_{\bfr_j})-1)|^2}=  \ddot{\gamma}^{(q)}_{\bfr,\bmu} /(P_{\balpha_1}\cdots P_{\balpha_{2q}}).
\end{multline}
Now we will calculate $ \beta^{(2q)}_{\br,\ddot{\bm}P_{\balpha}} $.
From \eqref{Th2-30a1}, \eqref{Th2-30a}, \eqref{Beg-3}  and \eqref{Beg-6}, we have
\begin{multline}   \label{Th2-35a}
 \beta^{(2q)}_{\br,\ddot{\bm}P_{\balpha}} =
 E_{s} \Big( e(\wp_{\br,\ddot{\bm}P_{\balpha} }) \prod_{j=1}^{2q}
 \psi_{ \br_j}( \ddot{\bm}_jP_{\balpha_j},\bx)  \Big), \quad
     M_{i,\br} \equiv
	\big( P_{\br}/p_i^{r_i} \big)^{-1} \mod p_i^{r_i},\\
 \wp_{\br,\bm}=\sum_{j=1}^{2q}  \frac{-m_j}{P_{\br_j}} \cX_{\br_j}, \quad
  \cX_{\br} \equiv \sum_{i=1}^s  M_{i,\br}
   P_{\br}  p_i^{-r_i} \cX_{i,r_i} \;(\mod  P_{\br}), \quad \cX_{i,r_i}=
  \sum_{1 \leq j \leq r}  x_{i,j} p_i^{j-1}.
\end{multline}
Using  \eqref{Th2-35}, we define
\begin{multline*}  
\tilde{\beta}^{(q,\varrho)}_{\bfr,\bmu,\balpha} :=   \beta^{(2q)}_{\br,\ddot{\bm}P_{\balpha}} =
  E_{s} \Big( \tilde{\wp}  \prod_{j=1}^q
 \psi_{ \bfr_j +\balpha_{\varrho_1(j)}}(\mu_jP_{\balpha_{\varrho_1(j)}},\bx)
 \psi_{ \bfr_j +\balpha_{\varrho_2(j)}}(-\mu_jP_{\balpha_{\varrho_2(j)}},\bx) \Big),\\
  \with \quad \tilde{\wp}:=e\big(-\sum_{j=1}^{q}  \frac{\mu_j}{P_{\bfr_j}}
    ( \cX_{\bfr_j + \balpha_{\varrho_1(j)}}  -  \cX_{\bfr_j + \balpha_{\varrho_2(j)}}    \big) \quad \ad \quad
    \bfr_j=\br_j-\balpha_{\sigma_1(j)}, \; \mu_j = \ddot{m}_{\sigma_1(j)}.
\end{multline*}
First, we will prove that $\tilde{\wp}=1 $. To do this, it's enough to make sure that
\begin{equation}   \nonumber
     \cX_{\bfr_j + \tilde{\balpha}} / P_{\bfr_j} \equiv
     \cX_{\bfr_j} / P_{\bfr_j}   \mod 1
     \quad {\rm for \; all} \quad
     \tilde{\balpha} =(\tilde{\alpha}_1, ..., \tilde{\alpha}_s ), \;\; \tilde{\alpha}_i \geq 0, \; i \in [1,s], \; j \in[1,q] .
\end{equation}
By \eqref{Th2-35a}, we get
\begin{multline}
 \prod_{1 \leq k \leq s, k \neq i} p_k^{\tilde{\alpha}_k } \;\; M_{i, \bfr_j + \tilde{\balpha} }
   \equiv   \prod_{1 \leq k \leq s, k \neq i} p_k^{\tilde{\alpha}_k } \;\;
    (P_{\bfr_j + \tilde{\balpha} } p_i^{ -\fr_{i,j}-\tilde{\alpha}_i } )^{-1} \\
     \equiv   \prod_{1 \leq k \leq s, k \neq i} p_k^{-\fr_{i,j} }
    \equiv    M_{i, \bfr_j } \mod p_i^{ \fr_{j,i}}
\end{multline}
and
\begin{multline*}
  \cX_{\bfr_j + \tilde{\balpha}} / P_{\bfr_j} \equiv
  \sum_{i=1}^s M_{i,\bfr_j + \tilde{\balpha}}
   P_{\bfr_j + \tilde{\balpha}}  p_i^{- \fr_{j,i} - \tilde{\alpha}_i}  \cX_{i,\fr_{j,i}  + \tilde{\alpha}_i}/ P_{\bfr_j}    \\
\equiv   \sum_{i=1}^s M_{i,\bfr_j + \tilde{\balpha}}
   P_{\tilde{\balpha}}  p_i^{- \fr_{j,i}  - \tilde{\alpha}_i}  \cX_{i, \fr_{j,i}  + \tilde{\alpha}_i}
  \equiv   \sum_{i=1}^s  \Big( M_{i,\bfr_j + \tilde{\balpha}}   \prod_{1 \leq k \leq s, k \neq i} p_k^{\tilde{\alpha}_k } \Big) \;
     p_i^{- \fr_{j,i}  }  \cX_{i,\fr_{j,i} }   \\
   \equiv \sum_{i=1}^s M_{i,\bfr_j}
     p_i^{-  \fr_{j,i}  }  \cX_{i,\fr_{j,i} }
      \equiv    \cX_{\bfr_j } / P_{\bfr_j} \mod 1.
\end{multline*}
Hence $\tilde{\wp}=1 $ and
\begin{equation}  \label{Th2-35b}
\tilde{\beta}^{(q,\varrho)}_{\bfr,\bmu,\balpha} =   \beta^{(2q)}_{\bfr+\balpha,\ddot{\bm}P_{\balpha}} =
  E_{s} \Big( \prod_{j=1}^q
 \psi_{ \bfr_j +\balpha_{\varrho_1(j)}}(\mu_jP_{\balpha_{\varrho_1(j)}},\bx)
 \psi_{ \bfr_j +\balpha_{\varrho_2(j)}}(-\mu_jP_{\balpha_{\varrho_2(j)}},\bx) \Big).
\end{equation}
In view of \eqref{Th2-30a}, we have that $ \psi_{ \bfr_j +\balpha_{\varrho_1(j)}}(\mu_jP_{\balpha_{\varrho_1(j)}},\bx) $ depend only on \\ $(x_{1,\fr_{1,j} + \alpha_{1,j} }, \ldots, x_{s,\fr_{s,j} + \alpha_{s,j} }) $.

From \eqref{Th2-46a} and  \eqref{Beg-22}, we get  $|\fr_{i,j_1}-  \fr_{i, j_2}| \geq V_1  $ for $j_1 \neq j_2, \; i=1,...,s$  with $ V_1=[\log^3_2 n] $.\\
Hence $|\fr_{i,j_1}+  \alpha_{i,j_3}-  \fr_{i, j_2} -
 \alpha_{i,j_4} | \geq V_1/2  $ for $j_1 \neq j_2, \; i=1,...,s$.

Therefore
 expectation and multiplication can be interchanged:
\begin{multline}  \label{Th2-45}
  \tilde{\beta}^{(q,\varrho)}_{\bfr,\bmu,\balpha} =  \prod_{j=1}^{q} \hat{\beta}^{(q,\varrho,j)}_{\bfr,\bmu,\balpha},\\
\with \quad   \hat{\beta}^{(q,\varrho,j)}_{\bfr,\bmu,\balpha}:=  E_{s} \big(  \psi_{ \bfr_j +\balpha_{\varrho_1(j)}}(\mu_jP_{\balpha_{\varrho_1(j)}},\bx)
 \psi_{ \bfr_j +\balpha_{\varrho_2(j)}}(-\mu_jP_{\balpha_{\varrho_2(j)}},\bx) \big).
\end{multline}
Using \eqref{Th2-46a}, \eqref{Th2-46b}, \eqref{Th2-46c} and \eqref{Th2-35b}),  we obtain
\begin{equation} \label{Th2-46}
 Z_0 =  \tilde{\kappa}_{\bfr,\balpha}\;
     \ddot{\gamma}^{(q)}_{\bfr,\bmu}  \;
       \tilde{\beta}^{(q,\varrho)}_{\bfr,\bmu,\balpha}/ (P_{\balpha_1} \cdots  P_{\balpha_{2q}}).
\end{equation}
So, we replaced $\gamma^{(2q)}_{\br,\ddot{\bm}P_{\balpha}}  $ and
 $\beta^{(2q)}_{\br,\ddot{\bm}P_{\balpha}}$ by $ \ddot{\gamma}^{(q)}_{\bfr,\bmu}$ and $       \tilde{\beta}^{(q,\varrho)}_{\bfr,\bmu,\balpha}   $. \\ \\

{\bf Step 5.} Now we apply Lemma 7.
Using \eqref{Th2-33}, \eqref{Th2-46}    and Lemma 7,  we obtain
\begin{equation} \label{Th2-34a1}
   \ddot{\varsigma}_{2q} 	= \frac{1}{q!} \sum_{(\varrho_1,\varrho_2) \in \tilde{\Xi}_q}
  \sum_{\substack{\alpha_{i,j} \in [0, \log_2 n]\\ i \in[1,s], j \in [1,2q]}} \;\; \sum_{\substack{\fr_{i,j} \in [ V_1/2,n]\\ i \in[1,s], j \in [1,q]}} \;\;
\sum_{\substack{ \mu_j \in I^{\prime}_{n}  \\ j \in [1,q]}}
\frac{\tilde{\kappa}_{\bfr,\balpha}     \ddot{\gamma}^{(q)}_{\bfr,\bmu}
 \tilde{\beta}^{(q,\varrho)}_{\bfr,\bmu,\balpha}}{P_{\balpha_1} \cdots  P_{\balpha_{2q}}}=\breve{\varsigma}_{2q} +O(n^{qs-\frac{9}{10}}),
\end{equation}
where
\begin{equation} \label{Th2-34a3}
\breve{\varsigma}_{2q}= \frac{1}{q!} \sum_{(\varrho_1,\varrho_2) \in \tilde{\Xi}_q}
    \sum_{\substack{\alpha_{i,j} \in [0, \log_2 n]\\ i \in[1,s], j \in [1,2q]}}
 \sum_{\substack{\fr_{i,j} \in [V_1/2,n]\\ i \in[1,s],j \in [1,q]}} \;\;
\sum_{\substack{\mu_j \in I^{\prime}_{n}  \\ j \in [1,q]}} \prod_{j=1}^{q}  \frac{2 P^{-1}_{\balpha_{\varrho_1(j)}}P^{-1}_{\balpha_{\varrho_2(j)}} \tilde{\kappa}_{\bfr,\balpha} }{|P_{\bfr_j}(1-e(\mu_j/P_{\bfr_j}))|^2} \hat{\beta}^{(q,\varrho,j)}_{\bfr,\bmu,\balpha}.
\end{equation}
Taking into account that
$$
\tilde{\kappa}_{\bfr,\balpha}  \in \{ 0,1 \} , \; \hat{\beta}^{(q,\varrho,j)}_{\bfr,\bmu,\balpha} \ll 1 \; \ad \;  \mu_j^2 \ll |P_{\bfr_j}(1-e(\mu_j/P_{\bfr_j}))|^2  ,
$$
we get that
the part of \eqref{Th2-34a3}, satisfying the condition \\ $\min_{i,j} (\fr_{i,j}+\alpha_{i,j}) < V_1 $ is equal to $O(n^{qs-\frac{9}{10}})$.
The same estimate is true for the cases
\begin{equation*}
    \max_{i,j} \fr_{i,j}+\alpha_{i,j} > n , \quad   \quad \min_{i,j_1,j_2, j_1\neq j_2} |\fr_{i,j_1}- \fr_{i,j_2}| < V_1.
\end{equation*}
%
Now applying \eqref{Th2-46a}, we get
\begin{multline} \nonumber
 \sum_{(\varrho_1,\varrho_2) \in \tilde{\Xi}_q}
    \sum_{\substack{\alpha_{i,j} \in [0, \log_2 n]\\ i \in[1,s], j \in [1,2q]}} \;
 \sum_{\substack{\fr_{i,j} \in [V_1/2,n],  P_{\bfr_j} \leq 2^{n+\log_2^3 n} \\  i \in[1,s],j \in [1,q]}} \;\sum_{\substack{\mu_j \in I^{\prime}_{n} \\ j \in [1,q]}} \; \prod_{j=1}^{q}  \frac{2 P^{-1}_{\balpha_{\varrho_1(j)}}P^{-1}_{\balpha_{\varrho_2(j)}} |\tilde{\kappa}_{\bfr,\balpha}-1| }{|P_{\bfr_j}(1-e(\mu_j/P_{\bfr_j}))|^2}    \\
  \times   \hat{\beta}^{(q,\varrho,j)}_{\bfr,\bmu,\balpha} \ll n^{qs-9/10} .
\end{multline}
Therefore
\begin{equation}   \label{Th2-34b}
   \breve{\varsigma}_{2q}=\hat{\varsigma}_{2q} \; + \; O( n^{qs-9/10}).
\end{equation}
with
\begin{multline} \nonumber
\hat{\varsigma}_{2q}= \frac{1}{q!} \sum_{(\varrho_1,\varrho_2) \in \tilde{\Xi}_q} \;
    \sum_{\substack{\alpha_{i,j} \in [0, \log_2 n],i \in[1,s]  \\    P_{\bfr_j} \leq 2^{n+\log_2^3 n}, j \in [1,2q]}} \;
\sum_{\substack{\fr_{i,j} \in [V_1/2,n],  \\  i \in[1,s],j \in [1,q]}}
   \;\;\sum_{\substack{\mu_j \in I^{\prime}_{n} \\ j \in [1,q]}} \;\;
\prod_{j=1}^{q}  \frac{2 P^{-1}_{\balpha_{\varrho_1(j)}}P^{-1}_{\balpha_{\varrho_2(j)}} }{|P_{\bfr_j}(1-e(\mu_j/P_{\bfr_j}))|^2} \\
  \times \hat{\beta}^{(q,\varrho,j)}_{\bfr,\bmu,\balpha}.
\end{multline}\\

{\bf Step 6.} In this last step, we show the connection between $\hat{\varsigma}_{2q}$ and $\varsigma_{2}^q $ and complete the proof of the lemma.

Changing the order of the summation, we obtain
\begin{multline}  \label{Th2-34abc}
\hat{\varsigma}_{2q}=  \frac{1}{q!} \sum_{(\varrho_1,\varrho_2) \in \tilde{\Xi}_q} \;
 \sum_{\substack{\fr_{i,j} \in [V_1/2,n], P_{\bfr_j} \leq 2^{n+\log_2^3 n}\\ i \in[1,s], j \in [1,q]   }}
\sum_{\substack{\mu_j \in I^{\prime}_{n}  \\ j \in [1,q]}} \prod_{j=1}^{q} \frac{2}{|P_{\bfr_j}(1-e(\mu_j/P_{\bfr_j}))|^2}
G_{\varrho,\bfr_j,\bmu_j}, \\
\with \quad G_{\varrho,\bfr_j,\bmu_j}:=
 \sum_{\substack{\alpha_{i,\varrho_1(j)}, \alpha_{i,\varrho_2(j)} \in [0, \log_2 n]\\ i \in[1,s]}}
 P^{-1}_{\balpha_{\varrho_1(j)}}P^{-1}_{\balpha_{\varrho_2(j)}} \hat{\beta}^{(q,\varrho,j)}_{\bfr,\bmu,\balpha}.
\end{multline}
By \eqref{Th2-45}, we have
\begin{multline*}
G_{\varrho,\bfr_j,\bmu_j}=   \sum_{\substack{\alpha_{i,\varrho_1(j)}, \alpha_{i,\varrho_2(j)}\\ \in [0, \log_2 n], i \in[1,s]}}  \frac{  E_{s} \big(  \psi_{ \bfr_j +\balpha_{\varrho_1(j)}}(\mu_jP_{\balpha_{\varrho_1(j)}},\bx)
 \psi_{ \bfr_j +\balpha_{\varrho_2(j)}}(-\mu_jP_{\balpha_{\varrho_2(j)}},\bx) \big)}{ P_{\balpha_{\varrho_1(j)}}P_{\balpha_{\varrho_2(j)}}}
 \\
   = \sum_{ \ddot{\alpha}_i,\dddot{\alpha}_i  \in [0,   \log_2 n], i \in[1,s] }
 P^{-1}_{\ddot{\balpha}} P^{-1}_{\dddot{\balpha}}
  E_{s} \Big(
 \psi_{ \bfr_j +\ddot{\balpha}}(\mu_jP_{ \ddot{\balpha}},\bx)
 \psi_{ \bfr_j +\dddot{\balpha}}(-\mu_jP_{ \dddot{\balpha}},\bx) \Big)=: \hat{G}_{\bfr_j,\bmu_j}.
\end{multline*}
Therefore, $G_{\varrho,\bfr_j,\bmu_j}$ does not depend on $(\varrho_1,\varrho_2)$.
 From  \eqref{Th2-33a} and \eqref{Th2-34abc}, we obtain
\begin{multline*}
\hat{\varsigma}_{2q} =\frac{(2q)!}{q! 2^{q}}   \sum_{\substack{\fr_{i,j} \in [V_1/2,n], P_{\bfr_j} \leq 2^{n+\log_2^3 n}\\ i \in[1,s], j \in [1,q]}} \;\;
\sum_{\substack{\mu_j \in I^{\prime}_{n}  \\ j \in [1,q]}}
  \prod_{j=1}^{q}  \frac{2}{|P_{\bfr_j}(1-e(\mu_j/P_{\bfr_j}))|^2} \hat{G}_{\bfr_j,\bmu_j} \\
 = \frac{(2q)!}{ q!2^{q}} \prod_{j=1}^{q}\Bigg(   \sum_{\substack{\fr_{i,j} \in [V_1/2,n], i \in[1,s]\\
 P_{\bfr_j} \leq 2^{n+\log_2^3 n}}} \;\;
\sum_{\mu_j \in I^{\prime}_{n} }
    \frac{2 \hat{G}_{\bfr_j,\bmu_j} }{|P_{\bfr_j}(1-e(\mu_j/P_{\bfr_j}))|^2}\Bigg).
\end{multline*}
Applying \eqref{Th2-34}, \eqref{Th2-34a1} and \eqref{Th2-34b},  we get
\begin{equation*}
\varsigma_{2q} =  \frac{(2q)!}{ q!2^{q}} \prod_{j=1}^{q}\Bigg(   \sum_{\substack{\fr_{i,j} \in [V_1/2,n], i \in[1,s]\\
 P_{\bfr_j} \leq 2^{n+\log_2^3 n}}} \;\;
\sum_{\mu_j \in I^{\prime}_{n} }
    \frac{2 \hat{G}_{\bfr_j,\bmu_j} }{|P_{\bfr_j}(1-e(\mu_j/P_{\bfr_j}))|^2}\Bigg) + O(n^{qs-\frac{9}{10}}).
\end{equation*}
Using this statement for $q=1$, we obtain an expression for $\varsigma_{2} $ and that
\begin{equation*}
   \varsigma_{2q}=\frac{(2q)!}{ q!2^{q}} \big(\varsigma_{2}  +  O(n^{s-9/10}) \big)^q  \; + \; O(n^{qs-9/10}).
\end{equation*}
In view of \eqref{Le3-2}, we have $\varsigma_{2}  \ll n^s $. Hence
\begin{equation*}
   \varsigma_{2q} =\frac{(2q)!}{ q!2^{q}}\varsigma_{2}^q \; + \; O(n^{qs-9/10}) .
\end{equation*}
Therefore, Lemma 19 is proved. \qed \\ \\
\subsubsection{\textcolor{blue}{ End of the proof of Theorem 2}}

 In the following Lemma 20, we recall the connection between Halton's sequences and Hammersley's point sets. Next, applying Lemma 16, we obtain the assertion of Lemma 20. \\ \\
{\bf Lemma 20.} {\it Let $\bar{\bx}=(x_1,...,x_s,x_{s+1})$ and $\bx =(x_1,...,x_s) $. Then}
\begin{equation}
\imath_1:=  D(\bar{\bx}, \cH_{s+1,N} ) =   D(\bx, (H_s(k))_{k=0}^{[Nx_{s+1}]-1}  ) + \epsilon_1,\label{Th2-23}
\end{equation}
\begin{equation}
 D(\bar{\bx}, \cH^{sym}_{s+1,N} ) =   D(\bx, (H_s(k))_{k=-[Nx_{s+1}]}^{[Nx_{s+1}]-1}  ) + 4\epsilon_2, \quad |\epsilon_i| \leq 1 ,\; i=1,2. \label{Th2-23a}
\end{equation}

\begin{equation} \label{Th2-1-1}
E_{s+1}(\sD_{T_{s},\sh,1}(0,[Nx_{s+1}]  )) \ll n^{ \sh  s /2 -1/5 }, \;\;\;  s \geq 3,
\end{equation}

\begin{equation}  \label{Th2-1-2}
  \left\| \fD_{T_{\fs}}(0,[Nx_{s+1}]  ) \right\|_{s+1,2q}  \ll n^{  s/2  -1/(10q) }, \; s \geq 3, s >\fs,
\end{equation}
\begin{equation}   \label{Th2-1-3}
 E_{s+1}(\sD_{T_{2},\sh,1}(-[Nx_{s+1}],2[Nx_{s+1}] ))  \ll n^{ \sh  s /2 -1/5 }, \;\;\; s=\fs=2
\end{equation}
and
\begin{equation}  \label{Th2-1-4}
  \left\|\fD_{T_{\fs}}(-[Nx_{s+1}],2[Nx_{s+1}] )  \right\|_{s+1,2q}  \ll n^{  s/2  -1/(10q) },
 \; s=2,\fs=1 .
\end{equation} \\
{\bf Proof.} %
We will prove \eqref{Th2-23}. The proof of  \eqref{Th2-23a} is similar.
From \eqref{In1} and \eqref{In6a},     we have
\begin{equation*}
\imath_1= \card \{ 0 \leq k <N \;| \; \phi_i(k) <x_i, \; i=1,..,s,\; k/N < x_{s+1}\}  -  x_1 \cdots x_s x_{s+1} N.
\end{equation*}
Hence
\begin{multline*}
 \card \{ 0 \leq k <[Nx_{s+1}] \;| \; \phi_i(k) <x_i, \; i=1,..,s\}  -  x_1 \cdots x_s[Nx_{s+1}] -1
   \leq
\imath_1 \\ \leq \card \{ 0 \leq k <[Nx_{s+1}] \;| \; \phi_i(k) <x_i, \; i=1,..,s\}  -  x_1 \cdots x_s[Nx_{s+1}] +1
\end{multline*}
and \eqref{Th2-23} follows.\\

Let's consider  \eqref{Th2-1-1}.
 In view of \eqref{In3}, we have
\begin{equation} \label{In33}
     \left\| f(\bx) \right\|_{s,p}= \Big( E_s( |f(\bx)|^p)\Big)^{1/p}, \qquad  E_s(f(\bx) )=  \int_{[0,1)^s}  f(\bx) d \bx .
\end{equation}
Applying Lemma 16, we get
\begin{equation}    \nonumber
 |E_{s+1}(\sD_{T_{s},\sh,1}(0,[Nx_{s+1}]  ))|  \leq \max_{x_{s+1}}  |E_{s}(\sD_{T_{s},\sh,1}(0,[Nx_{s+1}]  ))| \ll n^{ \sh  s /2 -1/5 }.
\end{equation}
Hence \eqref{Th2-1-1} is proved.\\

 Let's consider \eqref{Th2-1-2}.
From \eqref{Lem2-6} and \eqref{Lem2-7}, we obtain
\begin{equation}    \nonumber
\fD_{T_{\fs}}^{2q}( 0,[Nx_{s+1}] ) =\sD_{T_{\fs},2q,1}(0,[Nx_{s+1}]) \quad \for \quad s >\fs .
\end{equation}
From \eqref{Th2-1-1}, we have
\begin{equation}    \nonumber
E_{s+1}(\fD_{T_{\fs}}^{2q}( 0,[Nx_{s+1}] )) =E_{s+1}(\sD_{T_{\fs},2q,1}(0,[Nx_{s+1}])) \ll n^{ q  s  -1/5 }.
\end{equation}
Therefore \eqref{Th2-1-2} is proved.

 Let's consider the case $s=2$.  Using Lemma 16, we get \eqref{Th2-1-3}.\\
According to \eqref{Lem2-6} and \eqref{Lem2-7}, we obtain
\begin{equation*}
  \fD_{T_{\fs}}^{2q}(-[Nx_{s+1}],2[Nx_{s+1}] )
 =\sD_{T_{\fs},2q,1}(-[Nx_{s+1}],2[Nx_{s+1}] )\;\; \for \;\;2= s >\fs=1.
\end{equation*}
 Now  \eqref{Th2-1-4} follows from Lemma 18 and \eqref{In33}.

Hence, Lemma 20 is proved. \qed  \\

 In the following Lemma 21, we show that the case $\varpi_{\br,\bm,2} =1$
 (see definitions in \eqref{Lem2-6}, \eqref{Lem2-7} and \eqref{Beg-28}) is the main un the proof of Theorem 2. In other words, Theorem 2 follows directly from Lemma 19.\\ \\
{\bf Lemma 21.} {\it With notations as above}
\begin{multline}\label{Th2-24a}
\imath_2:=E_{s+1}( D^{\sh}(\bar{\bx}, \cH_{s+1,N} ) )  = E_{s+1}(\sD_{T_{s},\sh,2}(0,[Nx_{s+1}]  ))
 + O( n^{ \sh  s /2 -1/10 }),\; \; s\geq 3,\\
E_{s+1}( D^{\sh}(\bar{\bx}, \cH_{s+1,N}^{sym} ) )
 = E_{s+1}(\sD_{T_{s},\sh,2}(-[Nx_{s+1}],2[Nx_{s+1}]-1  ))\\
 + O( n^{ \sh  s /2 -1/10})\quad  \for \quad s=2.
\end{multline} \\ \\
{\bf Proof.} %
Let's consider  the first estimate.
Let
\begin{equation} \label{Th2-prof}
  \dddot{\rD}([Nx_{s+1}]):= D(\bar{\bx}, \cH_{s+1,N} )  -\fD_{T_{s}}(0,[Nx_{s+1}]).
\end{equation}
By \eqref{Beg-28}, we have
\begin{equation}  \label{Th2-profa}
   \tilde{\rD}(0,[Nx_{s+1}])= D(\bx, (H(k))_{k=0}^{[Nx_{s+1}]-1}  )  -\fD_{T_{s}}(0,[Nx_{s+1}])
\end{equation}
and
\begin{equation} \label{Th2-profb}
|  \tilde{\rD}(0,[Nx_{s+1}])   |  \leq
\sum_{\fs=1}^{s-1}  \sum_{T_{\fs} \subseteq\{1,...,s\} }
  | \fD_{T_{\fs}}(0,[Nx_{s+1}]) | +O( \log^{3s} n) .
\end{equation}
In view of \eqref{Th2-23}, \eqref{Th2-profa} and \eqref{Th2-profb}, we have
\begin{equation*}
         |\dddot{\rD}([Nx_{s+1}])| \leq  |\tilde{\rD}([Nx_{s+1}])| +1.
\end{equation*}
Applying  Minkowski's inequality and  Lemma 20, we obtain for $s \geq 3$:
\begin{equation}  \nonumber
\left\| \dddot{\rD}([Nx_{s+1}])    \right\|_{s+1,2\nu}  \leq
\sum_{\fs=1}^{s-1}  \sum_{T_{\fs} \subseteq\{1,...,s\} }
   \left\| \fD_{T_{\fs}}(0,[Nx_{s+1}])   \right\|_{s+1,2\nu}  +O( \log^s n) .
\end{equation}
Taking into account that $s \geq 3$,\; $s > \fs $ is true in this inequality, from Lemma 20 we obtain that
\begin{equation}  \label{Th2-24-1}
\left\| \dddot{\rD}([Nx_{s+1}]) \right\|_{s+1,2\nu}    \ll n^{\frac{s}{2}- \frac{1}{10\nu} } \quad \nu=1,2, \ldots \;.
\end{equation}
From \eqref{Le13-2}, we derive
\begin{equation}  \label{Th2-24-2}
E_{s+1}( \fD_{T_{s}}^{2\nu}(0,[Nx_{s+1}]  )) \leq \sup_{0 \leq x_{s+1} \leq 1} E_{s}( \fD_{T_{s}}^{2\nu}(0,[Nx_{s+1}]  )) \ll    n^{ \nu s   }, \; \nu=1,2,...\; .
\end{equation}
By  \eqref{Th2-prof}, we have
\begin{multline*}
\imath_2 = E_{s+1}(( \fD_{T_{\fs}}(0,[Nx_{s+1}]  ) + \dddot{\rD}([Nx_{s+1}]) )^{\sh}) = E_{s+1}( \fD_{T_{\fs}}^{\sh}(0,[Nx_{s+1}]  ) ) \\
 + \epsilon_1
 2^{\sh}\sum_{1 \leq \nu   \leq \sh} E_{s+1} \Big( \big|\fD_{T_{s}}(0,[Nx_{s+1}]  )\big|^{\sh-\nu}  \big| \dddot{\rD}([Nx_{s+1}]  )\big|^{\nu} \Big).
\end{multline*}
Using \eqref{Lem2-7}, \eqref{Th2-24-1}, \eqref{Th2-24-2},  Cauchy-Shwarz's inequality and Lemma 20, we get
\begin{multline*}
\imath_2  =  E_{s+1}( \fD_{T_{\fs}}^{\sh}(0,[Nx_{s+1}]  ) )  \\
 +  \epsilon_2
 2^{\sh}\sum_{1 \leq \nu   \leq \sh}
 \Big( E_{s+1} ( |\fD_{T_{s}}(0,[Nx_{s+1}]  )|^{2\sh-2\nu}  )
    E_{s+1}(|\dddot{\rD}^{2\nu}([Nx_{s+1}]  )|^{2\nu}) \Big)^{1/2}  \\
=  E_{s+1}\big( \fD_{T_{\fs}}^{\sh}(0,[Nx_{s+1}]  ) \big)
  + O\Big(  \sum_{1 \leq \nu   \leq \sh} \big( E_{s+1} ( |\fD_{T_{s}}(0,[Nx_{s+1}]  )|^{2\sh-2\nu}    n^{\nu s -1/5 } ) \big)^{1/2}  \Big)
\end{multline*}
\begin{multline*}
  =   E_{s+1}\big( \fD_{T_{\fs}}^{\sh}(0,[Nx_{s+1}]  ) \big)
 +  O\big( \sum_{1 \leq \nu   \leq \sh}  n^{ \sh s /2 -1/10 }\big) \\
 = E_{s+1}\big(\sD_{T_{s},\sh,2}(0,[Nx_{s+1}]  )\big)  + O( n^{ \sh  s /2 -1/10 }), \; |\epsilon_i| \leq 1, \; i=1,2.
\end{multline*}
Hence, the first estimate is proved. The proof of the second estimate is similar. We need only use
\eqref{Th2-1-3} and  \eqref{Th2-1-4} instead of \eqref{Th2-1-1} and  \eqref{Th2-1-2}.
Therefore, Lemma 21 is proved. \qed \\ \\
{\textcolor{blue}{\bf Proof of Lemma A and completion of the proof of Theorem 2.}} \\

Let's consider Lemma A for odd $\sh=2q+1$. Let $s \geq 3$.
 From \eqref{Lem2-6},  we get $ E_{s+1}( \sD_{T_{s},\sh,2}(0,[Nx_{s+1}]  ))=0$. According to Lemma 21, we have
\begin{equation*}
   E_{s+1}( D^{\sh}(\bar{\bx}, \cH_{s+1,N} ) )   = O( n^{ \sh  s /2 -1/10 }).
\end{equation*}
By \eqref{In2} and Roth's inequality  \eqref{In4}, we obtain
\begin{equation}  \label{Th2-R}
 \underline{ \lim}_{n \to \infty} n^{-s/2} \left\|   D( \bar{\bx}, \cH_{s+1,N} )  \right\|_{s+1,2}  >0 .
\end{equation}
Therefore,  \eqref{Beg-0} is proved for $s \geq 3$ and odd $\sh$.
Similarly, for $s=2$ from \eqref{Th2-24a} we obtain that $E_{s+1}( D^{\sh}(\bar{\bx}, \cH_{s+1,N}^{sym} ) ) O( n^{ \sh  s /2 -1/10 }) $ and assertion \eqref{Beg-0} follows. \\

 Let's consider Lemma A with even $\sh=2q $. Let $s \geq 3$.
 In view of Theorem 1, we have $ E_{s+1} (D^{2}(\bar{\bx}, \cH_{s+1,N} )   \ll n^s $. Now using Lemma 19 and Lemma 21, we get
\begin{multline*}
E_{s+1} (D^{2q}(\bar{\bx}, \cH_{s+1,N} )  )  = E_{s+1}(\sD_{T_{s},2q,2}(0,[Nx_{s+1}]  ))
 + O( n^{  q s  -1/10 })   \\
 = \frac{(2q)!}{ 2^{q} q!}\Big(E_{s+1} ( \sD_{T_{s,2,2}}( 0,[Nx_{s+1}])  ) \Big)^q
 + O( n^{  q s  -1/10 }) \\
 = \frac{(2q)!}{ 2^{q} q!}\Big(  E_{s+1} (D^{2}(\bar{\bx}, \cH_{s+1,N} ) + O( n^{   s  -1/10 })  )    \Big)^q
 + O( n^{  q s  -1/10 }) \\
 = \frac{(2q)!}{ 2^{q} q!} (E_{s+1} (D^{2}(\bar{\bx}, \cH_{s+1,N} )  ))^q
 +  O(n^{qs-1/10}).
\end{multline*}
Applying \eqref{Th2-R}, we obtain  \eqref{Beg-0} for $s \geq 3$ and even $\sh $.
The proof for the case $s=2$ follows from   \eqref{Th2-31} and \eqref{Th2-24a}.
Hence, Lemma A and Theorem 2 are proved.  \qed  \\ \\
\subsection{\textcolor{blue}{Proof of Theorem 3}}

  We need the following simple variant of the
 {\it Continuous Mapping Theorem} (see \cite[Theorem 3.2.4., p.101]{Du}).\\ \\
{\bf Theorem B.} {\it Let $g$ be a continuous function. If $X_N \stackrel{w}{\rightarrow}   X$, then } $g(X_N) \stackrel{w}{\rightarrow}   g(X)$.

By \cite[p.31]{Bil}, a simple condition of uniform integrability of a sequence of functions
 $X_n$ is that $\sup_n E|X_n|^{1+\epsilon} < \infty$. According to \cite[Theorem 3.5, p.31]{Bil}, we have\\ \\
{\bf Theorem C.}
 {\it If $X_N$ are uniformly integrable and $X_N \stackrel{w}{\rightarrow}   X$,
  then $X$ is integrable and $E(X_N) \to E(X)$ .   } \\ \\
 We will consider the case $s \geq 3$.  The proof for the case $s \geq 2$ is similar.
Let  $ Y_N:= D(\bx, (\cH_{s+1,N} ))/D_{s+1,2}(\cH_{s+1,N}) $.
By Theorem 2,  $Y_N \stackrel{w}{\rightarrow}   \cN(0,1)=:Y$. We take the continuous function $g(x)=|x|^p$.
Using Theorem B, we get $g(Y_N) \stackrel{w}{\rightarrow}   g(Y)$.
   Bearing in mind Theorem 1, we get that the functions $g(Y_N)\; (N=1,2,...)$ are uniformly integrable.

   Now using Theorem C, we get the assertion of Theorem 3. \qed \\ \\
{\bf Remark 3.}
Let $b \geq 2$ be integer, and let $y=.y_{1} y_{2},\dots y_{j} \dots =\sum\nolimits_{j=1}^{\infty} y_{j} / b^j$  be the
$b$-expansion of numbers $y$.
We define von Neumann-Kakutani's $b-$adic
adding machine :
\begin{equation} \label{R3-1}
    \TT_{b}(y) := (y_k+1)/{b}^k + \sum\nolimits_{j \geq k+1}  y_{j} / q^{j}, \qquad \;\;
		\TT_{b}^n(y)=\TT_{b}(\TT_{b}^{n-1}(y)),
\end{equation}
$n=2,3,\dots,\TT_{b}^0(y)=y,$ where $k= \min \{j \;| \; y_j \ne b-1 \}$.

Let $\cP=(p_1,...,p_{s})$ and let  $\TT_{\cP}(\by) =(\TT_{p_1}(y_1), ..., \TT_{p_s}(y_s))$.
As is known, the  sequence $(\TT_{p_i}^k(y))_{k \geq 0}$ coincides for $y=0$  with the van
der Corput sequence in base ${p_i}$ $(1 \leq i \leq s)$ (see e.g., \cite[\S 2.5]{FKP}). Hence $T^k_{\cP}(0) =H_s(k)$, $k=0,1,...$ .
It is easy to show that  $(\TT_{\cP}^k(\by))_{k \geq 0}$ is the low discrepancy ergodic transformation for all $\by$.

In view of \eqref{Beg-2-0} and  \eqref{R3-1}, we obtain
\begin{equation}  \nonumber
   | \TT_{\cP}(\by)  - \TT_{\cP}([\by]_{2\bn}) | \leq s2^{-2n}, \quad \with \quad \bn=(n,n,...,n), n=[\log_2 N] +1.
\end{equation}
By \eqref{In1} and \eqref{In6}, we  have
\begin{equation}  \nonumber
  D(\bx, ( \TT_{\cP}^k(\by) )_{k=0}^{N-1}  ) =  D(\bx, ( \TT_{\cP}([\by]_{2\bn}) )_{k=0}^{N-1}  ) + \epsilon s,  \;\;\; |\epsilon| \leq 1.
\end{equation}
We define  $Y_{2 \bn}$ from the condition $ [\by]_{2\bn} = H_s(Y_{2 \bn}) = T^{Y_{2 \bn} }_{\cP}(0) $.

Applying \eqref{In6} and \eqref{R3-1},  we get
\begin{equation*}  \nonumber
   T^k_{\cP}( [\by]_{2\bn} ) =  T^k_{\cP}( T^{Y_{2 \bn} }_{\cP}(0)  ) = T^{k+Y_{2 \bn} }_{\cP}(0)   =    H_s(k +Y_{2 \bn}).
\end{equation*}
Using  \eqref{Beg-18},  we obtain
\begin{multline}  \label{Rem}
  D(\bx, ( \TT_{\cP}^k(\by) )_{k=0}^{N-1}  ) =   D(\bx, ( \TT_{\cP}([\by]_{2\bn}) )_{k=0}^{N-1}  ) + \epsilon s =   D(\bx, ( H(k) )_{k=Q}^{Q+N-1}  ) + \epsilon s,\\
 \with  \quad Q = Y_{2 \bn} .
\end{multline}
Bearing in mind that the $O$ constant in Theorem 1 is independent of $Q$, we have \\ \\
{\bf Theorem 1$^{'}$.} {\it Let $s \geq 2, \; \by \in [0,1)^s\; p \geq 1$. Then}
\begin{equation} \nonumber
  D_{s,p} ((  \TT_{\cP}^k(\by) )_{k=0}^{N-1})  =O( \log^{s /2} N),
\end{equation}
where the $O$ constant is independent of $\by$. \\

We put
 ${\cT}_{\cP,s+1,N}=
	 (\TT_{\cP}^k(\by), k/N)_{k=0}^{N-1}$. For $s=2$, we put \\
${\cT}^{sym}_{\cP,3,N}=
	 ( \{{\rm sign}(k)\TT_{p_1}^{|k|}(\by) \}, \{{\rm sign}(k)\TT_{p_2}^{|k|}(\by) \}, \{ k/N\})_{k=-N}^{ N-1}$, where  ${\rm sign}(z) =1$ if $z \geq 0$ and  ${\rm sign}(z) =-1$ if $z < 0$, with $\cP=(p_1,p_2)$.

In view of \eqref{In6a}, Theorem 2 and  Theorem 3, we get \\  \\ \\
{\bf Theorem 2$^{\bf{'}}$.} {\it Let $s \geq 2$, $ \by \in [0,1)^s$, $\bar{\bx}$ be a uniformly distributed random variable in $[0,1]^{s+1}$.  Then }
\begin{equation}  \nonumber 
 \frac{ D(\bar{\bx}, {\cT}_{\cP,s+1,N} ) }{ \left\| D(\bar{\bx}, {\cT}_{\cP,s+1,N} ) \right\|_{s+1,2}}
	\stackrel{w}{\rightarrow} \cN(0,1) \; \; \for \; s \geq 3,\quad \;
  \frac{ D(\bar{\bx},{\cT}^{sym}_{\cP,3,N} ) }{\left\| D(\bar{\bx},{\cT}^{sym}_{\cP,3,N}  ) \right\|_{3,2}}
	\stackrel{w}{\rightarrow} \cN(0,1).
\end{equation}\\

Hence $ D(\bar{\bx}, {\cT}_{\cP,s+1,N} )$  satisfy to $s+1$-parametric {\it temporal central
limit theorem} (see definition in \cite{DS1}).\\
\\
{\bf Theorem 3$^{'}$.} {\it Let $s \geq 2$, $ \by \in [0,1)^s$  and $p>0$. Then}
\begin{equation} \nonumber
  \frac{ D_{s+1,p}( {\cT}_{\cP,s+1,N} )}{ D_{s+1,2}( {\cT}_{\cP,s+1,N} )}
	\stackrel{N \rightarrow \infty}{\longrightarrow}   \kappa_p^{1/p}, \; s \geq 3,
			\qquad \quad \frac{ D_{3,p}( {\cT}^{sym}_{\cP,3,N} )}{D_{3,2}( {\cT}^{sym}_{\cP,3,N} )}
	\stackrel{N \rightarrow \infty}{\longrightarrow}   \kappa_p^{1/p}.
\end{equation}

\addcontentsline{toc}{section}{\textcolor{blue}{References}}
\bibliographystyle{plain}

%

%
%
{\bf Address}: Department of Mathematics, Bar-Ilan University, Ramat-Gan, 5290002, Israel \\
{\bf E-mail}: mlevin@math.biu.ac.il\\
 https://orcid.org/0000-0003-1268-7172
\end{document}